\documentclass[leqno]{article} 
\usepackage{amsmath} 
\usepackage{amssymb} 
\usepackage[all]{xy} 
\usepackage{amsthm} 
\usepackage{authblk}
\newtheorem{thm}{Theorem}[section] 
 
\newtheorem{lem}[thm]{Lemma}

\theoremstyle{remark} 
\newtheorem{rem}{\bf {\slshape Remark}}[section] 
\numberwithin{equation}{section} 
\newtheorem{eg}[thm]{\bf {\slshape Example}}
\newtheorem{df}{\bf  Definition}

\title{Degeneration scheme of \\ 4-dimensional  Painlev\'e-type equations} 

\author{ 
Hiroshi Kawakami\thanks{\texttt{kawakami@gem.aoyama.ac.jp}}} 
\author{Akane Nakamura\thanks{\texttt{a-naka@josai.ac.jp}}}
\author{Hidetaka Sakai\thanks{\texttt{sakai@ms.u-tokyo.ac.jp}}}
\affil{${}^{*}$College of Science and Engineering, Aoyama gakuin university,
	}
\affil{${}^{\dagger}$ Faculty of Science, Josai University}
\affil{${}^{\ddagger}$Graduate School of Mathematical Sciences,  
	University of Tokyo.
	}
\date{} 

\begin{document} 

\maketitle

\begin{abstract} 
This is a continuation of the recent work \cite{S3} by one of the authors.
According to \cite{S3}, there are four 4-dimensional Painlev\'e-type equations derived from isomonodromic deformation of Fuchsian equations: they are the Garnier system in two variables, the Fuji-Suzuki system, the Sasano system, and the sixth matrix Painlev\'e system.
In this paper we degenerate these four source equations, and systematically obtain other 4-dimensional Painlev\'e-type equations.
If we only consider Painlev\'e-type equations whose associated linear equations are of unramified type, there are 22 types of 4-dimensional Painlev\'e-type equations: 9 of them are partial differential equations, 13 of them are ordinary differential equations.
Some well-known equations such as Noumi-Yamada systems are included in this list.
They are written as Hamiltonian systems, and their Hamiltonians are simply written by using Hamiltonians of the classical Painlev\'e equations.



\end{abstract} 
\section{Introduction}\label{intro}
\subsection{Background}
The Painlev\'e equations were discovered in the early twentieth century by Painlev\'e
through his study of classification of nonlinear ordinary differential equations with the Painlev\'e property.
Solutions of a nonlinear differential equation may have singular points that are 
not determined only by fixed singular points of the equation.
In this case, positions of singular points depend on initial values.  
We call such singular points movable singularities. 
As in the case of differential equations satisfied by elliptic functions,
when movable singularities of an equation are at most poles, we say that the equation
 has Painlev\'e property.
 If we restrict our attention to second order algebraic differential equations in the normal form and exclude the cases when they are solved by quadrature or solved by solutions of linear equations or elliptic functions, there are no other equations than Painlev\'e equations having the Painlev\'e property \cite{Pn2,Gm}. 
Painlev\'e transcendents
are studied from various viewpoints, as special functions defined by nonlinear differential equations.

Usually, the Painlev\'e equations are classified into six types.
However, if we consider rational surfaces called spaces of initial values \cite{O1},
the third Painlev\'e equation falls into three types, distinguished by the values of parameters.
Thus, it is natural to consider that there are 8 types of Painlev\'e equations \cite{S1}.
These surfaces are similar to rational elliptic surfaces but are slightly different.
It corresponds to the fact that autonomous limits of the Painlev\'e equations lead to the differential equations satisfied by elliptic functions.

Furthermore  the Painlev\'e equations have another origin: the isomonodromic deformation.
Here we briefly explain the idea of isomonodromic deformation.
Suppose that a Fuchsian equation is given on $X=\mathbb{P}^1\setminus\{u_1,\ldots,u_n\}$,
where $u_i$'s are singular points.
Considering an analytic continuation of its fundamental solutions, we obtain
a monodromy representation of the fundamental group of $X$.
We can define a correspondence from a Fuchsian equation to the  representations of the
fundamental group of $X$.
Riemann considered the problem whether a Fuchsian equation exists for any representation of the fundamental group of $X$. 
This problem is called the Riemann problem.
Many mathematicians dealt with the problem and now we can say that
the Riemann problem is positively solved under suitable conditions. 

On the other hand the correspondence  is not one-to-one.
More precisely, we can deform a Fuchsian equation without changing its monodromy.
Therefore it is important to describe the family of Fuchsian equations which have the same monodromy group.
This problem is called the Fuchs problem.
R. Fuchs treated a single second order Fuchsian equation whose solution has four singular points, normalized to $0,1,t$, and $\infty$.
He considered $t$ as a deformation parameter and obtained the condition that
the monodromy of the Fuchsian equation is independent of $t$.
The condition is equivalent to a nonlinear differential equation,
which is exactly the sixth Painlev\'e equation.
This is the first example of isomonodromic deformation \cite{F}.
Since then, some researchers generalized the  work of R. Fuchs.
Garnier considered the isomonodromic deformation of a second order Fuchsian equation
with arbitrary number of singular points.
This yields nonlinear partial differential equations, which are called Garnier systems.
Moreover Schlesinger consider the isomonodromic deformation of a system of linear differential equations of
the following form:
\[
\frac{dY}{dx}=\left(\sum_{i=1}^{n}\frac{A_i}{x-u_i}\right)Y.
\]
Here $A_i$'s are matrices independent of $x$.
The corresponding nonlinear equation is called the Schlesinger system.
Thus the Schlesinger system include the sixth Painlev\'e equation and the Garnier system as special cases.
A general theory of isomonodromic deformation that can treat non-Fuchsian equations
is established by Jimbo, Miwa, and Ueno \cite{JMU}.

From now on, we use the term ``Painlev\'e-type equations" for nonlinear equations derived by generalized isomonodromic deformations of linear differential equations.
In the present paper, linear differential equations are written in
the form of system, not in single equations.
When the order of Painlev\'e-type equation is two, we know that all the equations with Painlev\'e property have the corresponding isomonodromy problem of linear equations.
However, we still do not know whether all the equations with Painlev\'e property can be derived from isomonodromic deformation.

It is known that in general, Painlev\'e-type equations are written in ``multi-time" Hamiltonian systems:
\[
\frac{\partial q_j}{\partial t_{i}}=\frac{\partial H_{t_{i}}}{\partial p_j}, \quad
\frac{\partial p_j}{\partial t_{i}}=-\frac{\partial H_{t_{i}}}{\partial q_j}.
\]
There are as many Hamiltonians as independent variables.

Let us see the explicit forms of Hamiltonians for the Painlev\'e equations \cite{O2}:
\begin{align*} 
& t(t-1)H_{\rm VI}\left({\alpha , \beta\atop\gamma, \delta 
 };t;q,p\right)=\;q(q-1)(q-t)p^2\\ 
 & \hspace{13em}+\{ \delta q(q-1)-(2\alpha +\beta +\gamma +\delta )q(q-t)+\gamma
 (q-1)(q-t)\} p\\ 
 & \hspace{13em}+\alpha (\alpha +\beta )(q-t),\\ 
& tH_{\rm V}\left({\alpha , \beta \atop \gamma };t;q,p\right)=\;p(p+t)q(q-1) 
+\beta pq+\gamma p-(\alpha +\gamma )tq,\\ 
& H_{\rm IV}\left(\alpha , \beta;t;q,p\right)=\; 
pq(p-q-t)+\beta p+\alpha q,\quad
 tH_{\rm III}(D_6)\left(\alpha , \beta ;t;q,p\right)=\; 
p^2q^2-(q^2-\beta q-t)p-\alpha q,\\ 
& tH_{\rm III}(D_7)\left(\alpha;t;q,p\right)=\; 
p^2q^2+\alpha qp+tp+q,\quad
 tH_{\rm III}(D_8)\left(t;q,p\right)=\; 
p^2q^2+qp-q-\frac{t}{q},
\\ 
& H_{\rm II}\left(\alpha;t;q,p\right)=\; 
p^2-(q^2+t)p-\alpha q,\quad \hspace{3em}
 H_{\rm I}\left(t;q,p\right)=\; 
p^2-q^3-tq. 
\end{align*}
These are non-autonomous Hamiltonian system with $t$ as independent variable, $p$ and $q$ as canonical variables.
Corresponding autonomous systems of them can be solved by using elliptic functions.
All these 8 equations are ordinary differential equations
and expressed by one Hamiltonian.
\begin{rem}\label{rem:P_V}
As for Painlev\'e equation of the fifth type, it is sometimes convenient to use another Hamiltonian for calculations: 
\begin{equation}
t\tilde{H}_{\rm V}\left({\alpha ,\beta \atop \gamma};t;q,p\right)= 
q(q-1)^2p^2+\{(1-q)(\alpha+(\beta+2\gamma)q)+tq\}p-\gamma(\beta+\gamma)(1-q).
\end{equation}
The canonical transformation
\[
q \to 1-\frac{1}{q},\quad p \to q(p q-\gamma)
\]
changes the above $\tilde{H}_{\rm V}$ into a biquadratic polynomial
\[
H_{\rm V}\left({\beta+\gamma,\alpha+\beta \atop -\beta};t;q,p \right)-\frac{\alpha \gamma}{t}+\gamma.
\]
\qed
\end{rem}

%

Non-Fuchsian equations are derived from Fuchsian equations by degeneration such as confluence of singular points.
These limit procedures induce degenerations of the Painlev\'e equations \cite{G,JM}.

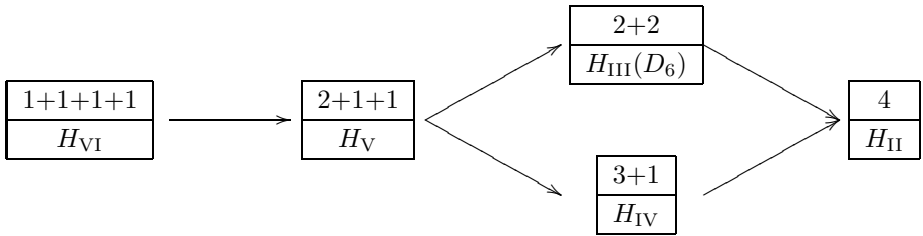
\bigskip
\begin{xy}
{(3,0) *{\begin{tabular}{|c|}
\hline
1+1+1+1\\
\hline
$H_{\rm VI}$\\
\hline
\end{tabular}
}},
{(40,0) *{\begin{tabular}{|c|}
\hline
2+1+1\\
\hline
$H_{\rm V}$\\
\hline
\end{tabular}
}},
{\ar (15,0);(31,0)},
{\ar (49,0);(67,10)},
{\ar (49,0);(67,-10)},
{(77,10) *{\begin{tabular}{|c|}
\hline
2+2\\
\hline
$H_{\rm III}(D_6)$\\
\hline
\end{tabular}}},
{\ar (86,10);(104,0)},
{\ar (86,-10);(104,0)},
{(77,-10) *{\begin{tabular}{|c|}
\hline
3+1\\
\hline
$H_{\rm IV}$\\
\hline
\end{tabular}}},
{(110,0) *{\begin{tabular}{|c|}
\hline
4\\
\hline
$H_{\rm II}$\\
\hline
\end{tabular}}},\end{xy}
\medskip
\noindent

\medskip

Here the symbol in the upper half of each box stands for the Poincar\'e rank of singular points of
corresponding linear equation.
For example, ``2+1+1" indicates that
 there are one irregular singular point of the Poincar\'e rank 1 and two regular singular points of Poincar\'e rank 0
 (that is, regular singularity).
  See subsection 1.3.

In the scheme above, we only considered confluence processes.
If we also consider degenerations of Jordan canonical forms of principal parts of coefficient matrices
of linear equations, the degeneration scheme becomes as follows:

\bigskip
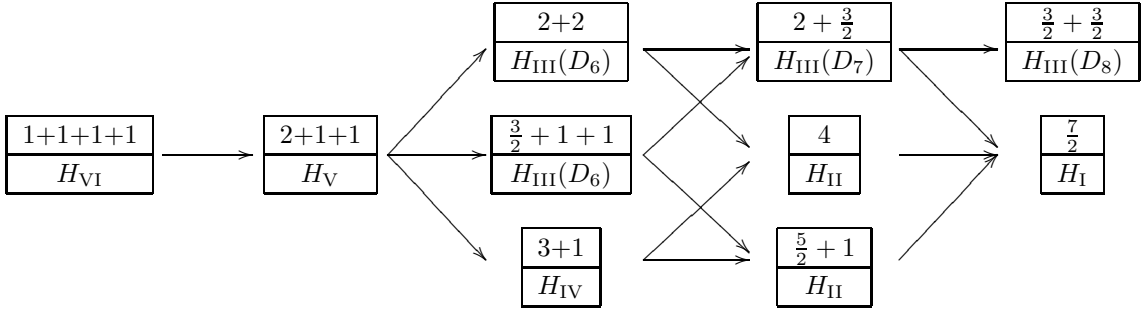
\begin{xy}
{(3,0) *{\begin{tabular}{|c|}
\hline
1+1+1+1\\
\hline
$H_{\rm VI}$\\
\hline
\end{tabular}
}},
{(35,0) *{\begin{tabular}{|c|}
\hline
2+1+1\\
\hline
$H_{\rm V}$\\
\hline
\end{tabular}
}},
{\ar (14,0);(26,0)},
{\ar (44,0);(57,14)},
{\ar (44,0);(57,0)},
{\ar (44,0);(57,-14)},
{(67,15) *{\begin{tabular}{|c|}
\hline
2+2\\
\hline
$H_{\rm III}(D_6)$\\
\hline
\end{tabular}}},
{(67,0) *{\begin{tabular}{|c|}
\hline
$\frac{3}{2}+1+1$\\
\hline
$H_{\rm III}(D_6)$\\
\hline
\end{tabular}}},
{\ar (78,14);(92,14)},
{\ar (78,14);(92,1)},
{\ar (78,0);(92,13)},
{\ar (78,0);(92,-13)},
{\ar (78,-14);(92,-1)},
{\ar (78,-14);(92,-14)},
{(67,-15) *{\begin{tabular}{|c|}
\hline
3+1\\
\hline
$H_{\rm IV}$\\
\hline
\end{tabular}}},
{(102,15) *{\begin{tabular}{|c|}
\hline
$2+\frac{3}{2}$\\
\hline
$H_{\rm III}(D_7)$\\
\hline
\end{tabular}}},
{(102,0) *{\begin{tabular}{|c|}
\hline
4\\
\hline
$H_{\rm II}$\\
\hline
\end{tabular}}},
{\ar (112,14);(125,1)},
{\ar (112,14);(125,14)},
{\ar (112,0);(125,0)},
{\ar (112,-14);(125,0)},
{(102,-15) *{\begin{tabular}{|c|}
\hline
$\frac{5}{2}+1$\\
\hline
$H_{\rm II}$\\
\hline
\end{tabular}}},
{(135,15) *{\begin{tabular}{|c|}
\hline
$\frac{3}{2}+\frac{3}{2}$\\
\hline
$H_{\rm III}(D_8)$\\
\hline
\end{tabular}}},
{(135,0) *{\begin{tabular}{|c|}
\hline
$\frac{7}{2}$\\
\hline
$H_{\rm I}$\\
\hline
\end{tabular}}},\end{xy}

\medskip
\noindent
This scheme was obtained by Ohyama and Okumura \cite{OO}.
Newly added parts require Puiseux series for formal solutions of corresponding linear equations.
We say that such equations are of ramified type, and distinguish from the unramified case.
We only treat unramified cases in this paper.

\subsubsection*{Higher dimensional Painlev\'e-type equations}
The theory of Painlev\'e equations is generalized to higher order nonlinear differential equations, and some important equations are proposed and investigated.
For examples, Gordoa, Joshi, and Pickering proposed a higher order generalization of the second and the fourth Painlev\'e equations \cite{GJP}.
According to Koike's work, they turned out to be restrictions of independent variables of the classically known degenerate Garnier systems \cite{Kk}.
However, some equations such as Noumi-Yamada systems or the Sasano system still remain beyond understanding by the framework of classical Painlev\'e equations or Garnier systems.

These equations seem to be sporadic, and not well-organized compared to the case of 2-dimensional phase space.
This article aims to make a theory of clear classification when the phase spaces are 4-dimensional. 
For that purpose, we approach this problem from the theory of isomonodromic deformation of linear differential equations.
That is to say, we try to make a list of all the 4-dimensional Painlev\'e-type equations by associated
linear differential equations. 

%

Let us recall what have done in the previous paper by the one of the authors \cite{S3}.
In that paper, four 4-dimensional Painlev\'e-type equations are derived by isomonodromic deformation of Fuchsian equations.
They are the Garnier system in two variables, the Fuji-Suzuki system, the Sasano system, and the sixth matrix Painlev\'e system.
The point is that these four equations are all the 4-dimensional Painlev\'e-type equations that we can obtain from isomonodromic deformation of Fuchsian equations.
This can be verified by the classification theory of Fuchsian equations, which has developed in recent years.

We treat Fuchsian equations in the form of first-order systems called Schlesinger normal forms:
\begin{equation}\label{eq:linear_system}
\frac{dY}{dx}=\left(\sum_{i=1}^n\frac{A_i}{x-u_i}\right)Y,
\end{equation}
where $A_{i}$'s are $m \times m$ matrices.
Usually, Fuchsian equations are roughly classified by the ranks and the number of singularities.
In addition to the ordinary rough classification, it is important to take  multiplicities of eigenvalues of $A_{i}$'s into account.
We classify Fuchsian equations using ``spectral type".

For the system above, the spectral type is an $(n+1)$-tuple of partitions of $m$:
\[
 m^1_1m^1_2\ldots m^1_{l_1}, m^2_1\ldots m^2_{l_2}, \ldots ,
 m^n_1\ldots m^n_{l_n}, m^{\infty}_1\ldots m^{\infty}_{l_\infty},
 \qquad \sum_{j=1}^{l_i}m_j^i=m~~\mbox{for}~~ i=1,\ldots, n, \infty.
\]
This means that there are $m^i_j \ (1\leq j\leq l_i)$ multiple 
eigenvalues of $A_{i}$.
Here we have assumed that $A_i$'s and $A_\infty :=-\sum_{i=1}^n A_i$ are diagonalizable.

Katz introduced two transformations called middle convolution and addition \cite{Kt,DR}.
The Painlev\'e-type equations are invariant under Katz's two operations~\cite{HF}.
We define rigidity index of the system as
\begin{equation*}
{\rm idx}=-(n-1)m^{2}+\sum_{i=0}^{n}\sum_{j=1}^{l_{i}}(m_{j}^{i})^{2}.
\end{equation*}
Katz's two operations preserve rigidity index.
Furthermore, if we fix the rigidity index, Fuchsian equations are reduced to finite types of equations by the two transformations \cite{Os}.
In particular, we have
\begin{thm}[Oshima]
Any irreducible Fuchsian equation with rigidity index $-2$ is equivalent to one of the 13-types of equations listed below, via Katz's transformations:\\
{\rm
11,11,11,11,11\\
21,21,111,111\quad 31,22,22,1111\quad 22,22,22,211\\
211,1111,1111\quad 221,221,11111\quad 32,11111,11111\quad 222,222,2211\quad 33,2211,111111\quad
44,2222,22211\\
 44,332,11111111\quad 55,3331,22222\quad 66,444,2222211.}
\end{thm}
The number of accessory parameters of the system coincides with $2-{\rm idx}$.
Accessory parameters are parameters of the system that are not determined by characteristic exponents at singular points.
When the rigidity index equals to $-2$, the number of accessory parameters is 4, and the phase spaces for the corresponding Painlev\'e-type equations are 4-dimensional.
Fuchsian equations with only three singular points do not admit non-trivial isomonodromic deformation.
We have four types of Fuchsian equations whose number of singular points is equal to or more than four:
\begin{equation*}
11,11,11,11,11\quad
21,21,111,111\quad 31,22,22,1111\quad 22,22,22,211.\\
\end{equation*}
Thus, there are just four non-trivial Painlev\'e-type equations corresponding to Fuchsian equations.
They are the Garnier system in two variables, the Fuji-Suzuki system, the Sasano system, and the sixth matrix Painlev\'e system, respectively.
These deformation equations are expressed in forms of Hamiltonians.
We will explain how to derive Hamiltonians in section 4.1.

\subsection{Non-Fuchsian linear equations}
Let us move on to our main interest in this paper: the non-Fuchsian cases.
We use two key ideas in this paper.
One of them is the notion of ``spectral type" for non-Fuchsian case.
It helps us to distinguish Painlev\'e type equations.
The other key idea is expressing Hamiltonians of Painlev\'e-type equations by Hamiltonians of classical Painlev\'e equations.
This idea is effective to identify equations. 

We make use of the classification of Fuchsian equations and their Painlev\'e-type equations,
then take  degeneration to obtain other Painlev\'e-type equations.
As a result, we obtained 22 types of equations, as we shall see in the next section.
Among these 22 types, 9 of them are partial differential systems, and 13 are ordinary differential systems.
It sometimes happens that different degenerations yield the same Painlev\'e-type equation.
In these cases, corresponding linear equations transform to one another by the Laplace transform.
We shall see such interesting topics in the last section.

In this subsection, we define the spectral type for a system of linear differential equations with irregular singularities.
The spectral type will be defined through the ``canonical form'' at each singular point.
In the irregular singular case, what corresponds to each singular point is not merely a partition,
but a {\it refining sequence of partitons} (Definition~\ref{def:rsp}).

Here we recall the canonical form of a system of linear differential equations at a singular point.
Linear equations that we treat in the present paper are those with rational function coefficients.
Such a system is generally written as follows:
\begin{equation}
\frac{dY}{dx}=
\left(\sum_{\nu=1}^n\sum_{k=0}^{r_{\nu}}\frac{A_{\nu}^{(k)}}{(x-u_{\nu})^{k+1}}
+\sum_{k=1}^{r_{\infty}}A_{\infty}^{(k)}x^{k-1}
\right)Y, \quad
A_j^{(k)} \in M_m(\mathbb{C}).
\end{equation}
This system has singular points at $x=u_1,\ldots,u_n$ and $\infty$.
Choosing one of the singular points and taking the local coordinates $z=x-u_\nu$ or $z=1/x$
 we consider the system at $z=0$:
\begin{equation}\label{eq:Laurent}
\frac{dY}{dz}=\left(
\frac{A_0}{z^{r+1}}+\frac{A_1}{z^{r}}+\cdots+A_{r+1}+A_{r+2} z+\cdots
\right)Y.
\end{equation}
We denote the field of formal Laurent series in $z$ by $\mathbb{C}((z))$,
and the field of Puiseux series $\cup_{p > 0}\mathbb{C}((z^{\frac1p}))$ by $\mathcal{K}_z$.
\begin{thm}[Hukuhara, Turrittin, Levelt]
	For any
	\begin{equation}\label{eq:Laurent2}
	\frac{dY}{dz}=\left(\frac{A_0}{z^{r+1}}+\frac{A_1}{z^{r}}+\cdots\right)Y, \quad
	A_j \in M_m(\mathbb{C}),
	\end{equation}
	there exists $P(z) \in \mathrm{GL}_m(\mathcal{K}_z)$ such that the transformation of dependent variable $Y=P(z)Z$
	brings the system into the following canonical form:
	\begin{equation}
	\frac{dZ}{dz}=
	\left(\frac{D_0}{z^{l_0}}+\frac{D_1}{z^{l_1}}+\cdots+\frac{D_{s-1}}{z^{l_{s-1}}}
	+\frac{\Theta}{z^{l_s}}\right)Z,
	\end{equation}
	where
	\begin{itemize}
		\item $l_0, \ldots, l_{s}\ (l_0 > l_1 > \cdots > l_{s-1}>l_s = 1)$ are rational numbers,
		
		\item $D_0, \ldots, D_{s-1}$ are diagonal matrices,
		
		\item $\Theta$ is a (not necessarily diagonal) Jordan matrix which commutes with all $D_j$'s.
	\end{itemize}
	Here $l_0, \ldots, l_{s}$ are uniquely determined only by the original system (\ref{eq:Laurent2}).
	
	If the equation
	\begin{equation}
	\frac{dW}{dz}=
	\left(\frac{\tilde{D}_0}{z^{l_0}}+\frac{\tilde{D}_1}{z^{l_1}}+\cdots+\frac{\tilde{D}_{s-1}}{z^{l_{s-1}}}
	+\frac{\tilde{\Theta}}{z^{l_s}}\right)W
	\end{equation}
	is another canonical form corresponding to the same system,
	there exist a constant matrix $g \in \mathrm{GL}_m(\mathbb{C})$ and a natural number $k \in \mathbb{Z}_{\ge 1}$ such that
	\begin{equation}
	\tilde{D}_j=g^{-1}D_j g, \quad \exp(2\pi i k \tilde{\Theta})=g^{-1}\exp(2\pi i k \Theta)g
	\end{equation}
	hold.
\end{thm}
We call the number $l_0-1$ the Poincar\'e rank of the singular point.
When there is a rational number $l_j$ that is not an integer, the singular point is called a ramified irregular singular point.
A linear differential equation is said to be of ramified type if it has a ramified irregular singular point.
In the present paper, as we have mentioned, we consider only linear equations of unramified type, namely,
we treat linear equations which do not have ramified irregular singular points.

At an unramified (irregular) singular point, the canonical form can be written in the form
\begin{equation}\label{eqn:loc_norm_form}
\frac{dZ}{dz}=
\left(\frac{T_0}{z^{l+1}}+\frac{T_1}{z^{l}}+\cdots+\frac{T_{l-1}}{z^2}+\frac{\Theta}{z}\right)Z,
\end{equation}
where some of $T_j$'s may be the zero matrix $O$.

In this paper, non-semisimple $\Theta$ does not appear.
Thus, by writing the diagonal entries of $T_j$'s and $\Theta$ as $t^j_k$ and $\theta_k\ (k=1, \ldots, m)$
respectively, we express the canonical form as follows:
\begin{equation}
\begin{array}{c}
x=u_i \ \text{(or $\infty$)} \\
\overbrace{\begin{array}{ccccc}
	t^0_1  & t^1_1 & \ldots & t^{l-1}_1  & \theta_1\\
	\vdots & \vdots    &        & \vdots & \vdots \\
	t^0_m  & t^1_m & \ldots & t^{l-1}_m  & \theta_m
	\end{array}}
\end{array}.
\end{equation}
The table of above expression of canonical forms at each singularity is called the Riemann scheme.

The canonical form (\ref{eqn:loc_norm_form})
can be easily solved.
In fact, the following matrix
\begin{equation}\label{eqn:exp_poly}
\exp\left(-\frac{T_{0}}{l z^l}-\dots-\frac{T_{l-1}}{z}\right)z^{\Theta}
\end{equation}
is a fundamental solution matrix of (\ref{eqn:loc_norm_form}).
The degree of the polynomial in the exponential function with respect to $z^{-1}$ is 
the Poincar\'e rank, which we have defined above. 


In the rest of this subsection, we briefly outline how the canonical form at an irregular singular point is computed,
and how the canonical form is described by the refining sequence of partitions.


We first consider the case $r=0$ in (\ref{eq:Laurent}), namely when the order of the pole is $1$.
\subsubsection*{The case of the simple pole}
We start with the following Lemma concerning so-called the Sylvester equation.
We will use this Lemma very effectively to solve differential equations formally.
\begin{lem}\label{thm:Sylvester_eq}
	Let $A \in M_m(\mathbb{C})$, $B \in M_n(\mathbb{C})$, and $C \in M_{m,n}(\mathbb{C})$.
	The matrix equation
	\begin{equation}
	AX-XB=C
	\end{equation}
	with respect to $m \times n$ matrix $X$ has a solution for any $C$
	if and only if
	\[
	\{ \text{eigenvalue of $A$} \} \cap \{ \text{eigenvalue of $B$} \}= \emptyset.
	\]
\end{lem}
\begin{proof}
	First, we prove the following:
	\begin{center}
		$AX-XB=O$ has a non-trivial solution
		$\iff$
		$\{ \text{eigenvalue of $A$} \} \cap \{ \text{eigenvalue of $B$} \} \neq \emptyset$.
		$(*)$
	\end{center}
	\noindent
	$(\Rightarrow)$ Assume that $\{ \text{eigenvalue of $A$} \} \cap \{ \text{eigenvalue of $B$} \}=\emptyset$.
	Then the characteristic polynomials $\Phi_A(x), \Phi_B(x)$ of $A, B$ are coprime.
	Thus there exist polynomials $f, g$ such that
	\[
	f(x)\Phi_A(x)+g(x)\Phi_B(x)=1.
	\]
	By the Cayley-Hamilton theorem, we have $g(A)\Phi_B(A)=I$.
	
	Let $X$ be a solution of $AX-XB=O$.
	Then we have
	\begin{align*}
	X=g(A)\Phi_B(A)X=g(A)X\Phi_B(B)=O.
	\end{align*}
	
	\noindent
	$(\Leftarrow)$ Let $\lambda$ be a common eigenvalue of $A$ and $B$,
	then there exist non-zero column vectors $v$ and $u$ such that
	\[
	Av=\lambda v, \quad
	{}^t\! u B=\lambda {}^t\! u.
	\]
	
	If we put $X=v \cdot{}^t\! u$, then $X \ne O$ and
	\[
	AX-XB=Av \cdot{}^t\! u-v \cdot{}^t\! u B=\lambda v\cdot{}^t\! u-\lambda v \cdot{}^t\! u=O.
	\]
	
	Now we prove the statement of the Lemma.
	Consider the following linear map
	\[
	\varphi : M_{m,n}(\mathbb{C}) \to M_{m,n}(\mathbb{C}), \quad
	\varphi(X)=AX-XB.
	\]
	Then the equation $AX-XB=C$ has a solution for any $C$
	
	$\iff$ $\varphi$ is surjective $\iff$ $\varphi$ is injective
	$\stackrel{(*)}{\iff}$ $\{ \text{eigenvalue of $A$} \} \cap \{ \text{eigenvalue of $B$} \}=\emptyset$.
\end{proof}

We consider the following linear ordinary differential system
\begin{align}
\frac{dY}{dx}=AY, \quad
A=\frac1x \left( A_0+A_1x+A_2x^2\cdots \right), \quad
A_j \in M_m(\mathbb{C}).
\end{align}
Here $A$ is allowed to be a divergent series.
We think of $A$ merely as a formal power series.

Changing the dependent variable so that $Y=PZ$ with the following formal power series
\begin{align}
P=I+\sum_{j=1}^\infty P_j x^j,
\end{align}
we obtain the equation for $Z$:
\begin{align}
\frac{dZ}{dx}&=\left( P^{-1}AP-P^{-1}\frac{dP}{dx} \right)Z.
\end{align}
We put
\begin{align}\label{eq:Bdef}
B:=P^{-1}AP-P^{-1}\frac{dP}{dx}.
\end{align}
Here $B$ can also be expanded as 
\begin{equation}
B=\frac1x \left( B_0+B_1x+\cdots \right).
\end{equation}
Rewriting (\ref{eq:Bdef}) as $AP-P'=PB$ and equating like powers of $x$,
we obtain the relations among $A_j, B_j, P_j$:
\begin{align}
B_0&=A_0, \\
(A_0-k I)P_k-P_k B_0&=B_k-A_k-\sum_{j=1}^{k-1}(A_{k-j}P_j-P_j B_{k-j}) \quad (k \ge 1). \label{eq:red_reg}
\end{align}
We regard the equation~(\ref{eq:red_reg}) as an equation for $P_k$.
If no two distinct eigenvalues of $A_0$ differ by an integer, then, by virtue of Lemma~\ref{thm:Sylvester_eq},
(\ref{eq:red_reg}) has a solution irrespective of the right-hand side.
In particular, we can choose $P_k \ (k=1,2,\ldots)$ so that $B_k=O \ (k=1,2,\ldots)$.
Then the result is
\begin{align}
\frac{dZ}{dx}&=\frac{A_0}{x}Z.
\end{align}

Next we consider the case when some eigenvalues of $A_0$ differ by integers.
Let $\lambda, \lambda+k \ (k \in \mathbb{Z}_{>0})$ be such a pair of eigenvalues.
We choose the gauge so that
\begin{align}
&A_0=
\begin{pmatrix}
A^0_{11} & O \\
O & A^0_{22}
\end{pmatrix}
\end{align}
where $A^0_{11}$ is a Jordan matrix, $A^0_{22}$ is a Jordan matrix with eigenvalue $\lambda+k$.
Note that the coefficient matrix $A$ reads
\begin{equation}
A=
\frac1x
\begin{pmatrix}
A^0_{11}+O(x) & O(x) \\
O(x) & A^0_{22}+O(x)
\end{pmatrix}.
\end{equation}
By transforming the dependent variable as $Y=SZ$ with
\begin{align}
S=
\begin{pmatrix}
I & O \\
O & x I
\end{pmatrix}
\end{align}
(shearing transformation),
the coefficient matrix of the system for $Z$ is
\begin{align}
S^{-1}AS-S^{-1}S'&=\frac1x
\begin{pmatrix}
A^0_{11}+O(x) & O(x^2) \\
O(1) & A^0_{22}-I+O(x)
\end{pmatrix}  \\
&=\frac1x
\begin{pmatrix}
A^0_{11} & O \\
* & A^0_{22}-I
\end{pmatrix}+O(1).
\nonumber
\end{align}
This indicates that the difference of eigenvalues of the leading term decreases by one.
By repeating this procedure, there will be no non-zero integral difference of eigenvalues of the leading term.

The above discussion is summarized as follows:
\begin{thm}\label{thm:first_red}
	For any system of the following form:
	\begin{align}\label{eq:first_kind}
	\frac{dY}{dx}&=\frac1x \left( A_0+A_1x+\cdots \right)Y, \quad
	A_j \in M_m(\mathbb{C}),
	\end{align}
	there exists a formal power series $P$ such that the transformation $Y=PZ$ 
	changes the equation~{\rm (\ref{eq:first_kind})} to
	\begin{align}
	\frac{dZ}{dx}=\frac{C}{x}Z
	\end{align}
	for some $C \in M_m(\mathbb{C})$.
	In particular, when no two distinct eigenvalues of $A_0$ differ by an integer, we can choose $C=A_0$.
\end{thm}

\subsubsection*{The case of pole order $> 1$}

Now we compute the canonical form of the following system with $r \ge 1$
\begin{equation}
\frac{dY}{dx}=AY, \quad
A=\frac{1}{x^{r+1}}\left( A_0+A_1x+\cdots+A_r x^r+\cdots \right).
\end{equation}
We choose the gauge so that $A_0$ has the form
\begin{equation}
A_0=
\begin{pmatrix}
A^0_{11} & O \\
O & A^0_{22}
\end{pmatrix},
\end{equation}
where $A^0_{11}$ and $A^0_{22}$ are Jordan matrices,
and assume that they satisfy
\begin{equation}\label{eq:ass_ev}
\{ \text{eigenvalue of $A^0_{11}$} \} \cap \{ \text{eigenvalue of $A^0_{22}$} \}=\emptyset.
\end{equation}

Transform the dependent variable $Y \to PY$ with
\begin{align}
P=I+\sum_{j=1}^\infty P_j x^j, \quad
P_j=
\begin{pmatrix}
O & P^j_{12} \\
P^j_{21} & O
\end{pmatrix},
\end{align}
where $P_j$'s are partitioned in the same manner as $A_0$.
The coefficient becomes
\begin{align}
B:=P^{-1}AP-P^{-1}\frac{dP}{dx}
:=\frac{B_0}{x^{r+1}}+\frac{B_1}{x^r}+\cdots+\frac{B_r}{x}+\cdots.
\end{align}
Rewriting this as $AP-P'=PB$ and equating like powers of $x$, we obtain the relations among $A_j, B_j, P_j$:
\begin{align}
B_0&=A_0, \\
A_0P_k-P_kB_0&=B_k-A_k-\sum_{j=1}^{k-1}(A_{k-j}P_j-P_j B_{k-j})+(k-r)P_{k-r} \quad (k \ge 1). \label{eq:red_irreg}
\end{align}
When $j<0$, we put $P_j=O$.
For simplification, we set
\begin{equation}
C_k=-A_k-\sum_{j=1}^{k-1}(A_{k-j}P_j-P_j B_{k-j})+(k-r)P_{k-r}.
\end{equation}
Since
\begin{equation}
A_0P_k-P_kB_0=
\begin{pmatrix}
O & A^0_{11}P^k_{12}-P^k_{12}A^0_{22} \\
A^0_{22}P^k_{21}-P^k_{21}A^0_{11} & O
\end{pmatrix},
\end{equation}
(1,1)-block and (2,2)-block of (\ref{eq:red_irreg}) reads
\begin{equation}\label{eq:red_irreg_diag}
O=B^k_{11}+C^k_{11}, \quad O=B^k_{22}+C^k_{22},
\end{equation}
(1,2)-block and (2,1)-block of (\ref{eq:red_irreg}) reads
\begin{equation}\label{eq:red_irreg_off-diag}
A^0_{11}P^k_{12}-P^k_{12}A^0_{22}=B^k_{12}+C^k_{12}, \quad
A^0_{22}P^k_{21}-P^k_{21}A^0_{11}=B^k_{21}+C^k_{21}.
\end{equation}
Here we have partitioned the matrices $B_k$ and $C_k$ into
\begin{equation}
B_k=
\begin{pmatrix}
B^k_{11} & B^k_{12} \\
B^k_{21} & B^k_{22}
\end{pmatrix}, \quad
C_k=
\begin{pmatrix}
C^k_{11} & C^k_{12} \\
C^k_{21} & C^k_{22}
\end{pmatrix}
\end{equation}
in the same manner as $A_0$.

Taking the assumption (\ref{eq:ass_ev}) and Lemma~\ref{thm:Sylvester_eq} into account,
we can choose $P^k_{12}, P^k_{21}$ so that $B^k_{12}$ and $B^k_{21}$ are the zero matrices.
Then we have
\begin{equation}
B_k=
\begin{pmatrix}
-C^k_{11} & O \\
O & -C^k_{22}
\end{pmatrix}.
\end{equation}

By repeating this procedure, we arrive at the following theorem.
\begin{thm}{(block diagonalization)}\label{thm:block_diag}
	For the following equation
	\begin{equation}\label{eq:irreg_sys}
	\frac{dY}{dx}=\frac{1}{x^{r+1}}\left( A_0+A_1x+\cdots \right)Y,
	\end{equation}
	there exists a formal power series $P$ such that the transformation of dependent variable $Y=PZ$ changes the equation~{\rm (\ref{eq:irreg_sys})} to
	the following form:
	\begin{align}
	&\frac{dZ}{dx}
	=
	\begin{pmatrix}
	B_1 & & \\
	& \ddots & \\
	& & B_n
	\end{pmatrix}Z, \quad
	B_k=\frac{1}{x^{r+1}}\left( B^k_0+B^k_1x+\cdots \right),
	\end{align}
	where the leading term $B^k_0$ of each block has only one eigenvalue.
\end{thm}
Therefore, equation (\ref{eq:irreg_sys}) is formally reduced to the following $n$ equations:
\begin{equation}
\frac{dZ_k}{dx}=B_k Z_k \quad (k=1,\ldots,n).
\end{equation}

\subsubsection*{Spectral type}\label{sec:spec_type}

Now we suppose $A_0$ is diagonalizable.
Then the leading term $B^k_0\ (k=1,\ldots,n)$ of each block $B_k$ is a scalar matrix.

Thus we focus on the next coefficient $B^k_1$.
If $B^k_1$ is also diagonalizable, we choose the gauge which diagonalizes $B^k_1$.
After that, we perform further block diagonalization corresponding to eigenvalues of $B^k_1$ in the same manner as
Theorem~\ref{thm:block_diag}.
Therefore the block $B_k$ again decomposes into smaller direct summands.
Notice that first two terms of each direct summands are scalar matrices.

In general, suppose that a direct summand of the following form
\begin{equation}
\frac{c_0 I}{x^{r+1}}+\frac{c_1 I}{x^{r}}+\cdots+\frac{c_{r-j}I}{x^{j+1}}+
\frac{B_{r-j+1}}{x^j}+\cdots
\end{equation}
appears and that $B_{r-j+1}$ is diagonalizable.
By means of a gauge transformation we have
\begin{equation}\label{eq:direct_sum}
\sim
\frac{c_0 I}{x^{r+1}}+\frac{c_1 I}{x^{r}}+\cdots+\frac{c_{r-j}I}{x^{j+1}}+
\frac{D_{r-j+1}}{x^j}+\cdots \quad
(D_{r-j+1}=d_1 I_{m_1} \oplus \cdots \oplus d_l I_{m_l}).
\end{equation}
By performing the block diagonalization corresponding to eigenvalues of $D_{r-j+1}$,
we can decompose (\ref{eq:direct_sum}) into smaller direct summands
\begin{equation}
\frac{c_0 I_{m_i}}{x^{r+1}}+\frac{c_1 I_{m_i}}{x^{r}}+\cdots+\frac{c_{r-j}I_{m_i}}{x^{j+1}}+
\frac{d_i I_{m_i}}{x^j}+\frac{*}{x^{j-1}}+\cdots
\quad (i=1,\ldots,l).
\end{equation}
Under the assumption that $B_{r-j+1}$ is diagonalizable for all $j=r, r-1, \ldots, 2$ and for all direct summands,
we can convert the equation (\ref{eq:irreg_sys}) to a direct sum of equations of the following form
\begin{equation}
\frac{dW}{dx}=\frac{1}{x^{r+1}}\left( C_0+C_1x+\cdots+C_{r-1}x^{r-1}+C_r x^r+\cdots \right)W \quad
\text{($C_0,\ldots,C_{r-1}$ are scalar matrices)}.
\end{equation}
Then by the Theorem~\ref{thm:first_red}, through a further gauge transformation with suitable formal power series,
we can eliminate the terms whose degree with respect to $x$ are greater than $-1$.
Thus we see that the equation (\ref{eq:irreg_sys}) formally transforms into the direct sum of the equations of the following form:
\begin{equation}
\frac{d\tilde{W}}{dx}=\frac{1}{x^{r+1}}\left( C_0+C_1x+\cdots+C_{r-1} x^{r-1}+\tilde{C}_r x^r \right)\tilde{W} \quad
\text{($C_0,\ldots,C_{r-1}$ are scalar matrices)}.
\end{equation}
The direct sum of the above form thus obtained is the canonical form in this case.
\begin{rem}
	When non-semisimple matrices appear in the course of the above calculation,
	we can reduce the case to the above case.
	To do so,
	we generally need to take an appropriate covering $x=\xi^k$ and to take shearing transformations
	with respect to $\xi$. 
	We have called the linear equation with such a singularity ramified type.
	\qed
\end{rem}

As we can see from the procedure to obtain the canonical form,
the canonical form has a nested structure;
the leftmost diagonal entries are divided into some groups, and in the next diagonal matrix,
each group is divided into some groups again, and so on. 
In this way, the structure of refining sequence of partitions of a canonical form emerges naturally.

\begin{df}[refining sequence of partitions]\label{def:rsp}
	Let $\lambda=\lambda_1\ldots \lambda_p$ and $\mu=\mu_1\ldots \mu_q$ be partitions of a natural number $m$:
	$\lambda_1+\cdots+\lambda_p=\mu_1+\cdots+\mu_q=m$．
	Here we assume that $\lambda_i$'s and $\mu_i$'s are not necessarily arranged in descending or ascending order.
	
	If there exist a disjoint decomposition $\{ 1,2,\ldots , p\}=I_1\coprod \cdots \coprod I_q$ of the index set of $\lambda$ such that $\mu_k=\sum_{j\in I_k}\lambda_j$ holds, then we call $\lambda$ a refinement of $\mu$.
	
	Let $[p_0,\ldots,p_r]$ be an $(r+1)$-tuple of partitions of $m$.
	When $p_{i+1}$ is a refinement of $p_i$ for all $i\ (i=0,\ldots,r-1)$, we call $[p_0,\ldots,p_r]$ a refining sequence of partitions.
\end{df}
\begin{eg}
	To the following canonical form
	\begin{align}
	{\scriptsize
		\begin{pmatrix}
		t^0_1 & & & & &  \\
		& t^0_1 & & & &  \\
		& & t^0_1 & & &  \\
		& & & t^0_2 & &  \\
		& & & & t^0_2 &  \\
		& & & & & t^0_3 
		\end{pmatrix}
		\frac{1}{x^3}+
		\begin{pmatrix}
		t^1_1 & & & & &  \\
		& t^1_1 & & & &  \\
		& & t^1_2 & & &  \\
		& & & t^1_3 & &  \\
		& & & & t^1_3 &  \\
		& & & & & t^1_4  \\
		\end{pmatrix}
		\frac{1}{x^2}+
		\begin{pmatrix}
		\theta_1 & & & & &  \\
		& \theta_2 & & & &  \\
		& & \theta_3 & & &  \\
		& & & \theta_4 & &  \\
		& & & & \theta_5 &  \\
		& & & & & \theta_6  \\
		\end{pmatrix}
		\frac{1}{x}
	},
	\end{align}
	we attach the following refining sequence of partitions
	\[
	[321, 2121, 111111].
	\]
	
	Here we introduce a convenient notation to express a refining sequence of partitions.
	First, write the rightmost partition
	\[
	111111.
	\]
	Second, put the numbers that are grouped together in the central partition in parentheses
	\[
	(11)(1)(11)(1).
	\]
	Finally, put the numbers that are grouped together in the leftmost partition in parentheses
	\[
	((11)(1))((11))((1)).
	\]
	Thus we can express the sequence as $((11)(1))((11))((1))$.
	\qed
\end{eg}
\begin{eg}
Let us look at another example
\[
[642,42411,2223111].
\]
First, write the rightmost partition
\[
2223111.
\]
Second, put the numbers that are grouped together in the central partition in parentheses
\[
(22)(2)(31)(1)(1).
\]
Finally, put the numbers that are grouped together in the leftmost partition in parentheses
\[
((22)(2))((31))((1)(1)).
\]
Thus the above sequence can be written as $((22)(2))((31))((1)(1))$.
\qed
\end{eg}


The refining sequence of partitions of $m$, or its abbreviated symbol,
that expresses the multiplicity of $t^i_j$'s and $\theta_j$'s of the canonical form
is called the spectral type at the singular point. 
A tuple of spectral types of all singular points (separated by ``\,,\,") is called
a spectral type of the equation.

When we express only the Poincar\'e rank of each singular point, 
we attach the number ``Poincar\'e rank plus 1" to each singular point and connect them with ``+".
We call it the singularity pattern of the equation.
At an unramified singularity, the number ``Poincar\'e rank plus 1" is equal to
the number of partitions in the refining sequence of partitions.
In this article, we use singularity patterns and spectral types to specify linear equations.

\begin{rem}
	The theory of isomonodromic deformation has been well studied since the work of Jimbo, Miwa, and Ueno \cite{JMU}.
	In their article, eigenvalues of leading terms of linear equations are assumed to be distinct.
	To classify Painlev\'e-type equations, however, we introduce a notion of spectral type, and also consider cases when  such eigenvalues are not necessarily distinct.
	If we consider examples such as generalized hypergeometric functions, it is natural to include the cases some eigenvalues are the same. 
	\qed
\end{rem}

\subsection{Degeneration scheme}
Now, we list these 22 types of equations in a degeneration scheme.
There are 26 boxes for Hamiltonians.
However, four of the Hamiltonians appear twice in the scheme, so we have 22 types.
Two linear equations with the same Hamiltonians can be converted into one another by Laplace transform, see section 4. 

In the case of a two dimensional phase space, the source equation of degeneration scheme,
which governs deformation of Fuchsian equation, is the sixth Painlev\'e equation.
Meanwhile, in the 4-dimensional case we have four source equations,
namely, the Garnier system in two variables, the $A_5^{(1)}$-type Fuji-Suzuki system,
the $D_6^{(1)}$-type Sasano system, and the sixth matrix Painlev\'e system \cite{S3}.
Accordingly, there are four series of degenerations.

Let us see the degeneration scheme.

\bigskip
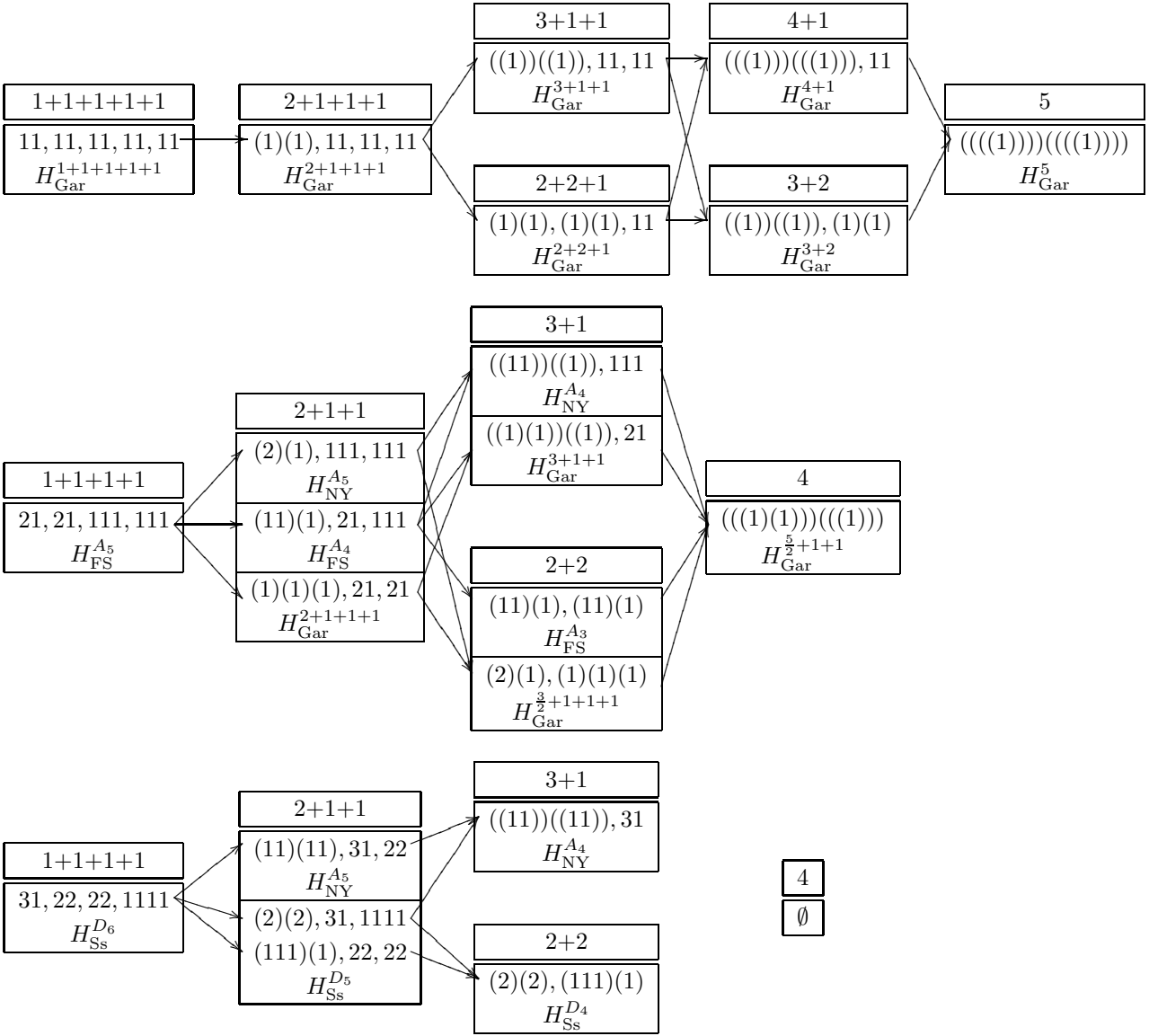
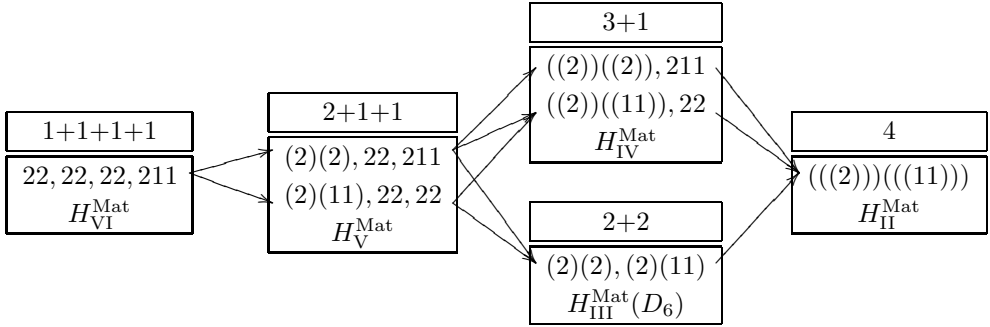
\begin{xy}
{(0,0) *{\begin{tabular}{|c|}
\hline
1+1+1+1+1\\
\hline\hline
$11,11,11,11,11$\\
$H_{\mathrm{Gar}}^{1+1+1+1+1}$\\
\hline
\end{tabular}
}},
{(35,0) *{\begin{tabular}{|c|}
\hline
2+1+1+1\\
\hline\hline
$(1)(1),11,11,11$\\
$H_{\mathrm{Gar}}^{2+1+1+1}$\\
\hline
\end{tabular}
}},
{\ar (12,0);(22,0)},
{(70,12) *{\begin{tabular}{|c|}
\hline
3+1+1\\
\hline\hline
$((1))((1)),11,11$\\
$H_{\mathrm{Gar}}^{3+1+1}$\\
\hline
\end{tabular}}},
{\ar (48,0);(56,12)},
{\ar (48,0);(56,-12)},
{(70,-12) *{\begin{tabular}{|c|}
\hline
2+2+1\\
\hline\hline
$(1)(1),(1)(1),11$\\
$H_{\mathrm{Gar}}^{2+2+1}$\\
\hline
\end{tabular}}},
{(105,12) *{\begin{tabular}{|c|}
\hline
4+1\\
\hline\hline
$(((1)))(((1))),11$\\
$H_{\mathrm{Gar}}^{4+1}$\\
\hline
\end{tabular}}},
{\ar (84,12);(90,12)},
{\ar (84,12);(90,-12)},
{\ar (84,-12);(90,12)},
{\ar (84,-12);(90,-12)},
{(105,-12) *{\begin{tabular}{|c|}
\hline
3+2\\
\hline\hline
$((1))((1)),(1)(1)$\\
$H_{\mathrm{Gar}}^{3+2}$\\
\hline
\end{tabular}}},
{\ar (120,12);(126,0)},
{\ar (120,-12);(126,0)},
{(140,0) *{\begin{tabular}{|c|}
\hline
5\\
\hline\hline
$((((1))))((((1))))$\\
$H_{\mathrm{Gar}}^5$\\
\hline
\end{tabular}}},\end{xy}

\bigskip
\begin{xy}
{(0,0) *{\begin{tabular}{|c|}
\hline
1+1+1+1\\
\hline\hline
$21,21,111,111$\\
$H_{\mathrm{FS}}^{A_5}$\\
\hline
\end{tabular}
}},
{\ar (12,-1);(22,10)},
{\ar (12,-1);(22,-1)},
{\ar (12,-1);(22,-12)},
{(35,0) *{\begin{tabular}{|c|}
\hline
2+1+1\\
\hline\hline
$(2)(1),111,111$\\
$H_{\mathrm{NY}}^{A_5}$\\
\hline
$(11)(1),21,111$\\
$H_{\mathrm{FS}}^{A_4}$\\
\hline
$(1)(1)(1),21,21$\\
$H_{\mathrm{Gar}}^{2+1+1+1}$\\
\hline
\end{tabular}
}},
{\ar (48,10);(56,22)},
{\ar (48,10);(56,-23)},
{\ar (48,-1);(56,22)},
{\ar (48,-1);(56,10)},
{\ar (48,-1);(56,-12)},
{\ar (48,-11);(56,10)},
{\ar (48,-11);(56,-23)},
{(70,18) *{\begin{tabular}{|c|}
\hline
3+1\\
\hline\hline
$((11))((1)),111$\\
$H_{\mathrm{NY}}^{A_4}$\\
\hline
$((1)(1))((1)),21$\\
$H_{\mathrm{Gar}}^{3+1+1}$\\
\hline
\end{tabular}}},
{(70,-18) *{\begin{tabular}{|c|}
\hline
2+2\\
\hline\hline
$(11)(1),(11)(1)$\\
$H_{\mathrm{FS}}^{A_3}$\\
\hline
$(2)(1),(1)(1)(1)$\\
$H_{\mathrm{Gar}}^{\frac{3}{2}+1+1+1}$\\
\hline
\end{tabular}}},
{(105,0) *{\begin{tabular}{|c|}
\hline
4\\
\hline\hline
$(((1)(1)))(((1)))$\\
$H_{\mathrm{Gar}}^{\frac{5}{2}+1+1}$\\
\hline
\end{tabular}}},
{\ar (84,22);(91,-1)},
{\ar (84,10);(91,-1)},
{\ar (84,-12);(91,-1)},
{\ar (84,-25);(91,-1)},
\end{xy}

\bigskip
\begin{xy}
{(0,0) *{\begin{tabular}{|c|}
\hline
1+1+1+1\\
\hline\hline
$31,22,22,1111$\\
$H_{\mathrm{Ss}}^{D_6}$\\
\hline
\end{tabular}
}},
{(35,0) *{\begin{tabular}{|c|}
\hline
2+1+1\\
\hline\hline
$(11)(11),31,22$\\
$H_{\mathrm{NY}}^{A_5}$\\ \hline
$(2)(2),31,1111$\\
$(111)(1),22,22$\\
$H_{\mathrm{Ss}}^{D_5}$\\
\hline
\end{tabular}
}},
{\ar (12,0);(22,8)},
{\ar (12,0);(22,-3)},
{\ar (12,0);(22,-8)},
{(70,12) *{\begin{tabular}{|c|}
\hline
3+1\\
\hline\hline
$((11))((11)),31$\\
$H_{\mathrm{NY}}^{A_4}$\\
\hline
\end{tabular}}},
{(70,-12) *{\begin{tabular}{|c|}
\hline
2+2\\
\hline\hline
$(2)(2),(111)(1)$\\
$H_{\mathrm{Ss}}^{D_4}$\\
\hline
\end{tabular}}},
{\ar (47,8);(57,12)},
{\ar (47,-3);(57,12)},
{\ar (47,-3);(57,-12)},
{\ar (47,-8);(57,-12)},
{(105,0) *{\begin{tabular}{|c|}
\hline
4\\
\hline\hline
$\emptyset$\\
\hline
\end{tabular}}},
\end{xy}

\bigskip
\begin{xy}
{(0,0) *{\begin{tabular}{|c|}
\hline
1+1+1+1\\
\hline\hline
$22,22,22,211$\\
$H_{\rm VI}^{\mathrm{Mat}}$\\
\hline
\end{tabular}
}},
{\ar (12,0);(23,3)},
{\ar (12,0);(23,-4)},
{(35,0) *{\begin{tabular}{|c|}
\hline
2+1+1\\
\hline\hline
$(2)(2),22,211$\\
$(2)(11),22,22$\\
$H_{\rm V}^{\mathrm{Mat}}$\\
\hline
\end{tabular}
}},
{\ar (47,3);(58,14)},
{\ar (47,3);(58,8)},
{\ar (47,3);(58,-12)},
{\ar (47,-4);(58,8)},
{\ar (47,-4);(58,-12)},
{(70,12) *{\begin{tabular}{|c|}
\hline
3+1\\
\hline\hline
$((2))((2)),211$\\
$((2))((11)),22$\\
$H_{\rm IV}^{\mathrm{Mat}}$\\
\hline
\end{tabular}}},
{(70,-12) *{\begin{tabular}{|c|}
\hline
2+2\\
\hline\hline
$(2)(2),(2)(11)$\\
$H_{{\rm III}}^{\mathrm{Mat}}(D_6)$\\
\hline
\end{tabular}}},
{\ar (82,14);(93,0)},
{\ar (82,8);(93,0)},
{\ar (82,-12);(93,0)},
{(105,0) *{\begin{tabular}{|c|}
\hline
4\\
\hline\hline
$(((2)))(((11)))$\\
$H_{\rm II}^{\mathrm{Mat}}$\\
\hline
\end{tabular}}},
\end{xy}
\bigskip

Here the symbol in the each upper box indicates singularity pattern.
There are one or more than one spectral types corresponding to each singularity pattern.
In the each lower box, we write the spectral types and corresponding Hamiltonian.
Explicit forms of Hamiltonians will be given in the next section.

We will see an example of degeneration process in section  \ref{sec:comment}.
All the data to calculate degeneration are listed in the appendix.

%

In this paper, we derive the list of Hamiltonians from those of Fuchsian-types using degeneration.
However, we
 can 
 also directly classify Painlev\'e-type equations corresponding to non-Fuchsian equations.
Indeed, there are also several interesting works attempting to construct a classification theory for non-Fuchsian equations \cite{A,B,B2,Hr,Kwk,Y}.
Especially, in hindsight, the work of Hiroe and Oshima \cite{HO} is very important to us.
From their results, it is proved that our list contains all the Painlev\'e-type equations corresponding to unramified linear equations.
In other words, irreducible unramified linear equations which admit non-trivial isomonodromic deformation can be transformed to one of the equations in our list by integral transformations.
One of the reasons we used degeneration is that it is the most efficient way to derive Hamiltonian explicitly, once we have Hamiltonians for Fuchsian cases.
Moreover, degeneration is often useful for analysis.



This article is organized as follows.
In the next section, we introduce 22 types of Hamiltonian systems.
In the third section, we show Lax pairs for Painlev\'e-type equations.
In the fourth section, we treat four topics:
 derivation of Hamiltonians in the case of Fuchsian type, procedure of degeneration, the Laplace transform, Fuchsian equations with three singular points. 
Detailed calculations of degenerations are written in the appendix.

We only deal with equations of unramified-type.
More equations are derived from ramified linear equations \cite{Kwk2}.

\bigskip 

\noindent 
{\bf Acknowledgments} 

\noindent 
We are deeply grateful to Prof.\ Oshima, Prof.\ Hiroe, whose enormous support and insightful comments were invaluable during the course of our study. 
We are also indebted to Prof.\ Okamoto, Prof.\ Ohyama, Prof.\ Suzuki, Prof.\ Tsuda, Prof.\ Jimbo, Prof.\ Nishioka, Prof.\ Takemura, Prof.\ Yamakawa, who gave us invaluable comments and warm encouragements.
This work was partially supported by Grant-in-Aid no.24540205,  
the Japan Society for the Promotion of Science. 

\section{List of Hamiltonians}

There are 22 Hamiltonian systems that we treat in this paper.
However these do not include the case that associated linear
equations have singularities of ramified type.
We will deal with only the unramified cases in this paper.
We expect a further research including ramified case.
Let us see the
expressions of the 22 Hamiltonians.

In the first place we look at the Garnier system and degenerate
Garnier systems, which are classical systems found in the early 20th century.
On the degeneration scheme of this family, a detailed study is now well
known \cite{Km}, though we introduce new expressions for some of these
systems.

The first system is the Garnier system, which is obtained from
a deformation of a Fuchsian equation with 5 regular singular points \cite{G}.
In the original paper by Garnier, the dependent variables are the
positions of apparent singular points.
However the equation does not have Painlev\'e property, that is, the
solution has movable algebraic singularities.
H. Kimura and Okamoto used  symmetric functions of apparent singularities as dependent variables, so that the Hamiltonian system enjoys the Painlev\'e property \cite{KO}.
The following Hamiltonian coincides with theirs, although they do not use the Hamiltonian of the
sixth Painlev\'e equation:
\begin{align}
 &t_i(t_i-1)H_{\mathrm{Gar}, t_i}^{1+1+1+1+1}
\left({\alpha, \beta \atop \gamma_1, \gamma_2, \delta}
;{t_1\atop t_2};{q_1,p_1 \atop q_2,p_2}\right)\\
&=t_i(t_i-1)H_{\rm VI}\left({\delta,\ \beta \atop
\gamma_{i}, \alpha+\gamma_{i+1}+1} ; t_i ; q_i,p_i\right)
+(2q_ip_i+q_{i+1}p_{i+1}-\beta-2\delta)q_1q_2p_{i+1}\nonumber\\
&\quad-\frac{1}{t_i-t_{i+1}}
\{ t_i(t_i-1)(p_iq_i+\gamma_{i})p_iq_{i+1}-t_i(t_{i+1}-1)(2p_iq_i+\gamma_{i})p_{i+1}q_{i+1}\nonumber\\
&\quad+t_{i+1}(t_i-1)({p_{i+1}}^2q_{i+1}+\gamma_{i+1}(p_{i+1}-p_i))q_i\} ,
\qquad (i\in \mathbb{Z}/2\mathbb{Z}).\nonumber
\end{align}
Relation between canonical variables and parameters of
associated linear equations will be explained later.

The second is the degenerate Garnier system obtained by a confluence of
two regular singular points;
\begin{align}
H_{\mathrm{Gar},\ t_1}^{2+1+1+1}\left({\alpha,\beta\atop\gamma,\delta};
{t_{1} \atop t_{2}};{q_{1},p_{1}\atop q_{2},p_{2}}\right)&=
H_{\rm{V}}\left({-\alpha-\gamma,-\beta-\gamma-\delta-1 \atop \alpha+\beta+\gamma};t_{1};q_{1},p_{1}\right) \label{eqn:2+1+1+1_1}\\
&\quad+\frac{p_{1}}{t_{1}}[\gamma(q_{1}-q_{2})+p_{2}q_{2}(q_{2}-1)]\nonumber\\
&\quad+\frac{1}{t_{1}-t_{2}}\left((q_{1}-q_{2})p_{1}-\beta
 \right)\left((q_2-q_1)p_2-\gamma\right) ,\nonumber\\
H_{\mathrm{Gar},\ t_2}^{2+1+1+1}\left({{\alpha,\beta\atop\gamma,\delta}};
{t_{1}\atop t_{2}};{q_{1},p_{1} \atop q_{2},p_{2}}\right)&=
H_{\rm{V}}\left({-\alpha-\beta,-\beta-\gamma-\delta-1 \atop \alpha+\beta+\gamma};t_{2};q_{2},p_{2}\right)\label{eqn:2+1+1+1_2}\\
&\quad+\frac{p_{2}}{t_{2}}[\beta(q_{2}-q_{1})+p_{1}q_{1}(q_{1}-1)]\nonumber\\
&\quad+\frac{1}{t_{2}-t_{1}}\left((q_{1}-q_{2})p_{1}-\beta\right)\left((q_{2}-q_{1})p_{2}-\gamma\right)
 .\nonumber
\end{align}
\begin{rem}
These canonical variables are different form the hitherto
 known \cite{Km}.
The original Hamiltonians were written as
\begin{align*}
{s_1}^2H_1&={\lambda_1}^2(\lambda_1-s_1){\mu_1}^2+2{\lambda_1}^2\lambda_2\mu_1\mu_2+
\lambda_1\lambda_2(\lambda_2-s_2){\mu_2}^2\\
&\quad-\{(\kappa_0+\theta_2-1){\lambda_1}^2+\kappa_1\lambda_1(\lambda_1-s_1)+\eta(\lambda_1-s_1)+
\eta s_1 \lambda_2\}\mu_1\\
&\quad-\{(\kappa_0+\kappa_1-1)\lambda_1\lambda_2+\theta_2\lambda_1(\lambda_2-s_2)-\eta(s_2-1)\lambda_2\}\mu_2
+\kappa \lambda_1,\\
s_2(s_2-1)H_2&={\lambda_1}^2\lambda_2{\mu_1}^2+2\lambda_1\lambda_2(\lambda_2-s_2)\mu_1\mu_2\\
&\quad+\left\{\lambda_2(\lambda_2-1)(\lambda_2-s_2)+\frac{s_2(s_2-1)}{s_1}\lambda_1\lambda_2\right\} {\mu_2}^2\\
&\quad-\{(\kappa_0+\kappa_1-1)\lambda_1\lambda_2+\theta_2\lambda_1(\lambda_2-s_2)
-\eta(s_2-1)\lambda_2\}\mu_1\\
&\quad-\biggl\{(\kappa_0-1)\lambda_2(\lambda_2-1)+\kappa_1\lambda_2(\lambda_2-s_2)
+\theta_2(\lambda_2-1)(\lambda_2-s_2)\biggr.\\
&\left.\quad+\frac{s_2(s_2-1)}{s_1}(\theta_2\lambda_1+\eta
 \lambda_2)\right\}\lambda_2+\kappa \lambda_2 .
\end{align*}
This expression is not symmetric in the canonical variables, while the Hamiltonians (\ref{eqn:2+1+1+1_1})--(\ref{eqn:2+1+1+1_2}) are
  symmetric.
Besides, they are expressed simply by using the Hamiltonian of the fifth Painlev\'e
 equation.
The other degenerate Garnier systems below are expressed in a similar fashion.
\qed
\end{rem}

The next Hamiltonians are ones of the degenerate Garnier system
associated to a linear equation with two irregular singular points and one regular singular point:
\begin{align}
&t_1H^{2+2+1}_{\mathrm{Gar},t_1}\left({\alpha,\beta \atop \gamma};{t_1 \atop t_2};
{q_1,p_1 \atop q_2,p_2}\right)\\
&=
t_1H_{\rm{V}}
\left({\alpha+\beta+\gamma,\beta-\alpha
\atop -\beta-\gamma};t_1;q_1,p_1\right)\nonumber\\
&\quad+q_1q_2(p_1q_1-\alpha)+p_2q_2(\alpha+p_1-2p_1q_1)
-\frac{t_2}{t_1}p_1(p_2-q_1) ,\nonumber\\
&t_2H^{2+2+1}_{\mathrm{Gar},t_2}\left({\alpha,\beta \atop \gamma};{t_1 \atop t_2};
{q_1,p_1 \atop q_2,p_2}\right)\\
&=
t_2H_{\rm{III}}(D_6)
\left({-\alpha-\beta-\gamma,-\beta}
;t_2;q_2,p_2\right)
-(p_1q_1-\alpha)q_2(q_1-1)+\frac{t_2}{t_1}p_1(p_2-q_1) .\nonumber
\end{align}

The following system is associated to a linear equation with two regular
singular points and one irregular singular point:
\begin{align}
&H^{3+1+1}_{\mathrm{Gar},t_1}
\left({\alpha,\beta \atop \gamma};{t_1 \atop t_2};
{q_1,p_1 \atop q_2,p_2}\right)\\
&=
H_{\rm{IV}}\left({\alpha,\gamma};t_1;q_1,p_1\right)
+p_2q_2p_1+\frac{1}{t_1-t_2}\{p_1(q_1-q_2)-\alpha\}\{p_2(q_2-q_1)
-\beta\} ,\nonumber\\
&H^{3+1+1}_{\mathrm{Gar},t_2}
\left({\alpha,\beta \atop \gamma};{t_1 \atop t_2};
{q_1,p_1 \atop q_2,p_2}\right)\\
&=
H_{\rm{IV}}\left({\beta,\gamma};t_2;q_2,p_2\right)
+p_1q_1p_2+\frac{1}{t_2-t_1}\{p_1(q_1-q_2)-\alpha\}\{p_2(q_2-q_1)
-\beta\} .\nonumber
\end{align}

The following system is the degenerate Garnier system associated to a
linear equation with only two irregular singular points:
\begin{align}
&H^{3+2}_{\mathrm{Gar},t_1}
\left({\alpha,\beta};{t_1 \atop t_2};
{q_1,p_1 \atop q_2,p_2}\right)
=
H_{\rm{III}}(D_6)\left({-\beta,\alpha+1};t_1;q_1,p_1\right)
-p_1-\frac{q_1q_2}{t_1}(q_2-p_2+t_2)+p_1p_2-q_2,\\
&H^{3+2}_{\mathrm{Gar},t_2}
\left({\alpha,\beta};{t_1 \atop t_2};
{q_1,p_1 \atop q_2,p_2}\right)
=
H_{\rm{IV}}\left({\alpha,\beta};t_2;q_2,p_2\right)
-p_1q_1(p_2-2q_2-t_2)-q_1q_2+t_1p_1.
\end{align}

The next system is associated to a linear equation with one regular
singular point and one irregular singular point:
\begin{align}
&H^{4+1}_{\mathrm{Gar},t_1}
\left({\alpha, \beta};{t_1 \atop t_2};
{q_1,p_1 \atop q_2,p_2}\right)
=
H_{\rm{II}}\left({-\beta};t_1;q_1,p_1\right)
+p_2q_2(q_1-q_2+t_2)+p_1p_2+\alpha q_2,\\
&H^{4+1}_{\mathrm{Gar},t_2}
\left({\alpha, \beta};{t_1 \atop t_2};
{q_1,p_1 \atop q_2,p_2}\right)\\
&=
-{p_2}^2q_2-t_2p_2{q_2}^2+{t_2}^2p_2q_2+\alpha t_2q_2-\beta p_2
+p_1p_2(q_1-2q_2+t_2)+q_1q_2(p_2q_2-\alpha)+\alpha p_1+t_1p_2q_2\nonumber.
\end{align}
\begin{rem}
The second Hamiltonian $H_{\mathrm{Gar},t_1}^{4+1}$ are not represented by a Hamiltonian of the Painlev\'e equation.
When we change the variables as
\[
t_2 \to t_2^{\frac{2}{3}}, \quad q_2 \to t_2^{-\frac{1}{3}}p_2,\quad p_2 \to -t_2^{\frac13}q_2,
\]
the second Hamiltonian can be expressed by the Hamiltonian of the fourth Painlev\'e equation:
\begin{align}
H^{4+1}_{\mathrm{Gar},t_2}
\left({\alpha, \beta};{t_1 \atop t_2};
{q_1,p_1 \atop q_2,p_2}\right)&=\frac23\big[ H_{\rm IV}\left(\beta,
 \alpha;t_2;q_2,p_2\right)-\frac{p_2q_2}{2t_2}-p_1q_2(q_1-2t_2^{-\frac13}p_2+t_2^{\frac23})\\
&\quad -t_2^{-\frac23}q_1p_2(p_2q_2+\alpha)
+\alpha t_2^{-\frac13}p_1-t_1t_2^{-\frac13}p_2q_2\big].\nonumber
\end{align}
Here the fourth Painlev\'e Hamiltonian appears.
\qed
\end{rem}

Concerning a linear equation with only one irregular singular point,
which is obtained by a confluence of all singularities, we have the
following system:
\begin{align}
&H_{\mathrm{Gar},t_1}^5
\left({\alpha};{t_1 \atop t_2};
{q_1,p_1 \atop q_2,p_2}\right)
=
-q_{1}(q_1 p_1-\alpha)
+q_2(q_1(p_2+q_2)-2p_1+t_1)
+p_1(p_2-2t_2),\\
&H_{\mathrm{Gar},t_2}^5
\left({\alpha};{t_1 \atop t_2};
{q_1,p_1 \atop q_2,p_2}\right)=
H_{\rm IV}\left(-1,\alpha;2t_2;q_2,p_2\right)
+q_1q_2(q_1q_2-2p_1+t_1)
+p_1(p_1-p_2q_1-t_1).
\end{align}

We have 7 systems which have two independent variables so far.
In Kimura's original paper \cite{Km}, there are 8 systems, but the last system denoted as $(9/2)$
is obtained from a deformation of a linear equation with singular point of ramified type;
so we put aside this one.

In the degeneration scheme, further two degenerate Garnier systems appear.
These two systems are associated to linear systems of rank two with
singular points of ramified type but these are also associated to unramified linear systems of rank three.
These rank three systems are obtained by confluence of singularities from
Fuji-Suzuki system.
With respect to degenerate Garnier systems of ramified type, all of them are
obtained by H. Kawamuko and there are 9 systems in all \cite{Kw}.

The first one of them is associated to a linear equation with two
irregular singular points, which is obtained from degereration of linear system of rank three associated to Fuji-Suzuki system:
\begin{align}
&t_1H_{\mathrm{Gar},t_1}^{\frac{3}{2}+1+1+1}
\left({\alpha,\beta\atop \gamma};{t_1 \atop t_2};
{q_1,p_1 \atop q_2,p_2}\right)\\
&=
t_1H_{\rm III}(D_6)\left({-\alpha,\gamma-\alpha};t_1;q_1,p_1\right)
+q_1(q_1p_1p_2-\alpha p_2)
+\frac{t_1}{t_1-t_2}(p_1(q_1-q_2)-\alpha)(p_2(q_2-q_1)-\beta),\nonumber\\
&t_2H_{\mathrm{Gar},t_2}^{\frac{3}{2}+1+1+1}
\left({\alpha,\beta \atop \gamma};{t_1 \atop t_2};
{q_1,p_1 \atop q_2,p_2}\right)\\
&=
t_2H_{\rm III}(D_6)\left({-\beta,\gamma-\beta};t_2;q_2,p_2\right)
+q_2(q_2p_1p_2-\beta p_1)
+\frac{t_2}{t_2-t_1}(p_1(q_1-q_2)-\alpha)(p_2(q_2-q_1)-\beta).\nonumber
\end{align}

The second one is associated to a linear equation possessing only one
irregular singular point, which is obtained by confluence of all
singular points from a linear equation associated to Fuji-Suzuki system:
\begin{align}
&H_{\mathrm{Gar},t_1}^{\frac{5}{2}+1+1}
\left({\alpha,\beta};{t_1 \atop t_2};
{q_1,p_1 \atop q_2,p_2}\right)\\
&=
H_{\rm II}\left(-\alpha;t_1;q_1,p_1\right)
+p_1p_2
+\frac{1}{t_1-t_2}(p_1(q_1-q_2)-\alpha)(p_2(q_2-q_1)-\beta),\nonumber\\
&H_{\mathrm{Gar},t_2}^{\frac{5}{2}+1+1}
\left({\alpha,\beta};{t_1 \atop t_2};
{q_1,p_1 \atop q_2,p_2}\right)\\
&=
H_{\rm II}\left(-\beta;t_2;q_2,p_2\right)
+p_1p_2
+\frac{1}{t_2-t_1}(p_1(q_1-q_2)-\alpha)(p_2(q_2-q_1)-\beta).\nonumber
\end{align}
\begin{rem}
The Garnier equations corresponding to ramified linear equations in Kawamuko's list \cite{Kw} are as follows: 
$H_{\rm Gar}^{\frac{3}{2}+1+1+1}$, $H_{\rm Gar}^{\frac{3}{2}+\frac{3}{2}+1}$, $H_{\rm Gar}^{\frac{3}{2}+2+1}$, $H_{\rm Gar}^{\frac{5}{2}+1+1}$, $H_{\rm Gar}^{\frac{5}{2}+\frac{3}{2}}$, $H_{\rm Gar}^{\frac{3}{2}+3}$, $H_{\rm Gar}^{\frac{5}{2}+2}$, $H_{\rm Gar}^{\frac{7}{2}+1}$, $H_{\rm Gar}^{\frac{9}{2}}$.
Among these 9 ramified Garnier equations, two of them also correspond to unramified linear equations of rank three.
\qed
\end{rem}	

We have already seen the 9 partial differential systems with two independent
variables so far, and we will see the 13 ordinary differential systems below.

We have three more systems which are obtained from the isomonodromic
deformation of Fuchsian equations, apart from the Garnier system.
These are the Fuji-Suzuki system, the Sasano system, and the system which we call the sixth matrix Painlev\'e
system.
We begin with the Fuji-Suzuki system.
It was obtained from the similarity reduction of a Drinfel'd-Sokolov
hierarchy  \cite{FS1}:
\begin{align}
&t(t-1)H_{\mathrm{FS}}^{A_5}\left({\alpha, \beta, \gamma \atop \delta, \epsilon, \omega};t;{q_1, p_1 \atop q_2, p_2}\right)\\
&=t(t-1)H_{\rm{VI}}\left({\alpha,\beta+\delta \atop \beta+\gamma,\epsilon-\omega+1};t;q_1,p_1\right)
+t(t-1)H_{\rm{VI}}\left({\beta,\alpha+\delta \atop \alpha+\gamma,\epsilon-\alpha+1};t;q_2,p_2\right)\nonumber\\
&\quad+(q_1-t)(q_2-1)\{(p_1q_1-\alpha)p_2+p_1(p_2q_2-\beta)\}.\nonumber
\end{align}
\begin{rem}
Before the result of Fuji and Suzuki, Tsuda calculated
 isomonodromic deformations with respect to a class of Fuchsian
 equations including 21,21,111,111 from a viewpoint of reduction theory
 of his UC hierarchy.
Although his system of equations was not written
in the form of the Hamiltonian system, it was found that his system
 includes the Fuji-Suzuki's coupled sixth Painlev\'e system in the
 paper  \cite{S3}.
In the aftermath Tsuda gave a Hamiltonian expression to the whole systems in
 his class \cite{Ts}.
In this paper we call this system the Fuji-Suzuki system because they
were the first ones who gave the expression of the coupled sixth Painlev\'e system.
\qed
\end{rem}

A degeneration from the above Fuji-Suzuki system produces the following
Hamiltonian:
\begin{align}
&tH_{\mathrm{FS}}^{A_4}\left({\alpha,\beta,\gamma\atop \delta,\epsilon};t;{q_1,p_1 \atop q_2,p_2}\right)\\
&=tH_{\rm{V}}\left({\alpha+\beta+\delta+\epsilon,\alpha+\gamma-\delta-1\atop -\alpha-\epsilon};t;q_1,p_1\right)
+tH_{\rm{V}}\left({\alpha+\epsilon,\alpha+\gamma-\epsilon-1\atop-\alpha};t;q_2,p_2\right)\nonumber\\
&\quad+p_1(q_2-1)(p_2(q_1+q_2)-\epsilon).\nonumber
\end{align}

Further degeneration also produces the following Hamiltonian:
\begin{align}
tH_{\mathrm{FS}}^{A_3}\left({\alpha,\beta \atop \gamma,\delta};t;{q_1,p_1 \atop q_2,p_2}\right)
&=
tH_{\rm{III}}(D_6)\left(\alpha+\gamma,-\beta+\gamma;t;q_1,p_1\right)
+tH_{\rm{III}}(D_6)\left(\delta,-\beta+\delta;t;q_2,p_2\right)\\
&\quad+p_1q_2(p_2(q_1+q_2)+\delta).\nonumber
\end{align}
We call Hamiltonian systems of $H_{\mathrm{FS}}^{A_4}$ and $H_{\mathrm{FS}}^{A_3}$ degenerate Fuji-Suzuki systems.
The degenerate Fuji-Suzuki systems $H_{\mathrm{FS}}^{A_4}$, $H_{\mathrm{FS}}^{A_3}$ have particular solutions expressed by generalized hypergeometric functions~\cite{Sz,Ts2}.

The famous Noumi-Yamada systems are written in the form of Hamiltonian
system by using the following Hamiltonians:
\begin{align}\
t H_{\mathrm{NY}}^{A_5}\left({\alpha, \beta, \gamma \atop \delta, \epsilon};t;{q_1,p_1 \atop q_2,p_2}\right)&=
tH_{\rm{V}}\left({\alpha+\beta, \alpha+\gamma+\epsilon \atop -\alpha};t;q_1,p_1\right)
+tH_{\rm{V}}\left({\alpha+\gamma+\delta, \alpha+\gamma+\epsilon \atop -\alpha-\gamma};t;q_2,p_2\right)\\
&\quad+2 p_1p_2q_1(q_2-1),\nonumber \\
H_{\mathrm{NY}}^{A_4}\left({\alpha, \beta \atop \gamma, \delta};t;{q_1,p_1 \atop q_2,p_2}\right)&=
H_{\rm{IV}}(\beta, \alpha;t;q_1,p_1)+
H_{\rm{IV}}(\delta, \alpha+\gamma;t;q_2,p_2)+2q_1p_1p_2 .
\end{align}
These two systems are obtained by degenerations from both the Fuji-Suzuki
system and the Sasano system.
The Noumi-Yamada systems were found under search for nonlinear equations with the symmetry group $A_{n}^{(1)}$.
It is important that the Painlev\'e-type equations such as Noumi-Yamada systems with quite different origins are organized in our unified way in this paper. 
\begin{rem}
The Noumi-Yamada systems are well known in the following expression:
 \begin{align}
  NY^{A_4}~:~&\frac{df_i}{dt}=f_i(f_{i+1}-f_{i+2}+f_{i+3}-f_{i+4})+\alpha_i,
  \quad (i\in {\mathbb Z}/5{\mathbb Z}),\label{eq:NY_4}\\
  NY^{A_5}~:~&\frac{df_i}{dt}=f_i(f_{i+1}f_{i+2}-f_{i+2}f_{i+3}+f_{i+1}f_{i+4}-f_{i+2}f_{i+5}
  +f_{i+3}f_{i+4}-f_{i+4}f_{i+5})\label{eq:NY_5}\\
  &+(-1)^i(\alpha_{i+1}+\alpha_{i+2}+\alpha_{i+5})f_i+\alpha_i(f_i+f_{i+2}+f_{i+4}),
  \quad (i\in {\mathbb Z}/6{\mathbb Z}),\nonumber
 \end{align}
as systems of equations for the unknown function $f_0,\ldots ,
 f_l,(l=4~\mbox{or}~5)$ \cite{NY}.
These systems coincide with the Hamiltonians above by putting
$p_1=f_2$, $q_1=-f_1$, $p_2=f_4$, and $q_2=-f_1-f_3$.
Here the parameters are $\alpha=-\alpha_1$, $\beta=-\alpha_2$,
 $\gamma=-\alpha_3$, $\delta=-\alpha_4$, $\epsilon=-\alpha_5$ \cite{SY}.
 \qed
\end{rem}
\begin{rem}
Adler \cite{Ad} and Veselov-Shabat \cite{VS} also studied these equations independently prior to Noumi-Yamada \cite{NY}.
In their terminology, equations (\ref{eq:NY_4}), (\ref{eq:NY_5}) are Darboux chain with period $5$ and $6$, respectively. 
\qed
\end{rem}

We introduce another coupled Painlev\'e system called Sasano system.
Sasano systems were obtained by a generalization of space of initial
values for the Painlev\'e equations \cite{Ss},
and it is also derived from similarity reduction of Drinfel'd-Sokolov
hierarchy \cite{FS2}:
\begin{align}
H_{\mathrm{Ss}}^{D_6}\left({\alpha, \beta, \gamma \atop \delta, \epsilon, \zeta};t;
{q_1,p_1 \atop q_2,p_2}\right)&=H_{\rm{VI}}\left({\beta+\gamma+2\delta+\epsilon+\zeta, -\beta-\zeta \atop -\beta-2\gamma-2\delta-\epsilon,1-\alpha-\beta-2\delta-\epsilon-\zeta};t;q_1,p_1\right)\\
&\quad+H_{\rm{VI}}\left({\gamma+\delta, \epsilon \atop \zeta, 1-\alpha-\gamma};t;q_2,p_2\right)\nonumber\\
&\quad+\frac{2}{t(t-1)}(q_1-1)p_{2} q_{2}\{(q_{1}-t)p_{1}-\beta-\gamma-2\delta-\epsilon-\zeta)\} .\nonumber
\end{align}
This system is different from the Fuji-Suzuki system of type $A_5^{(1)}$
in its ``coupling term''.

Degenerations give the following Hamiltonians:
\begin{align}
tH_{\mathrm{Ss}}^{D_5}
 \left({\alpha,\beta,\gamma \atop \delta, \epsilon};t;
{q_1,p_1\atop q_2,p_2}\right)
&=tH_{\rm V}\left({\epsilon,\alpha-\beta\atop
\beta};t;q_1,p_1\right)\\
&\quad +tH_{\rm V}
 \left({-2\alpha-3\beta-\gamma-\delta-2\epsilon,-\alpha-\beta-2\delta\atop \alpha+2\beta+\delta+\epsilon};t;q_2,p_2\right)\nonumber\\
&
\quad
+2p_{2}q_{1}(p_{1}(q_{1}-1)-\beta-\epsilon)
, \nonumber \\
tH_{\mathrm{Ss}}^{D_4}
 \left({\alpha,\beta \atop \gamma, \delta};t;
{q_1,p_1\atop q_2,p_2}\right)
&=tH_{\rm III}(D_{6})
 \left(\alpha+\beta+\gamma, -\alpha-2\delta 
;t;q_1,p_1\right)\\
&\quad+tH_{\rm III}(D_{6})\left(-\gamma,
-\alpha-2\gamma;t;q_2,p_2\right)
+2p_{2}q_1(p_{1}q_{1}+\alpha+\beta+\gamma).\nonumber
\end{align}
In this article, we call these two systems degenerate Sasano systems.

At the end we introduce five Hamiltonian systems which we call matrix
Painlev\'e systems:
\begin{multline}
t(t-1)H^{\mathrm{Mat}}_{\mathrm{VI}}\left({\alpha , \beta , \gamma\atop\delta ,\omega};t;
{q_1, p_1\atop q_2, p_2}\right)=
\mathrm{tr}\left[Q(Q-1)(Q-t)P^2+
\{(\delta -(\alpha -\omega )K)Q(Q-1)\right.\\
\left. +\gamma (Q-1)(Q-t)-(2\alpha +\beta +\gamma +\delta )Q(Q-t)\}P
+\alpha (\alpha +\beta )Q\right],
\end{multline}
\begin{align}
tH^{\mathrm{Mat}}_{\mathrm{V}}\left({\alpha,\beta\atop\gamma ,\omega};t;{q_1, p_1\atop
 q_2, p_2}\right)&=
\mathrm{tr}[Q(Q-1)P(P+t)+\beta QP+\gamma P-(\alpha +\gamma ) tQ],\\
H^{\mathrm{Mat}}_{\mathrm{IV}}\left(\alpha,\beta , \omega ;t;{q_1, p_1\atop q_2, p_2}\right)
 &=\mathrm{tr}[QP(P-Q-t)+\beta P+\alpha Q],\\
tH^{\mathrm{Mat}}_{\mathrm{III}(D_6)}\left(\alpha,\beta , \omega ;t;{q_1, p_1\atop q_2, p_2}\right)
 &=\mathrm{tr}[Q^2P^2-(Q^2-\beta Q- t)P-\alpha Q],\\
H^{\mathrm{Mat}}_{\mathrm{II}}\left(\alpha , \omega ;t;{q_1, p_1\atop q_2, p_2}\right)
 &=\mathrm{tr}[P^2-(Q^2+t)P-\alpha Q].
\end{align}
Here the matrices $P$ and $Q$ satisfy the relation $[P,Q]=(\alpha
-\omega )K \quad (K={\rm diag}(1,-1))$, and canonical variables can be
written as
\[
 p_1={\rm tr}P,\quad q_1=\frac12{\rm tr}Q,\quad
 p_2=-\frac{P_{12}}{Q_{12}},\quad q_2=-Q_{12}Q_{21}.
\]
We denote $(1,2)$-component of matrix $P$ as $P_{12}$, and so on.
Using these canonical variables, the above Hamiltonians can be rewritten as follows:
\begin{align*}
t(t-1)H_{\mathrm{VI}}^{\mathrm{Mat}}\left({\alpha , \beta , \gamma \atop \delta ,\omega};t;
{q_1, p_1\atop q_2, p_2}\right)
&=2t(t-1)H_{\mathrm{VI}}\left({\alpha, \beta \atop \gamma, \delta};t;q_1,\frac{p_1}{2}\right)
+(2q_1-1)({p_2}^2{q_2}^2+(p_2q_2-\alpha +\omega)^2)\\
&\quad-\frac{1}{2}p_1q_2((3q_1-t-1)p_1-4\alpha -2\beta)\\
&\quad-2(q_1-t)({q_1}^2-q_1-q_2)p_2(p_2q_2-\alpha +\omega)\\
&\quad+\{ \{3{q_1}^2-2(t+1)q_1-q_2+t\}p_1-(2\alpha +\beta)(2q_1-t)
-(1-t)\delta -\gamma\}\\
&\quad\quad\times (2p_2q_2-\alpha +\omega), \\
tH_{\mathrm{V}}^{\mathrm{Mat}}\left({\alpha,\beta\atop\gamma ,\omega};t;{q_1, p_1\atop q_2, p_2}\right)
&=2tH_{\mathrm{V}}\left({\alpha,-2\alpha+\beta+2\omega \atop	\alpha+\gamma-\omega};t;q_1,\frac{p_1}{2}\right)
   -2 p_2({q_1}^2-q_1-q_2)(p_2q_2-\alpha +\omega)\\
   &\quad-p_1q_2\left( \frac{p_1}{2}+t \right)
   +2p_2q_2\{(p_1+t)(2q_1-1)+\beta\}+(\alpha -\omega)(t-\beta),\\
H_{\mathrm{IV}}^{\mathrm{Mat}}\left({\alpha,\beta,\omega};t;{q_1, p_1\atop q_2, p_2}\right)
&=2H_{\mathrm{IV}}\left(2\alpha-\omega,-\alpha+\beta+\omega;t;q_1,\frac{p_1}{2}\right)
-2(p_2q_1+t)(p_2q_2-\alpha +\omega)\\
&\quad+2p_2q_2(p_1-2q_1)+p_1q_2-(\alpha -\omega)t,\\
tH_{{\mathrm{III}}(D_6)}^{\mathrm{Mat}}\left({\alpha,\beta,\omega};t;{q_1, p_1\atop q_2, p_2}\right)
&=2tH_{\mathrm{III}}(D_6)\left(\omega, -2\alpha+\beta+2\omega;t;q_1,\frac{p_1}{2}\right)
+2(2p_1q_1-2q_1+\beta)p_2q_2\\
&\quad+p_1q_2\left( 1-\frac{p_1}{2} \right)
-2 p_2\left( {q_1}^2-q_2\right)(p_2q_2-\alpha +\omega)-(\alpha-\omega)\beta,\\
H_{\mathrm{II}}^{\mathrm{Mat}}\left(\alpha , \omega
;t;{q_1, p_1\atop q_2, p_2}\right)
&=2H_{\mathrm{II}}\left(\omega ;t; q_1,\frac{p_1}{2}\right)
-4p_2q_1q_2+p_1q_2-2p_2\left( p_2q_2-\alpha +\omega\right).
\end{align*}

\begin{rem}
The matrix Painlev\'e systems were found independently by P. Boalch.
At the conference in Tokyo 2008 when one of the authors showed the Hamiltonian of matrix $P_{\rm VI}$ for the first time, P. Boalch talked about the correspondence between the moduli of the linear connections and the quiver varieties.
As an example of four-dimensional case, he showed a quiver with the same Dynkin diagram as the $H_{\rm VI}^{\mathrm{Mat}}$.
In the terminology of \cite{B2}, the matrix Painlev\'e equations are called ``higher Painlev\'e systems".
For the details, see \cite{B2,B3}.
\qed
\end{rem}

\section{Isomonodromic deformation}\label{sec:deformation}
In this section, we present both Painlev\'e-type equations and their associated linear equations.
To begin with, we recall isomonodromic deformation and its generalization.
Within the space of Fuchsian equations, we want to express submanifold whose points are equations with the same monodromy.
This submanifold is expressed by a differential equation. 
We call such an equation the isomonodromic deformation equation.
In non-Fuchsian case, we also require Stokes coefficients to be constant as well as the monodromy.
When we use the term isomonodromic deformation,
we mean this generalized isomonodromic deformation.
In this paper, we also use the term Painlev\'e-type equation in place of isomonodromic deformation equation.

As the coordinates for the submanifold, we can choose configurations of singularities of linear equations and coefficients of exponential parts of local solutions at irregular singularities  (\ref{eqn:exp_poly}),  see \cite{JMU}.
We call them deformation parameters.
Let us denote a deformation parameter by $t$.
We consider a linear equation
\[
\frac{dY}{dx}=A(x,t)Y
\]
whose coefficient is rational function in $x$.
Painlev\'e-type equation associated with this linear equation
is expressed by the compatibility condition of
\[
\frac{\partial Y}{\partial x}=A(x,t)Y, \hspace{5mm}
\frac{\partial Y}{\partial t}=B(x,t)Y,
\]
where $B$ is a matrix whose entries are rational functions in $x$.
This pair of linear equations is called the Lax pair.
The compatibility condition is given by
\[
\frac{\partial A(x,t)}{\partial t}-\frac{\partial B(x,t)}{\partial x}+[A(x,t),B(x,t)]=O.
\]  
From the compatibility condition, the entries of the matrix $B$ can be expressed by those of $A$.
This equation can be expressed as a Hamiltonian system.
As an example, we will derive a Hamiltonian of the sixth Painlev\'e equation from the linear equation in section \ref{sec:comment}, .
In this section, we give matrices $A, B$ and the Hamiltonian for each Painlev\'e-type equation.

We will divide the linear equations into four families, and describe
each of them in detail.

\subsection{Garnier system and degenerate Garnier systems}
In the first place, we will see a Lax pair
associated with the Garnier system of two variables.

The Garnier system is obtained from isomonodromic deformation of a
Fuchsian equation
of the second order with five regular singular points.
This Fuchsian equation is denoted by the singularity pattern $1+1+1+1+1$.
Local data that characterize linear equation are given by characteristic
exponents at each singular point, and this type of Fuchsian equation can be
reduced to a system with the following Riemann scheme by a suitable
gauge transformation:
\[
 \left(\begin{array}{ccccc}
  x=0 & x=1 & x=t_1 & x=t_2 & x=\infty \\
  0   &   0 &   0   &   0   & \theta^{\infty}_1 \\
\theta^0 & \theta^1 & \theta^{t_1} & \theta^{t_2} & \theta^{\infty}_2
       \end{array}\right) .
\]
The Fuchs relation is written as $\theta^{0}+\theta^{1}+\theta^{t_{1}}+\theta^{t_{2}}+\theta^{\infty}_{1}+\theta^{\infty}_{2}=0$.
We will see a parameterization of the Fuchsian system.
In this case the Lax pair is simply written by using
coefficients of the Fuchsian system.
\begin{equation}
\left\{
\begin{aligned}
\frac{\partial Y}{\partial x}&=\left(
\frac{A_0}{x}+\frac{A_1}{x-1}+\frac{A_{t_1}}{x-t_1}+\frac{A_{t_2}}{x-t_2}\right)Y, \\
\frac{\partial Y}{\partial t_1}&=-\frac{A_{t_1}}{x-t_1}Y, \qquad
\frac{\partial Y}{\partial t_2}=-\frac{A_{t_2}}{x-t_2}Y.
\end{aligned}
\right.
\end{equation}
Here $A_0$, $A_1$, and $A_t$ are given as follows:
\begin{align*}
A_{\xi}&=
\begin{pmatrix}
1 &  \\
 & u
\end{pmatrix}^{-1}
P^{-1}\hat{A_{\xi}}P
\begin{pmatrix}
1 & \\
 & u
\end{pmatrix} ,\quad
(\xi=0,1,t_1,t_2),\\
\hat{A}_0&=
\begin{pmatrix}
1 \\
0
\end{pmatrix}
\left(\theta^0,\ -1+\frac{q_1}{t_1}+\frac{q_2}{t_2}\right),\quad
\hat{A}_1=
\begin{pmatrix}
1 \\
p_1q_1+p_2q_2-\theta^\infty_2
\end{pmatrix}
\left(\theta^1+\theta^\infty_2-p_1q_1-p_2q_2,\ 1\right),\\
\hat{A}_{t_i}&=
\begin{pmatrix}
1 \\
t_ip_i \\
\end{pmatrix}
\left(\theta^{t_i}+p_iq_i,\ -\frac{q_i}{t_i}\right) ,
\qquad
\mbox{where}~~P=
 \begin{pmatrix}
 1 & 0 \\
 \frac{a}{\theta^{\infty}_1-\theta^{\infty}_2} & 1
 \end{pmatrix} .
\end{align*}
Here $a$ is the $(2,1)$-element of the matrix $\hat{A}_\infty:=-\hat{A}_0-\hat{A}_1-\hat{A}_{t_1}-\hat{A}_{t_2}$.

The compatibility conditions of these systems is expressed as a Hamiltonian system.
This is the partial differential system called Garnier system and
the system is given by the following Hamiltonians:
\begin{align}
 &t_i(t_i-1)H_{\mathrm{Gar}, t_i}^{1+1+1+1+1}
\left({\theta^0, \theta^1  \atop \theta^{t_1}, \theta^{t_2}, \theta^\infty_2}
;{t_1\atop t_2};{q_1,p_1 \atop q_2,p_2}\right)\\
&=t_i(t_i-1)H_{\rm VI}\left({\theta^\infty_2, \theta^1 \atop
\theta^{t_i}, \theta^0+\theta^{t_{i+1}}+1} ; t_i ; q_i,p_i\right)
+(2q_ip_i+q_{i+1}p_{i+1}-\theta^1-2\theta^\infty_2)q_1q_2p_{i+1}\nonumber\\
&\quad-\frac{1}{t_i-t_{i+1}}
\{ t_i(t_i-1)(p_iq_i+\theta^{t_i})p_iq_{i+1}-t_i(t_{i+1}-1)(2p_iq_i+\theta^{t_i})p_{i+1}q_{i+1}\nonumber\\
&\quad+t_{i+1}(t_i-1)({p_{i+1}}^2q_{i+1}+\theta^{t_{i+1}}(p_{i+1}-p_i))q_i\} ,
\qquad (i\in \mathbb{Z}/2\mathbb{Z}) .\nonumber
\end{align}
Furthermore the parameter $u$, which expresses gauge freedom, satisfy
the following equations:
\begin{align}
t_1(t_1-1)\frac{1}{u}\frac{\partial u}{\partial t_1}&=
q_1\{2p_1(t_1-q_1)+\theta^1+2\theta^\infty_2\}-2q_1p_2q_2+t_1\theta^{t_1},\\
t_2(t_2-1)\frac{1}{u}\frac{\partial u}{\partial t_2}&=
q_2\{2p_2(t_2-q_2)+\theta^1+2\theta^\infty_2\}-2q_2p_1q_1+t_2\theta^{t_2} .
\end{align}
In particular, the time evolution of $p_{i}, q_{i}\ (i=1,2)$ is independent
of $u$'s behavior.

In the second place, we will see the non-Fuchsian system obtained by a
confluence of two singular points from the above Fuchsian system.
This system is expressed by the singularity pattern $2+1+1+1$.

\subsubsection*{Singularity pattern: 2+1+1+1}
We have derived data below by calculating degenerations from $1+1+1+1+1$. 
To illustrate the procedure of degeneration,
we will see the degeneration of the Fuji-Suzuki system as an example, in section \ref{sec:comment}.

The Riemann scheme is given by
\[
\left(\begin{array}{cccc}
  x=0 & x=1 & x=t_2/t_1 & x=\infty \\
\begin{array}{c}0 \\ \theta^0 \end{array} &
\begin{array}{c}0 \\ \theta^1 \end{array} &
\begin{array}{c}0 \\ \theta^t \end{array} &
\overbrace{\begin{array}{cc}
     0    & \theta^\infty_1 \\
     -t_1 & \theta^\infty_2
        \end{array}}\ 
       \end{array}\right) ,
\]
and then the Fuchs-Hukuhara relation is written as
$\theta^0+\theta^1+\theta^t+\theta^\infty_1+\theta^\infty_2=0$.
The Lax pair is expressed as
\begin{equation}
\left\{
\begin{aligned}
\frac{\partial Y}{\partial x}&=\left(
\frac{A_0}{x}+\frac{A_1}{x-1}+\frac{A_t}{x-\frac{t_2}{t_1}}+A_\infty
 \right)Y ,\\
\frac{\partial Y}{\partial t_1}&=\left(
E_2x+B_1
+\frac{\frac{t_2}{{t_1}^2}A_t}{x-\frac{t_2}{t_1}}
\right)Y ,\qquad
\frac{\partial Y}{\partial
 t_2}=-\frac{\frac{1}{t_1}A_t}{x-\frac{t_2}{t_1}}Y .
\end{aligned}
\right.
\end{equation}
Here
\begin{align*}
A_{\xi}&=
\begin{pmatrix}
1 &  \\
 & u 
\end{pmatrix}^{-1}
\hat{A}_\xi
\begin{pmatrix}
1 & \\
 & u 
\end{pmatrix} ,
\quad(\xi=0,1,t) ,\\
\hat{A}_0&=
\begin{pmatrix}
1 \\
1
\end{pmatrix}
\left( \mu_1\lambda_1+\mu_2\lambda_2-\theta^1-\theta^t-\theta^\infty_1,
\ -\mu_1\lambda_1-\mu_2\lambda_2-\theta^\infty_2 \right) ,\\
\hat{A}_1&=
\begin{pmatrix}
1 \\
\lambda_1
\end{pmatrix}
\left( \theta^1-\mu_1\lambda_1,\ \mu_1 \right),\quad
\hat{A}_t=
\begin{pmatrix}
1 \\
\lambda_2 \\
\end{pmatrix}
\left( \theta^t-\mu_2\lambda_2,\ \mu_2 \right) ,\quad
A_\infty=
 \begin{pmatrix}
 0 &  \\
  & t_1
 \end{pmatrix} ,\\
A_\infty^{(0)}&=-(A_0+A_1+A_t) ,\quad
E_2=
 \begin{pmatrix}
 0 &  \\
  & 1
 \end{pmatrix} ,\quad
B_1=
\frac{1}{t_1}
\begin{pmatrix}
0 & (-A_\infty^{(0)})_{12} \\
(-A_\infty^{(0)})_{21} & 0
\end{pmatrix} .
\end{align*}
The Hamiltonians are given as 
\begin{align}\label{eq:Gar2111_1}
&\tilde{H}_{\mathrm{Gar}, t_1}^{2+1+1+1}
\left({\theta^\infty_2, \theta^1 \atop \theta^t, -\theta^0-1};{t_1 \atop t_2};
{\lambda_1,\mu_1\atop \lambda_2,\mu_2}\right)\\
&=\tilde{H}_{\rm V}\left({\theta^0+\theta^\infty_2,\ \theta^0+\theta^t+\theta^\infty_1 \atop
\theta^1};t_1;\lambda_1,\mu_1\right)
+\frac{\mu_2\lambda_2}{t_1}(1-\lambda_1)(\mu_1-(\mu_1\lambda_1-\theta^1))\nonumber\\
&\hspace{5mm}+\frac{1}{t_1-t_2}(\mu_1(\lambda_1-\lambda_2)-\theta^1)(\mu_2(\lambda_2-\lambda_1)-\theta^t),\nonumber\\
&\tilde{H}_{\mathrm{Gar}, t_2}^{2+1+1+1}
\left({\theta^\infty_2, \theta^1 \atop \theta^t, -\theta^0-1};{t_1 \atop t_2};
{\lambda_1,\mu_1\atop \lambda_2,\mu_2}\right)
\label{eq:Gar2111_2}
\\
&=\tilde{H}_{\rm V}\left({\theta^0+\theta^\infty_2,\ \theta^0+\theta^1+\theta^\infty_1 \atop
\theta^t}
;t_2;\lambda_2,\mu_2\right)+\frac{\mu_1\lambda_1}{t_2}(1-\lambda_2)(\mu_2-(\mu_2\lambda_2-\theta^t))\nonumber\\
&\hspace{5mm}+\frac{1}{t_2-t_1}(\mu_1(\lambda_1-\lambda_2)-\theta^1)(\mu_2(\lambda_2-\lambda_1)-\theta^t).\nonumber
\end{align}
Notice that $\tilde{H}_{\rm V}$ is different from $H_{\rm V}$ (see
 Remark \ref{rem:P_V}).

When we change canonical variables as $\lambda_{1}=1-\frac{1}{q_{1}}$, $\mu_{1}=q_{1}(p_{1}q_{1}-\theta^{1})$,
 $\lambda_{2}=1-\frac{1}{q_{2}}$, and $\mu_{2}=q_{2}(p_{2}q_{2}-\theta^{t})$,
the Hamiltonians are given as
\begin{align}
&H_{\mathrm{Gar}, t_1}^{2+1+1+1}
\left({\theta^\infty_2, \theta^1 \atop \theta^t, -\theta^0-1};{t_1 \atop t_2};
{q_1,p_1\atop q_2,p_2}\right)\\
&=H_{\rm V}\left({\theta^0+\theta^1+\theta^\infty_1,\ 2\theta^0+\theta^\infty_1+\theta^\infty_2 \atop
-\theta^0-\theta^\infty_1};t_1;q_1,p_1\right)
+\frac{p_1}{t_1}(p_2q_2(q_2-1)+\theta^t(q_1-q_2))\nonumber\\
&\hspace{5mm}-\frac{1}{t_1-t_2}(p_1(q_1-q_2)-\theta^1)(p_2(q_1-q_2)+\theta^t) ,\nonumber\\
&H_{\mathrm{Gar}, t_2}^{2+1+1+1}
\left({\theta^\infty_2, \theta^1 \atop \theta^t,-\theta^0-1};{t_1 \atop t_2};
{q_1,p_1\atop q_2,p_2}\right)\\
&=H_{\rm V}\left({\theta^0+\theta^t+\theta^\infty_1,\ 2\theta^0+\theta^\infty_1+\theta^\infty_2 \atop
-\theta^0-\theta^\infty_1};t_2;q_2,p_2\right)
+\frac{p_2}{t_2}(p_1q_1(q_1-1)+\theta^1(q_2-q_1))\nonumber\\
&\hspace{5mm}-\frac{1}{t_2-t_1}(p_1(q_1-q_2)-\theta^1)(p_2(q_1-q_2)+\theta^t)
 .\nonumber
\end{align}

The gauge parameter $u$ satisfy the equations:
\begin{align}
\frac{1}{u}\frac{\partial u}{\partial t_1}=\frac{1}{t_1}(p_1+\theta^\infty_1-\theta^\infty_2),\quad
\frac{1}{u}\frac{\partial u}{\partial t_2}=\frac{p_2}{t_2}.
\end{align}

\medskip
\noindent
{\bf Singularity pattern: 3+1+1}

The Riemann scheme is given as
\[
\left(
\begin{array}{ccc}
x=0  & x=t_2-t_1 & x=\infty	\\
    \begin{array}{c}
	0  \\ 
	\theta^0
	\end{array}
   &\begin{array}{c}
	0  \\ 
	\theta^1
	\end{array}
&\overbrace{
	\begin{array}{ccc}
	0 & 0 & \theta^{\infty}_1 \\
	-1& t_2 & \theta^{\infty}_2
	\end{array}}
\end{array}
\right),
\]
and the Fuchs-Hukuhara relation is written as
$\theta^0+\theta^1+\theta_1^\infty +\theta_2^\infty =0$.

The Lax pair is given as
\begin{equation}
\left\{
\begin{aligned}
\frac{\partial Y}{\partial x}&=\left(
\frac{A_0^{(0)}}{x}+\frac{A_1^{(0)}}{x-(t_2-t_1)}+A_\infty^{(-1)}+A_\infty^{(-2)}x
 \right)Y ,\\
\frac{\partial Y}{\partial t_1}&=\frac{A_1^{(0)}}{x-(t_2-t_1)}Y ,\qquad
\frac{\partial Y}{\partial t_2}=
\left( -\frac{A_1^{(0)}}{x-(t_2-t_1)}-A_\infty^{(-2)}x+B_1 \right)Y .
\end{aligned}
\right.
\end{equation}
Here
\begin{align*}
A_{\xi}^{(-k)}&=
\begin{pmatrix}
1 &  \\
 & u 
\end{pmatrix}^{-1}
\hat{A}_\xi
\begin{pmatrix}
1 & \\
 & u 
\end{pmatrix},\\
\hat{A}_0^{(0)}&=
\begin{pmatrix}
q_2 \\
1
\end{pmatrix}
\left( p_2,\ -p_2q_2+\theta^0 \right),\quad
\hat{A}_1^{(0)}=
\begin{pmatrix}
q_1 \\
1 \\
\end{pmatrix}
\left( p_1,\ -p_1q_1+\theta^1 \right) ,\quad
A_\infty^{(-2)}=
\begin{pmatrix}
0 &  \\
 & 1
\end{pmatrix},\\
\hat{A}_\infty^{(-1)}&=-
\begin{pmatrix}
0 & p_1q_1+p_2q_2+\theta^\infty_1 \\
1 & t_2
\end{pmatrix},\quad
B_1=
\begin{pmatrix}
0 & (-A_\infty^{(-1)})_{12} \\
(-A_\infty^{(-1)})_{21} & 0
\end{pmatrix}.
\end{align*}
The Hamiltonians are written as
\begin{align}
&H^{3+1+1}_{\mathrm{Gar},t_1}
\left({\theta^1, \theta^0 \atop \theta^\infty_1};{t_1 \atop t_2};
{q_1,p_1 \atop q_2,p_2}\right)\\
&=
H_{\rm{IV}}\left({\theta^1,\theta^\infty_1};t_1;q_1,p_1\right)
+p_2q_2p_1+\frac{1}{t_1-t_2}\{p_1(q_1-q_2)-\theta^1\}\{p_2(q_2-q_1)-\theta^0\},\nonumber\\
&H^{3+1+1}_{\mathrm{Gar},t_2}
\left({\theta^1, \theta^0 \atop \theta^\infty_1};{t_1 \atop t_2};
{q_1,p_1 \atop q_2,p_2}\right)\\
&=
H_{\rm{IV}}\left({\theta^0,\theta^\infty_1};t_2;q_2,p_2\right)
+p_1q_1p_2+\frac{1}{t_2-t_1}\{p_1(q_1-q_2)-\theta^1\}\{p_2(q_2-q_1)-\theta^0\}.\nonumber
\end{align}
The gauge parameter $u$ satisfies
\begin{equation}
\frac{1}{u}\frac{\partial u}{\partial t_1}=-p_1,\quad
\frac{1}{u}\frac{\partial u}{\partial t_2}=t_2-p_2.
\end{equation}

\medskip
\noindent
{\bf Singularity pattern: 2+2+1}

Rienmann scheme is given by
\[
\left(
\begin{array}{ccc}
x=0  & x=1 & x=\infty	\\
\overbrace{
	\begin{array}{cc}
	0 & 0\\
	\frac{t_2}{t_1} & \theta^{0}
	\end{array}}
	&\begin{array}{c}
	0\\ 
	\theta^{1}
	\end{array}
&\overbrace{
	\begin{array}{cc}
	0 & \theta^\infty_1 \\
	-t_1 & \theta^\infty_2
	\end{array}}
\end{array}
\right),
\]
and the Fuchs-Hukuhara relation is written as
$\theta^0+\theta^1+\theta^\infty_1+\theta^\infty_2=0$.

The Lax pair is given as
\begin{equation}
\left\{
\begin{aligned}
\frac{\partial Y}{\partial x}&=\left(
\frac{A_0^{(-1)}}{x^2}+\frac{A_0^{(0)}}{x}+\frac{A_1^{(0)}}{x-1}+A_\infty
 \right)Y ,\\
\frac{\partial Y}{\partial t_1}&=
\left(E_2x+B_1
+\frac{A_0^{(-1)}}{t_1x}
\right)Y ,\qquad
\frac{\partial Y}{\partial t_2}=-\frac{A_0^{(-1)}}{t_2x}Y .
\end{aligned}
\right.
\end{equation}
Here
\begin{align*}
A_{\xi}^{(-k)}&=
\begin{pmatrix}
1 &  \\
 & u 
\end{pmatrix}^{-1}
\hat{A}_\xi
\begin{pmatrix}
1 & \\
 & u 
\end{pmatrix},\\
\hat{A}_0^{(-1)}&=\frac{t_2}{t_1}
\begin{pmatrix}
1 \\
1
\end{pmatrix}
\left( 1-\mu_2,\ \mu_2 \right),\quad
\hat{A}_0^{(0)}=
\begin{pmatrix}
\mu_1\lambda_1-\theta^1-\theta^\infty_1 & -\mu_1\lambda_1-\mu_2\lambda_2-\theta^\infty_2 \\
\mu_1\lambda_1+(1-\mu_2)\lambda_2-\theta^1-\theta^\infty_1 & -\mu_1\lambda_1-\theta^\infty_2
\end{pmatrix},\\
\hat{A}_1^{(0)}&=
\begin{pmatrix}
1 \\
\lambda_1 \\
\end{pmatrix}
\left( -\mu_1\lambda_1+\theta^1,\ \mu_1 \right),\quad
A_\infty=
 \begin{pmatrix}
 0 &  \\
  & t_1
 \end{pmatrix},\quad
A_\infty^{(0)}=-(A_0^{(0)}+A_1^{(0)}),\\
E_2&=
 \begin{pmatrix}
 0 &  \\
  & 1
 \end{pmatrix},\quad
B_1=
\frac{1}{t_1}
\begin{pmatrix}
0 & (-A_\infty^{(0)})_{12} \\
(-A_\infty^{(0)})_{21} & 0
\end{pmatrix}.
\end{align*}

The Hamiltonians are expressed as
\begin{align}
&t_1\tilde{H}^{2+2+1}_{\mathrm{Gar},t_1}\left({\theta^0,\theta^1 \atop \theta^\infty_1,\theta^\infty_2};{t_1 \atop t_2};
{\lambda_1,\mu_1 \atop \lambda_2,\mu_2}\right)\\
&=
t_1\tilde{H}_{\rm{V}}\left({\theta^0+\theta^\infty_2, \theta^0+\theta^\infty_1 \atop \theta^1};t_1;\lambda_1,\mu_1\right)
+((\mu_1\lambda_1-\theta^1 )\lambda_1-\mu_1)\mu_2\lambda_2+\mu_1\lambda_2\nonumber\\
&\quad-\frac{t_2}{t_1}(\mu_1(\lambda_1-1)-\theta^1)(\mu_2(\lambda_1-1)+1),\nonumber\\
&t_2\tilde{H}^{2+2+1}_{\mathrm{Gar},t_{2}}\left({\theta^0,\theta^1 \atop \theta^\infty_1,\theta^\infty_2};{t_1 \atop t_2};
{\lambda_1,\mu_1 \atop \lambda_2,\mu_2}\right)\\
&=
t_2H_{\rm{III}}(D_6)
\left({\theta^\infty_2,-\theta^0};t_2;\lambda_2,\mu_2\right)
-\mu_1\lambda_1\lambda_2+\frac{t_2}{t_1}(\mu_1(\lambda_1-1)-\theta^1)(\mu_2(\lambda_1-1)+1).\nonumber
\end{align}

If we change the canonical variables as
\[
\lambda_1 = 1-\frac{1}{q_1},\quad \mu_1 = q_1(p_1q_1-\theta^1),
\]
then we obtain
\begin{align}
&t_1H^{2+2+1}_{\mathrm{Gar},t_1}\left({\theta^1, \theta^0 \atop -\theta^0-\theta^1-\theta^\infty_2};{t_1 \atop t_2};
{q_1,p_1 \atop q_2,p_2}\right)\\
&=
t_1H_{\rm{V}}
\left({-\theta^\infty_2,\theta^0-\theta^1
\atop \theta^1+\theta^\infty_2};t_1;q_1,p_1\right)
+q_1q_2(p_1q_1-\theta^1)+p_2q_2(\theta^1+p_1-2p_1q_1)
-\frac{t_2}{t_1}p_1(p_2-q_1),\nonumber\\
&t_2H^{2+2+1}_{\mathrm{Gar},t_{2}}\left({\theta^1, \theta^0 \atop -\theta^0-\theta^1-\theta^\infty_2};{t_1 \atop t_2};
{q_1,p_1 \atop q_2,p_2}\right)\\
&=
t_2H_{\rm{III}}(D_6)
\left({\theta^\infty_2,-\theta^0};t_2;q_2,p_2\right)
-(p_1q_1-\theta^1)q_2(q_1-1)+\frac{t_2}{t_1}p_1(p_2-q_1).\nonumber
\end{align}

The gauge $u$ satisfies
\begin{equation}
\frac{1}{u}\frac{\partial u}{\partial t_1}=\frac{p_1+\theta^\infty_1-\theta^\infty_2}{t_1},\quad
\frac{1}{u}\frac{\partial u}{\partial t_2}=-\frac{q_2}{t_2}.
\end{equation}

\medskip
\noindent
{\bf Singularity pattern: 3+2}

The Riemann scheme is given by
\[
\left(
\begin{array}{cc}
x=0  & x=\infty	\\
\overbrace{
    \begin{array}{cc}
	0 & 0  \\ 
	t_1 & \theta^0
	\end{array}}
&\overbrace{
	\begin{array}{ccc}
	0 & 0 & \theta^{\infty}_1 \\
	-1 & t_2 & \theta^{\infty}_2
	\end{array}}
\end{array}
\right) ,
\]
and then the Fuchs-Hukuhara relation is written as
$\theta^0+\theta^\infty_1+\theta^\infty_2=0$.

The Lax pair is expressed as
\begin{equation}
\left\{
\begin{aligned}
\frac{\partial Y}{\partial x}&=\left(
\frac{A_0^{(-1)}}{x^2}+\frac{A_0^{(0)}}{x}+A_\infty^{(-1)}+A_\infty^{(-2)}x
 \right)Y ,\\
\frac{\partial Y}{\partial t_1}&=-\frac{A_0^{(-1)}}{t_1x}Y ,\qquad
\frac{\partial Y}{\partial t_2}=(-A_\infty^{(-2)}x+B_1)Y ,
\end{aligned}
\right.
\end{equation}
where
\begin{align*}
A_{\xi}^{(-k)}&=
\begin{pmatrix}
1 &  \\
 & u 
\end{pmatrix}^{-1}
\hat{A}_\xi
\begin{pmatrix}
1 & \\
 & u 
\end{pmatrix},\\
\hat{A}_0^{(-1)}&=
\begin{pmatrix}
q_2 \\
1
\end{pmatrix}
\left( -q_1,\ q_1q_2+ t_1 \right),\quad
\hat{A}_0^{(0)}=
\begin{pmatrix}
-p_1q_1+p_2q_2 & -q_2(p_2q_2-\theta^0)+p_1(2q_1q_2+t_1)\\
p_2 & p_1q_1-p_2q_2+\theta^0 \\
\end{pmatrix},\\
A_\infty^{(-2)}&=
\begin{pmatrix}
0 &  \\
 & 1
\end{pmatrix},\quad
\hat{A}_\infty^{(-1)}=
\begin{pmatrix}
0 & p_1q_1-p_2q_2-\theta^\infty_1 \\
-1 & -t_2
\end{pmatrix},\quad
B_1=
\begin{pmatrix}
0 & (-A_\infty^{(-1)})_{12} \\
(-A_\infty^{(-1)})_{21} & 0
\end{pmatrix}.
\end{align*}
The Hamiltonians are given by
\begin{align}
H^{3+2}_{\mathrm{Gar},t_1}
\left({\theta^0, \theta^\infty_1};{t_1 \atop t_2};
{q_1,p_1 \atop q_2,p_2}\right)
&=
H_{\rm{III}}(D_6)\left({-\theta^{\infty}_1,\theta^0+1};t_1;q_1,p_1\right)
-p_1-\frac{q_1q_2}{t_1}(q_2-p_2+t_2)+p_1p_2-q_2,\\
H^{3+2}_{\mathrm{Gar},t_2}
\left({\theta^0, \theta^\infty_1};{t_1 \atop t_2};
{q_1,p_1 \atop q_2,p_2}\right)
&=
H_{\rm{IV}}\left({\theta^0,\theta^{\infty}_1};t_2;q_2,p_2\right)
-p_1q_1(p_2-2q_2-t_2)-q_1q_2+t_1p_1.
\end{align}
The gauge parameter $u$ satisfies
\begin{equation}
\frac{1}{u}\frac{\partial u}{\partial t_1}=-\frac{q_1}{t_1},\quad
\frac{1}{u}\frac{\partial u}{\partial t_2}=t_2-p_2.
\end{equation}

\medskip
\noindent
{\bf Singularity pattern: 4+1}

The Riemann scheme is given by
\[
\left(
\begin{array}{cc}
x=0 & x=\infty	\\
    \begin{array}{c}
	0  \\ 
	\theta^0
	\end{array}
&\overbrace{
	\begin{array}{cccc}
	0 & 0    & 0          & \theta^{\infty}_1 \\
	1 & 2t_2 & t_1+{t_2}^2  & \theta^{\infty}_2
	\end{array}}
\end{array}
\right) ,
\]
and then the Fuchs-Hukuhara relation is written as
$\theta^0 +\theta_1^\infty +\theta_2^\infty =0$.

The Lax pair is expressed as
\begin{equation}
\left\{
\begin{aligned}
\frac{\partial Y}{\partial x}&=\left(
\frac{A_0^{(0)}}{x}+A_\infty^{(-1)}+A_\infty^{(-2)}x+A_\infty^{(-3)}x^2
 \right)Y ,\\
\frac{\partial Y}{\partial t_1}&=(-E_2x+B_1)Y ,\qquad
\frac{\partial Y}{\partial t_2}=
\left( -E_2x^2+A_\infty^{(-2)}x+A_\infty^{(-1)}+T^{(\infty)}_1 \right)Y .
\end{aligned}
\right.
\end{equation}
Here
\begin{align*}
A_{\xi}^{(-k)}&=
\begin{pmatrix}
1 &  \\
 & u 
\end{pmatrix}^{-1}
\hat{A}_\xi^{(-k)}
\begin{pmatrix}
1 & \\
 & u 
\end{pmatrix},\qquad
\hat{A}_0^{(0)}=
\begin{pmatrix}
q_2 \\
1
\end{pmatrix}
\left( p_2,\ -p_2q_2+\theta^0 \right),\\
A_\infty^{(-3)}&=
\begin{pmatrix}
0 &  \\
 & -1
\end{pmatrix},\quad
\hat{A}_\infty^{(-2)}=
\begin{pmatrix}
0 & p_1 \\
1 & -2t_2
\end{pmatrix},\quad
\hat{A}_\infty^{(-1)}=
\begin{pmatrix}
-p_1 & p_1(q_1+t_2)-p_2q_2-\theta^\infty_1\\
-q_1+t_2 & p_1-t_1-{t_2}^2
\end{pmatrix},\\
B_1&=
\begin{pmatrix}
0 & \left(A_\infty^{(-2)}\right)_{12} \\
\left(A_\infty^{(-2)}\right)_{21} & 0
\end{pmatrix},\quad
T^{(\infty)}_1=
\begin{pmatrix}
0 & \\
 & t_1+{t_2}^2
\end{pmatrix}
,\quad
E_2=
\begin{pmatrix}
0 & \\
 & 1
\end{pmatrix}.
\end{align*}
The Hamiltonians are given by
\begin{align}
&H^{4+1}_{\mathrm{Gar},t_1}
\left({\theta^0, \theta^\infty_1};{t_1 \atop t_2};
{q_1,p_1 \atop q_2,p_2}\right)
=
H_{\rm{II}}\left({-\theta^\infty_1};t_1;q_1,p_1\right)
+p_2q_2(q_1-q_2+t_2)+p_1p_2+\theta^0 q_2,\\
&H^{4+1}_{\mathrm{Gar},t_2}
\left({\theta^0, \theta^\infty_1};{t_1 \atop t_2};
{q_1,p_1 \atop q_2,p_2}\right)\\
&=
-{p_2}^2q_2-t_2p_2{q_2}^2+{t_2}^2p_2q_2+\theta^0 t_2q_2-\theta^\infty_1 p_2
+p_1p_2(q_1-2q_2+t_2)+q_1q_2(p_2q_2-\theta^0)+\theta^0 p_1+t_1p_2q_2\nonumber.
\end{align}

The gauge parameter satisfies
\begin{equation}
\frac{1}{u}\frac{\partial u}{\partial t_1}=-q_1-t_2,\quad
\frac{1}{u}\frac{\partial u}{\partial t_2}=p_2-t_1-{t_2}^2.
\end{equation}

\medskip
\noindent
{\bf Singularity pattern: 5}

The Riemann scheme is given by
\[
\left(
\begin{array}{c}
x=\infty \\
\overbrace{
	\begin{array}{ccccc}
	0  & 0 & 0    & 0   & \theta^{\infty}_1 \\
	-1 & 0 & -2t_2 & -t_1 & \theta^{\infty}_2
	\end{array}}
\end{array}
\right) ,
\]
and Fuchs-Hukuhara relation is written as
$\theta^\infty_1+\theta^\infty_2=0$.

The Lax pair is expressed as
\begin{equation}
\left\{
\begin{aligned}
\frac{\partial Y}{\partial x}&=\left(
A_\infty^{(-4)}x^3+A_\infty^{(-3)}x^2+A_\infty^{(-2)}x+A_\infty^{(-1)}
 \right)Y ,\\
\frac{\partial Y}{\partial
 t_1}&=\left(A_\infty^{(-4)}x+A_\infty^{(-3)}\right)Y ,\qquad
\frac{\partial Y}{\partial
 t_2}=\left(A_\infty^{(-4)}x^2+A_\infty^{(-3)}x+A_\infty^{(-2)}+T_\infty^{(-2)}\right)Y ,
\end{aligned}
\right.
\end{equation}
where
\begin{align*}
A_\infty^{(-k)}&=
\begin{pmatrix}
1 &  \\
 & u 
\end{pmatrix}^{-1}
\hat{A}_\infty^{(-k)}
\begin{pmatrix}
1 & \\
 & u 
\end{pmatrix},\\
A_\infty^{(-4)}&=
\begin{pmatrix}
0 &  \\
 & 1
\end{pmatrix},\quad
\hat{A}_\infty^{(-3)}=
\begin{pmatrix}
0 & q_2 \\
-1 & 0
\end{pmatrix},\quad
\hat{A}_\infty^{(-2)}=
\begin{pmatrix}
-q_2 & -p_1 \\
-q_1 & q_2+2t_2 \\
\end{pmatrix},\\
\hat{A}_\infty^{(-1)}&=
\begin{pmatrix}
p_1-q_1q_2 & p_1q_1-(p_2-q_2-2t_2)q_2-\theta^\infty_1 \\
-p_2 & -p_1+q_1q_2+t_1
\end{pmatrix},\quad
T_\infty^{(-2)}=
\begin{pmatrix}
0 &  \\
 & -2t_2
\end{pmatrix}.
\end{align*}
The Hamiltonians are given by
\begin{align}
&H_{\mathrm{Gar},t_1}^5
\left({\theta^\infty_1};{t_1 \atop t_2};
{q_1,p_1 \atop q_2,p_2}\right)=
-q_{1}(q_1 p_1-\theta^{\infty}_{1})
+q_2(q_1(p_2+q_2)-2p_1+t_1)
+p_1(p_2-2t_2),\\
&H_{\mathrm{Gar},t_2}^5
\left({\theta^\infty_1};{t_1 \atop t_2};
{q_1,p_1 \atop q_2,p_2}\right)=
H_{\rm IV}\left(-1,\theta^\infty_1;2t_2;q_2,p_2\right)
+q_1q_2(q_1q_2-2p_1+t_1)
+p_1(p_1-p_2q_1-t_1).
\end{align}

The gauge parameter $u$ satisfies
\begin{equation}
\frac{1}{u}\frac{\partial u}{\partial t_1}=-q_1,\quad
\frac{1}{u}\frac{\partial u}{\partial t_2}=2t_2-p_2.
\end{equation}

\subsection{Fuji-Suzuki system}
The next family of Painlev\'e-type equations are the Fuji-Suzuki system and its degenerated systems.
The Fuji-Suzuki system is derived from the Fuchsian equation of type $21,21,111,111$.
They are given by the Lax pairs of $3\times 3$ matrices.

\subsubsection*{Singularity pattern: 1+1+1+1}
\noindent
\underline {Spectral type: 21,21,111,111}

The Riemann scheme is given by
\[
 \left(\begin{array}{cccc}
  x=0 & x=1 & x=t & x=\infty \\
  0 & 0 & 0 & \theta^{\infty}_1 \\
  \theta_1^0 & 0 & 0 & \theta^{\infty}_2 \\
  \theta_2^0 & \theta^1 & \theta^t & \theta^{\infty}_3
       \end{array}\right) ,
\]
and then the Fuchs relation is written as
$\theta_1^0 +\theta_2^0 +\theta^1 +\theta^t +\theta_1^\infty
+\theta_2^\infty +\theta_3^\infty =0$.

The Lax pair is expressed as
\begin{equation}
\left\{
\begin{aligned}
\frac{\partial Y}{\partial x}&=
\left(
\frac{A_0}{x}+\frac{A_1}{x-1}+\frac{A_t}{x-t}
\right)Y ,\\
\frac{\partial Y}{\partial t}&=-\frac{A_t}{x-t}Y .
\end{aligned}
\right.
\end{equation}
Here
\begin{align*}
A_{\xi}&=
U^{-1}P^{-1}\hat{A}_{\xi}PU,\quad
P=
\begin{pmatrix}
1 & 0 & 0 \\
\frac{a}{\theta^{\infty}_1-\theta^{\infty}_2} & 1 & 0 \\
\frac{1}{\theta^{\infty}_1-\theta^{\infty}_3}
\left(b+\frac{ac}{\theta^{\infty}_1-\theta^{\infty}_2}\right)&
\frac{c}{\theta^{\infty}_2-\theta^{\infty}_3} & 1
\end{pmatrix} ,\quad
U={\rm diag}(1,u,v),
\quad (\xi=0,1,t),\\
\hat{A}_0&=
\begin{pmatrix}
\theta^0_1 & \frac{q_1}{t}-1 & \frac{q_2}{t}-1 \\
0 & \theta^0_2 & p_1(q_2-q_1)+\theta^{\infty}_2+\theta^0_2 \\
0 & 0 & 0
\end{pmatrix},\quad
\hat{A}_t=
\begin{pmatrix}
1 \\
tp_1 \\
tp_2
\end{pmatrix}
\left(\theta^t +p_1 q_1+p_2 q_2,\ -\frac{q_1}{t},\ -\frac{q_2}{t}\right) ,\\
\hat{A}_1&=
\begin{pmatrix}
1 \\
p_1 q_1-\theta^{\infty}_2-\theta^0_2 \\
p_2 q_2-\theta^{\infty}_3
\end{pmatrix}
(-p_1 q_1-p_2 q_2-\theta^0_1-\theta^t-\theta^{\infty}_1,\ 1,\ 1),
\end{align*}
\[
\left\{
\begin{aligned}
a&=-tp_1(p_1q_1+p_2q_2+\theta^t)
+(p_1q_1-\theta^0_2-\theta^{\infty}_2)(p_1q_1+p_2q_2+\theta^0_1+\theta^t+\theta^{\infty}_1) ,\\
b&=-tp_2(p_1q_1+p_2q_2+\theta^t)
+(p_2q_2-\theta^{\infty}_3)(p_1q_1+p_2q_2+\theta^0_1+\theta^t+\theta^{\infty}_1) ,\\
c&=p_2(q_1-q_2)+\theta^{\infty}_3 .
\end{aligned}
\right.
\]

The Hamiltonian is given by
\begin{align}
&H_{\mathrm{FS}}^{A_5}\left({\theta^0_2+\theta^\infty_2,\theta^\infty_3,\theta^t
 \atop \theta^1,\theta^0_1,\theta^0_2};t;{q_1, p_1 \atop q_2, p_2}\right)\\
&=H_{\rm{VI}}\left({\theta^0_2+\theta^\infty_2,\theta^1+\theta^\infty_3
\atop \theta^t+\theta^\infty_3,\theta^0_1-\theta^0_2+1};t;q_1,p_1\right)
+H_{\rm{VI}}\left({\theta^\infty_3,\theta^0_2+\theta^1+\theta^\infty_2
\atop \theta^0_2+\theta^t+\theta^\infty_2,\theta^0_1-\theta^0_2-\theta^\infty_2+1};t;q_2,p_2\right)\nonumber\\
&\quad+\frac{1}{t(t-1)}(q_1-t)(q_2-1)\{(p_1q_1-\theta^0_2-\theta^\infty_2)p_2+p_1(p_2q_2-\theta^\infty_3)\}.\nonumber
\end{align}

The gauge parameters $u, v$ satisfy
\begin{align}
t(t-1)\frac{1}{u}\frac{du}{dt}&=
2p_1q_1(t-q_1)+p_2q_2(t-q_2)+(-\theta^0_1+\theta^0_2-\theta^t-\theta^\infty_1+\theta^\infty_2)q_1
+\theta^\infty_3 q_2+q_1p_2(1-q_2)+t\theta^t,\\
t(t-1)\frac{1}{v}\frac{dv}{dt}&=
q_1\{p_1(t-q_1)+\theta^0_2+\theta^\infty_2\}+q_2\{2p_2(t-q_2)-\theta^0_1-\theta^t-\theta^\infty_1+\theta^\infty_3\}
+p_1q_2(t-q_1)+t\theta^t.
\end{align}

\subsubsection*{Singularity pattern: 2+1+1}
\noindent
\underline {Spectral type: (2)(1),111,111}

The Riemann scheme is given by
\[
\left(\begin{array}{ccc}
  x=0 & x=1 & x=\infty \\
\begin{array}{c}0 \\ \theta^0_1 \\ \theta^0_2  \end{array}
& \overbrace{\begin{array}{cc}
     0 & 0 \\
     0 & 0 \\
     t & \theta^1
        \end{array}}\ 
& \begin{array}{c} \theta^{\infty}_1 \\ \theta^{\infty}_2  \\ \theta^{\infty}_3 \end{array} 
       \end{array}\right) ,
\]
and the Fuchs-Hukuhara relation is written as
$\theta_1^0 +\theta_2^0 +\theta^1 +\theta_1^\infty +\theta_2^\infty
+\theta_3^\infty =0$.

The Lax pair is expressed as
\begin{equation}
\left\{
\begin{aligned}
\frac{\partial Y}{\partial x}&=
\left(
\frac{A_1^{(-1)}}{(x-1)^2}+\frac{A_1^{(0)}}{x-1}+\frac{A_0^{(0)}}{x}
\right)Y ,\\
\frac{\partial Y}{\partial
 t}&=-\frac{1}{x-1}\left(\frac{A_1^{(-1)}}{t}\right)Y .
\end{aligned}
\right.
\end{equation}
Here
\begin{align*}
A_{\xi}^{(k)}&=
U^{-1}P^{-1}\hat{A}_{\xi}^{(k)}PU,\quad
P=
\begin{pmatrix}
1 & 0 & 0 \\
\frac{a}{\theta^{\infty}_1-\theta^{\infty}_2} & 1 & 0 \\
\frac{1}{\theta^{\infty}_1-\theta^{\infty}_3}
\left(b+\frac{ac}{\theta^{\infty}_1-\theta^{\infty}_2}\right)&
\frac{c}{\theta^{\infty}_2-\theta^{\infty}_3} & 1
\end{pmatrix},\quad U={\rm diag}(1,u,v) ,\\
\hat{A}_0^{(0)}&=
\begin{pmatrix}
\theta^0_1 & t(q_2-1) & t(q_1-1) \\
0 & \theta^0_2 & p_2(q_1-q_2)+\theta^0_2+\theta^{\infty}_2 \\
0 & 0 & 0
\end{pmatrix},\quad
\hat{A}_1^{(-1)}=
\begin{pmatrix}
1 \\
p_2/t \\
p_1/t
\end{pmatrix}
(p_1+p_2+t,-t,-t),\\
\hat{A}_1^{(0)}&=
\begin{pmatrix}
\theta^0_2+\theta^1+\theta^{\infty}_2+\theta^{\infty}_3 & t(1-q_2) & t(1-q_1) \\
-a & -\theta^0_2-\theta^{\infty}_2 & p_2(q_2-q_1)-\theta^0_2-\theta^{\infty}_2 \\
-b & p_1(q_1-q_2)-\theta^{\infty}_3 & -\theta^{\infty}_3
\end{pmatrix},
\end{align*}
\[
\left\{
\begin{aligned}
ta&=p_2(q_2-1)(p_1+p_2+t)-(\theta^0_2+\theta^{\infty}_2)(p_1+t)
-(2\theta^0_2+\theta^1+2\theta^{\infty}_2+\theta^{\infty}_3)p_2 ,\\
tb&=p_1(q_1-1)(p_1+p_2+t)-\theta^{\infty}_3(p_2+t)
-(\theta^0_2+\theta^1+\theta^{\infty}_2+2\theta^{\infty}_3)p_1 ,\\
c&=p_1(q_2-q_1)+\theta^{\infty}_3 .
\end{aligned}
\right.
\]

The Hamiltonian is given by
\begin{align}
&tH_{\mathrm{NY}}^{A_5}
 \left({\theta^\infty_1-\theta^\infty_3-1, \theta^\infty_3, -\theta^\infty_2 \atop
 \theta^0_2+\theta^\infty_2, \theta^0_1-\theta^0_2};t;
{q_1,p_1\atop q_2,p_2}\right)\\
&=tH_{\rm V}
 \left({\theta^\infty_1-1,
\theta^0_1-\theta^0_2+\theta^\infty_1-\theta^\infty_2-\theta^\infty_3-1 \atop
-\theta^\infty_1+\theta^\infty_3+1};t;q_1,p_1\right)\nonumber\\
&\quad+tH_{\rm V}\left({\theta^0_2+\theta^\infty_1-\theta^\infty_3-1,
\theta^0_1-\theta^0_2+\theta^\infty_1-\theta^\infty_2-\theta^\infty_3-1 \atop
-\theta^\infty_1+\theta^\infty_2+\theta^\infty_3+1};t;q_2,p_2\right)+2p_1p_2q_1(q_2-1).\nonumber
\end{align}
The gauge parameters $u,v$ satisfy
\begin{align}
t\frac{1}{u}\frac{du}{dt}&=p_1(1-2q_1)+(2p_2+t)(1-q_2)
-\theta^0_1+\theta^0_2-\theta^\infty_1+\theta^\infty_2+\theta^\infty_3,\\
t\frac{1}{v}\frac{dv}{dt}&=(2p_1+t)(1-q_1)+2p_2(1-q_2)
-\theta^0_1+\theta^0_2-\theta^\infty_1+\theta^\infty_2+\theta^\infty_3.
\end{align}

\medskip
\noindent
\underline{Spectral type: (11)(1),21,111}

The Riemann scheme is given by
\[
\left(
\begin{array}{ccc}
  x=0 & x=1 & x=\infty \\
\overbrace{\begin{array}{cc}
     0 & 0 \\
     0 & \theta^0_1 \\
     t & \theta^0_2
        	\end{array}}\ 
& \begin{array}{c}0 \\ 0 \\ \theta^1  \end{array}
& \begin{array}{c} \theta^{\infty}_1 \\ \theta^{\infty}_2  \\ \theta^{\infty}_3 \end{array} 
\end{array}
\right) ,
\]
and the Fuchs-Hukuhara relation is written as
$\theta_1^0 +\theta_2^0 +\theta^1 +\theta_1^\infty +\theta_2^\infty
+\theta_3^\infty =0$.

The Lax pair is expressed as
\begin{equation}
\left\{
\begin{aligned}
\frac{\partial Y}{\partial x}&=
\left(
\frac{A_0^{(-1)}}{x^2}+\frac{A_0^{(0)}}{x}+\frac{A_1^{(0)}}{x-1}
\right)Y ,\\
\frac{\partial Y}{\partial
 t}&=-\frac{1}{x}\left(\frac{A_0^{(-1)}}{t}\right)Y .
\end{aligned}
\right.
\end{equation}
Here
\begin{align*}
A_{\xi}^{(k)}&=
U^{-1}P^{-1}\hat{A}_{\xi}^{(k)}PU,\quad
P=
\begin{pmatrix}
1 & 0 & 0 \\
\frac{a}{\theta^{\infty}_1-\theta^{\infty}_2} & 1 & 0 \\
\frac{1}{\theta^{\infty}_1-\theta^{\infty}_3}
\left(b+\frac{ac}{\theta^{\infty}_1-\theta^{\infty}_2}\right)&
\frac{c}{\theta^{\infty}_2-\theta^{\infty}_3} & 1
\end{pmatrix},\quad U={\rm diag}(1,u,v),\\
\hat{A}_1^{(0)}&=
\begin{pmatrix}
1 \\
-p_1q_1 \\
-p_2q_2
\end{pmatrix}
(p_1q_1+p_2q_2+\theta^1,1,1) ,\quad
\hat{A}_0^{(-1)}=
\begin{pmatrix}
1 \\
0 \\
0
\end{pmatrix}
(t, -1/q_1, -1/q_2) ,\\
\hat{A}_0^{(0)}&=
\begin{pmatrix}
-p_1q_1-p_2q_2-\theta^1-\theta^\infty_1 & -1 & -1 \\
-tq_1(p_1q_1-\theta^0_1-\theta^\infty_2) & p_1q_1-\theta^\infty_2 & p_1q_1 \\
-tq_2(p_2q_2-\theta^\infty_3) & q_2(p_2q_2-\theta^\infty_3)/q_1 & p_2q_2-\theta^\infty_3
\end{pmatrix} ,
\end{align*}
\[
\left\{
\begin{aligned}
a&=tq_1(p_1q_1-\theta^0_1-\theta^\infty_2)
+p_1q_1(p_1q_1+p_2q_2+\theta^1) ,\\
b&=tq_2(p_2q_2-\theta^\infty_3)+p_2q_2(p_1q_1+p_2q_2+\theta^1) ,\\
c&=-q_2(p_2q_2-\theta^\infty_3)\left(\frac{1}{q_1}-\frac{1}{q_2}\right)+\theta^{\infty}_3 .
\end{aligned}
\right.
\]

The Hamiltonian is given by
\begin{align}
&tH_{\mathrm{FS}}^{A_4}
 \left({\theta^1, \theta^0_1, \theta^\infty_1 \atop \theta^\infty_2, \theta^\infty_3};t;
{q_1,p_1\atop q_2,p_2}\right)\\
&=tH_{\rm V}\left(
{\theta^0_1+\theta^1+\theta^\infty_2+\theta^\infty_3,
\theta^1+\theta^\infty_1-\theta^\infty_2-1 \atop -\theta^1-\theta^\infty_3}
;t;q_1,p_1\right)\nonumber\\
&\quad+tH_{\rm V}\left(
{\theta^1+\theta^\infty_3, \theta^1+\theta^\infty_1-\theta^\infty_3-1 \atop -\theta^1}
;t;q_2,p_2\right)
+p_1(q_2-1)\{p_2(q_1+q_2)-\theta^{\infty}_3\}.\nonumber
\end{align}

The gauge parameters $u,v$ satisfy
\begin{equation}
 tq_1\frac{1}{u}\frac{du}{dt}=q_2(p_2q_2-p_2-\theta^\infty_3)-q_1(2p_1+p_2+t)-\theta^1,\quad
 tq_2\frac{1}{v}\frac{dv}{dt}=-p_1(q_1+q_2)-(2p_2+t)q_2-\theta^1.
\end{equation}

\medskip
\noindent
\underline{Spectral type: (1)(1)(1),21,21}

The Riemann scheme is given by
\[
\left(
\begin{array}{ccc}
x=0 & x=1 & x=\infty \\
\begin{array}{c} 0 \\ 0 \\ \theta^0 \end{array}&
\begin{array}{c} 0 \\ 0 \\ \theta^1 \end{array}
&
\overbrace{\begin{array}{cc}
     0 & \theta^{\infty}_1 \\
     -t_1 & \theta^{\infty}_2 \\
     -t_2 & \theta^{\infty}_3
           \end{array}}
\end{array}
\right),
\]
and the Fuchs-Hukuhara relation is written as
$\theta^0 +\theta^1 +\theta_1^\infty +\theta_2^\infty
+\theta_3^\infty =0$.

The Lax pair is expressed as
\begin{equation}
\left\{
\begin{aligned}
\frac{\partial Y}{\partial x}&=
\left(A_\infty+\frac{A_0^{(0)}}{x}+\frac{A_1^{(0)}}{x-1}
\right)Y ,\\
\frac{\partial Y}{\partial t_1}&=(E_2x+B_1)Y,\qquad
\frac{\partial Y}{\partial t_2}=(E_3x+B_2)Y ,
\end{aligned}
\right.
\end{equation}
where
\begin{align*}
A_{\xi}^{(k)}&=
U^{-1}\hat{A}_{\xi}^{(k)}U,\qquad
B_i=
U^{-1}\hat{B}_iU,\quad
U={\rm diag}(1,u,v),\\
A_\infty&=
\begin{pmatrix}
0 &    &  \\
 & t_1 &  \\
 &    & t_2
\end{pmatrix} ,\qquad
\hat{A}_0^{(0)}=
\begin{pmatrix}
1 \\
\mu_1 \\
\mu_2
\end{pmatrix}(\mu_1\lambda_1+\mu_2\lambda_2+\theta^0,\ -\lambda_1,\ -\lambda_2),\\
\hat{A}_1^{(0)}&=
\begin{pmatrix}
1 \\
\mu_1\lambda_1-\theta^\infty_2 \\
\mu_2\lambda_2-\theta^\infty_3
\end{pmatrix}(-\mu_1\lambda_1-\mu_2\lambda_2+\theta^1+\theta^\infty_2+\theta^\infty_3,\
 1,\ 1), \quad
 E_{2}={\rm diag}(0,1,0), \quad  E_{3}={\rm diag}(0,0,1),
 \\
\hat{B}_1&=
\begin{pmatrix}
0 & \frac{(\hat{A}_0^{(0)}+\hat{A}_1^{(0)})_{12}}{t_1} & 0 \\
\frac{(\hat{A}_0^{(0)}+\hat{A}_1^{(0)})_{21}}{t_1} & 0 & \frac{(\hat{A}_0^{(0)}+\hat{A}_1^{(0)})_{23}}{t_1-t_2}\\
0 & \frac{(\hat{A}_0^{(0)}+\hat{A}_1^{(0)})_{32}}{t_1-t_2} & 0
\end{pmatrix},\quad
\hat{B}_2=
\begin{pmatrix}
0 & 0 & \frac{(\hat{A}_0^{(0)}+\hat{A}_1^{(0)})_{13}}{t_2}\\
0 & 0 & \frac{(\hat{A}_0^{(0)}+\hat{A}_1^{(0)})_{23}}{t_2-t_1}\\
\frac{(\hat{A}_0^{(0)}+\hat{A}_1^{(0)})_{31}}{t_2} & \frac{(\hat{A}_0^{(0)}+\hat{A}_1^{(0)})_{32}}{t_2-t_1} & 0
\end{pmatrix}.
\end{align*}

The Hamiltonian is given by
\begin{align}
&t_1\tilde {H}^{2+1+1+1}_{\mathrm{Gar},t_{1}}
\left({\theta^0, \theta^1 \atop \theta^\infty_2,\theta^\infty_3};{t_1 \atop t_2};
{\lambda_1,\mu_1 \atop \lambda_2,\mu_2}\right)
\label{eq:Gar2111_3}
\\
&=
t_1\tilde{H}_{\rm V}\left({\theta^0, \theta^1+\theta^\infty_3 \atop \theta^\infty_2};t_1;\lambda_1,\mu_1\right)
+(1-\lambda_1)\mu_2\lambda_2(\mu_1-\mu_1\lambda_1+\theta^\infty_2)\nonumber\\
&\quad+\frac{t_1}{t_1-t_2}(\mu_1(\lambda_1-\lambda_2)-\theta^\infty_2)(\mu_2(\lambda_2-\lambda_1)-\theta^\infty_3),\nonumber\\
&t_2\tilde{H}^{2+1+1+1}_{\mathrm{Gar},t_{2}}
\left({\theta^0, \theta^1 \atop \theta^\infty_1,\theta^\infty_2,\theta^\infty_3};{t_1 \atop t_2};
{\lambda_1,\mu_1 \atop \lambda_2,\mu_2}\right)
\label{eq:Gar2111_4}
\\
&=
t_2\tilde{H}_{\rm V}\left({\theta^0, \theta^1+\theta^\infty_2 \atop \theta^\infty_3};t_2;\lambda_2,\mu_2\right)
+(1-\lambda_2)\mu_1\lambda_1(\mu_2-\mu_2\lambda_2+\theta^\infty_3)\nonumber\\
&\quad+\frac{t_2}{t_2-t_1}(\mu_1(\lambda_1-\lambda_2)-\theta^\infty_2)(\mu_2(\lambda_2-\lambda_1)-\theta^\infty_3).\nonumber
\end{align}

When we change the canonical variables as
\[
\lambda_1 = 1-\frac{1}{q_1},\ \mu_1 = q_1(p_1q_1-\theta^\infty_2),\ 
\lambda_2 = 1-\frac{1}{q_2},\ \mu_2 = q_2(p_2q_2-\theta^\infty_3) ,
\]
then we obtain
\begin{align}
H_{\mathrm{Gar},\ t_1}^{2+1+1+1}
\left({\theta^0+\theta^\infty_1,\theta^\infty_2 \atop \theta^\infty_3,\theta^\infty_1-1};
{t_1 \atop t_2};{q_1,p_1\atop q_2,p_2}\right)&=
H_{\rm{V}}\left({\theta^1+\theta^\infty_2,\theta^0+\theta^1 \atop -\theta^1};t_1;q_1,p_1\right)\\
&\quad+\frac{p_1}{t_1}[\theta^\infty_3(q_1-q_2)+p_2q_2(q_2-1)]\nonumber\\
&\quad+\frac{1}{t_1-t_2}\left((q_1-q_2)p_1-\theta^\infty_2
 \right)\left((q_2-q_1)p_2-\theta^\infty_3\right) ,\nonumber\\
H_{\mathrm{Gar},\ t_2}^{2+1+1+1}
\left({\theta^0+\theta^\infty_1,\theta^\infty_2 \atop \theta^\infty_3,\theta^\infty_1-1};
{t_1 \atop t_2};{q_1,p_1\atop q_2,p_2}\right)&=
H_{\rm{V}}\left({\theta^1+\theta^\infty_3,\theta^0+\theta^1 \atop -\theta^1};t_2;q_2,p_2\right)\\
&\quad+\frac{p_2}{t_2}[\theta^\infty_2(q_2-q_1)+p_1q_1(q_1-1)]\nonumber\\
&\quad+\frac{1}{t_2-t_1}\left((q_2-q_1)p_2-\theta^\infty_3
 \right)\left((q_1-q_2)p_1-\theta^\infty_2\right) .\nonumber
\end{align}

The gauge parameters  $u,v$ satisfy
\begin{align}
\frac{t_1q_1}{u}\frac{\partial u}{\partial t_1}&=
2p_1q_1(q_1-1)-(2\theta^\infty_2+t_1)q_1-p_2q_2-\theta^1
+\frac{p_2q_2(t_1q_1-t_2q_2)+\theta^\infty_3 t_2 q_2}{t_{1}-t_{2}},\\
\frac{t_1q_2}{v}\frac{\partial v}{\partial t_1}&=
\frac{-p_1q_1(t_1q_1+(t_2-2t_1)q_2)+\theta^\infty_2 t_1q_1}{t_{1}-t_{2}}-(p_1+\theta^\infty_2)q_2,\\
\frac{t_2q_1}{u}\frac{\partial u}{\partial t_2}&=
\frac{p_2q_2(t_2q_2+(t_1-2t_2)q_1)-\theta^\infty_3 t_2q_2}{t_{1}-t_{2}}-(p_2+\theta^\infty_3)q_1,\\
\frac{t_2q_2}{v}\frac{\partial v}{\partial t_2}&=
2p_2q_2(q_2-1)-(2\theta^\infty_3+t_2)q_2-p_1q_1-\theta^{1}
-\frac{p_1q_1(t_2q_2-t_1q_1)+\theta^\infty_2 t_1 q_1}{t_{1}-t_{2}}.
\end{align}

\subsubsection*{Singularity pattern: 3+1}
\noindent
\underline{Spectral type: ((11))((1)),111}

The Riemann scheme is given as
\[
\left(
\begin{array}{cc}
  x=0 & x=\infty \\
\overbrace{\begin{array}{ccc}
   0 & 0  & 0 \\
   0 & 0  & \theta^0_1 \\
   1 & -t & \theta^0_2
        	\end{array}}
& \begin{array}{c} \theta^{\infty}_1 \\ \theta^{\infty}_2  \\ \theta^{\infty}_3 \end{array}
\end{array}
\right),
\]
and then the Fuchs-Hukuhara relation is written as
$\theta_1^0 +\theta_2^0 +\theta_1^\infty +\theta_2^\infty
+\theta_3^\infty =0$.

The Lax pair is expressed as
\begin{equation}
\left\{
\begin{aligned}
\frac{\partial Y}{\partial x}&=
\left(
\frac{A_0^{(-2)}}{x^3}+\frac{A_0^{(-1)}}{x^2}+\frac{A_0^{(0)}}{x}
\right)Y ,\\
\frac{\partial Y}{\partial t}&=\frac{A_0^{(-2)}}{x}Y ,
\end{aligned}
\right.
\end{equation}
where
\begin{align*}
A_{\xi}^{(k)}&=
U^{-1}P^{-1}\hat{A}_{\xi}^{(k)}PU,\quad
P=
\begin{pmatrix}
1 & 0 & 0 \\
\frac{a}{\theta^{\infty}_1-\theta^{\infty}_2} & 1 & 0 \\
\frac{1}{\theta^{\infty}_1-\theta^{\infty}_3}
\left(b+\frac{ac}{\theta^{\infty}_1-\theta^{\infty}_2}\right)&
\frac{c}{\theta^{\infty}_2-\theta^{\infty}_3} & 1
\end{pmatrix},\quad U={\rm diag}(1,u,v),\\
\hat{A}_0^{(-2)}&=
\begin{pmatrix}
1  \\
0  \\
0
\end{pmatrix}
(1,1,1),\quad
\hat{A}_0^{(-1)}=
\begin{pmatrix}
p_1+p_2-t & q_2  & q_1 \\
-p_2      & -p_2 & -p_2 \\
-p_1      & -p_1 & -p_1
\end{pmatrix},\quad
\hat{A}_0^{(0)}=
\begin{pmatrix}
-\theta^\infty_1 & 0 & 0 \\
-a & -\theta^\infty_2 & 0 \\
-b & -c & -\theta^\infty_3
\end{pmatrix} ,
\end{align*}
\[
a=p_2(p_2-q_2-t)+p_1p_2+\theta^0_1+\theta^\infty_2 ,\qquad
b=p_1(p_1-q_1-t)+p_1p_2+\theta^\infty_3 ,\qquad
c=p_1(q_2-q_1)+\theta^{\infty}_3 .
\]

The Hamiltonian is given by
\begin{align}
&H_{\mathrm{NY}}^{A_4}
 \left({ \theta^\infty_1-\theta^\infty_3-1, \theta^\infty_3 \atop 
-\theta^\infty_2,\theta^0_1+\theta^\infty_2};t;
  {q_1,p_1\atop q_2,p_2}\right)\\
&=
H_{\rm IV}\left({\theta^\infty_3,\theta^\infty_1-\theta^\infty_3-1}
 ;t;q_1,p_1\right)
+H_{\rm IV}\left({\theta^0_1+\theta^\infty_2,\theta^\infty_1-\theta^\infty_2-\theta^\infty_3-1}
 ;t;q_2,p_2\right)
+2p_1q_1p_2.\nonumber
\end{align}

The gauge parameters $u,v$ satisfy
\begin{align}
\frac{1}{u}\frac{du}{dt}=-p_1-2p_2+q_2+t,\quad
\frac{1}{v}\frac{dv}{dt}=q_1-2p_1-2p_2+t.
\end{align}

\medskip
\noindent
\underline{Spectral type: ((1)(1))((1)),21}

The Riemann scheme is given by
\[
\left(
\begin{array}{cc}
  x=0 & x=\infty \\
\begin{array}{c} 0 \\ 0 \\ \theta^0 \end{array}
&\overbrace{\begin{array}{ccc}
    0 & 0   & \theta^\infty_1 \\
   -1 & t_1 & \theta^\infty_2 \\
   -1 & t_2 & \theta^\infty_3
        	\end{array}}
\end{array}
\right) ,
\]
and the Fuchs-Hukuhara relation is written as
$\theta^0 +\theta_1^\infty +\theta_2^\infty +\theta_3^\infty =0$.

The Lax pair is expressed as
\begin{equation}
\left\{
\begin{aligned}
\frac{\partial Y}{\partial x}&=
\left(
A_\infty^{(-2)}x+A_\infty^{(-1)}+\frac{A_0^{(0)}}{x}
\right)Y ,\\
\frac{\partial Y}{\partial t_1}&=(-E_2x+B_1)Y ,\qquad
\frac{\partial Y}{\partial t_2}=(-E_3x+B_2)Y ,
\end{aligned}
\right.
\label{eq:LaxGar3+1+1}
\end{equation}
where
\begin{align*}
A_{\xi}^{(k)}&=
U^{-1}\hat{A}_{\xi}^{(k)}U,\quad
B_i=
U^{-1}\hat{B}_iU,\quad
U={\rm diag}(1,u,v),\\
\hat{A}_\infty^{(-2)}&=
\begin{pmatrix}
0 &  &  \\
 & 1 &  \\
 &  & 1
\end{pmatrix},\quad
\hat{A}_\infty^{(-1)}=
\begin{pmatrix}
0                       & -1   & -1   \\
-p_1q_1+\theta^\infty_2 & -t_1 &  0   \\
-p_2q_2+\theta^\infty_3 &  0   & -t_2
\end{pmatrix},\\
\hat{A}_0^{(0)}&=
\begin{pmatrix}
1 \\
p_1 \\
p_2
\end{pmatrix}
(p_1q_1+p_2q_2+\theta^0,\ -q_1,\ -q_2),\quad
E_2={\rm diag}(0,1,0),\quad E_3={\rm diag}(0,0,1),\\
\hat{B}_1&=
\begin{pmatrix}
0 & 1 & 0 \\
p_1q_1-\theta^\infty_2 & 0 & \frac{p_1(q_1-q_2)-\theta^\infty_2}{t_1-t_2} \\
0 & \frac{p_2(q_2-q_1)-\theta^\infty_3}{t_1-t_2} & 0
\end{pmatrix},\quad
\hat{B}_2=
\begin{pmatrix}
0 & 0 & 1 \\
0 & 0 & \frac{p_1(q_1-q_2)-\theta^\infty_2}{t_2-t_1} \\
p_2q_2-\theta^\infty_3 & \frac{p_2(q_2-q_1)-\theta^\infty_3}{t_2-t_1} & 0
\end{pmatrix}.
\end{align*}

The Hamiltonians are given by
\begin{align}
&H^{3+1+1}_{\mathrm{Gar},t_1}
\left({\theta^\infty_2,\theta^\infty_3 \atop \theta^0};{t_1 \atop t_2};
{q_1,p_1 \atop q_2,p_2}\right)\\
&=
H_{\rm{IV}}\left({\theta^\infty_2,\theta^0};t_1;q_1,p_1\right)
+p_2q_2p_1+\frac{1}{t_1-t_2}\{p_1(q_1-q_2)-\theta^\infty_2\}\{p_2(q_2-q_1)
-\theta^\infty_3\},\nonumber\\
&H^{3+1+1}_{\mathrm{Gar},t_2}
\left({\theta^\infty_2,\theta^\infty_3 \atop \theta^0};{t_1 \atop t_2};
{q_1,p_1 \atop q_2,p_2}\right)\\
&=
H_{\rm{IV}}\left({\theta^\infty_3,\theta^0};t_2;q_2,p_2\right)
+p_1q_1p_2+\frac{1}{t_2-t_1}\{p_1(q_1-q_2)-\theta^\infty_2\}\{p_2(q_2-q_1)
-\theta^\infty_3\}.\nonumber
\end{align}

The gauge parameters $u, v$ satisfy
\begin{align}
(t_1-t_2)\frac{1}{u}\frac{\partial u}{\partial t_1}&=
p_2(q_1-q_2)+(t_1-t_2)(q_1+t_1)+\theta^\infty_3,\quad
(t_1-t_2)\frac{1}{v}\frac{\partial v}{\partial t_1}=
p_1(q_2-q_1)+\theta^\infty_2,\\
(t_2-t_1)\frac{1}{u}\frac{\partial u}{\partial t_2}&=
p_2(q_1-q_2)+\theta^\infty_3,\quad
(t_2-t_1)\frac{1}{v}\frac{\partial v}{\partial t_2}=
p_1(q_2-q_1)+(t_2-t_1)(q_2+t_2)+\theta^\infty_2.
\end{align}

\subsubsection*{Singularity pattern: 2+2}
\noindent
\underline{Spectral type: (11)(1),(11)(1)}

The Riemann scheme is given as
\[
\left(
\begin{array}{cc}
  x=0 & x=\infty \\
\overbrace{\begin{array}{cc}
     0 & 0 \\
     0 & \theta^0_1 \\
     t & \theta^0_2
           \end{array}}
&
\overbrace{\begin{array}{cc}
     0 & \theta^{\infty}_1 \\
     0 & \theta^{\infty}_2 \\
     1 & \theta^{\infty}_3
           \end{array}}
\end{array}
\right) ,
\]
and the Fuchs-Hukuhara relation is written as
$\theta_1^0 +\theta_2^0 +\theta_1^\infty +\theta_2^\infty
+\theta_3^\infty =0$.

The Lax pair is expressed as
\begin{equation}
\left\{
\begin{aligned}
\frac{\partial Y}{\partial x}&=
\left(
\frac{A_0^{(-1)}}{x^2}+\frac{A_0^{(0)}}{x}+A_\infty
\right)Y ,\\
\frac{\partial Y}{\partial t}&=-\frac{A_0^{(-1)}}{tx}Y ,
\end{aligned}
\right.
\end{equation}
where
\begin{align*}
A_{\xi}^{(k)}&=
U^{-1}P^{-1}\hat{A}_{\xi}^{(k)}PU,\quad
P=
\begin{pmatrix}
1 & 0 & 0 \\
-p_1q_1 & 1 & 0 \\
-p_2q_2 &
\frac{p_2q_2(q_1-q_2)+\theta^\infty_1q_2}{(\theta^\infty_2-\theta^\infty_1)q_1}
& 1
\end{pmatrix},\quad U={\rm diag}(1,u,v),\\
\hat{A}_\infty&=
\begin{pmatrix}
-1        \\
p_1q_1  \\
p_2q_2 
\end{pmatrix} (1,0,0),\quad
\hat{A}_0^{(0)}=
\begin{pmatrix}
-p_1q_1-p_2q_2-\theta^\infty_3 & -1 & -1 \\
q_1(p_1q_1-\theta^0_1-\theta^\infty_2) & p_1q_1-\theta^\infty_2 & p_1q_1 \\
q_2(p_2q_2-\theta^\infty_1) & q_2(p_2q_2-\theta^\infty_1)/q_1 & p_2q_2-\theta^\infty_1
\end{pmatrix},\\
\hat{A}_0^{(-1)}&=\begin{pmatrix}
1  \\
0  \\
0 
\end{pmatrix}(t, t/q_1, t/q_2).
\end{align*}

The Hamiltonian is given by
\begin{align}
&tH_{\mathrm{FS}}^{A_3}
\left({ -\theta^0_1, -\theta^\infty_3\atop -\theta^\infty_2, -\theta^\infty_1};t;
  {q_1,p_1 \atop q_2,p_2}\right)\\
&=tH_{\rm III}(D_6)\left({-\theta^0_1-\theta^\infty_2, \theta^\infty_3-\theta^\infty_2}
;t;q_1,p_1\right)
+tH_{\rm III}(D_6)\left({-\theta^\infty_1, \theta^\infty_3-\theta^\infty_1}
;t;q_2,p_2\right) \nonumber\\
&\quad+p_1q_2\{ p_2(q_1+q_2)-\theta^\infty_1)\}.\nonumber
\end{align}

The gauge parameters $u,v$ satisfy
\begin{equation}
tq_1\frac{1}{u}\frac{du}{dt}=q_2(p_2q_2-\theta^\infty_1)+t,\qquad
\frac{1}{v}\frac{dv}{dt}=\frac{1}{q_2}.
\end{equation}

\medskip
\noindent
\underline{Spectral type: (2)(1),(1)(1)(1)}

The Riemann scheme is given by
\[
\left(
\begin{array}{cc}
  x=0 & x=\infty \\
\overbrace{\begin{array}{cc}
     0    & 0 \\
     0    & 0 \\
     1 & \theta^0
           \end{array}}
&
\overbrace{\begin{array}{cc}
     0 & \theta^{\infty}_1 \\
     -t_1 & \theta^{\infty}_2 \\
     -t_2 & \theta^{\infty}_3
           \end{array}}
\end{array}
\right) ,
\]
and the Fuchs-Hukuhara relation is written as
$\theta^0 +\theta_1^\infty +\theta_2^\infty
+\theta_3^\infty =0$.

The Lax pair is expressed as
\begin{equation}
\left\{
\begin{aligned}
\frac{\partial Y}{\partial x}&=
\left(
\frac{A_0^{(-1)}}{x^2}+\frac{A_0^{(0)}}{x}+A_\infty
\right)Y ,\\
\frac{\partial Y}{\partial t_1}&=(E_2x+B_1)Y ,\qquad
\frac{\partial Y}{\partial t_2}=(E_3x+B_2)Y ,
\end{aligned}
\right.
\end{equation}
where
\begin{align*}
A_{\xi}^{(k)}&=
U^{-1}\hat{A}_{\xi}^{(k)}U,\quad
B_i=
U^{-1}\hat{B}_iU,\quad
U={\rm diag}(1,u,v),\\
A_\infty&=
\begin{pmatrix}
0 &    &  \\
 & t_1 &  \\
 &    & t_2
\end{pmatrix},\quad
\hat{A}_0^{(-1)}=
\begin{pmatrix}
1 \\
p_1 \\
p_2
\end{pmatrix}(-p_1-p_2+1,\ 1,\ 1),\\
\hat{A}_0^{(0)}&=
\begin{pmatrix}
-\theta^\infty_1 & -q_1 & -q_2 \\
-p_1q_1(p_1+p_2-1)+\theta^\infty_2(p_2-1)+(\theta^\infty_2-\theta^\infty_1)p_1 & -\theta^\infty_2 & p_1(q_1-q_2)-\theta^\infty_2 \\
-p_2q_2(p_1+p_2-1)+\theta^\infty_3(p_1-1)+(\theta^\infty_3-\theta^\infty_1)p_2 & p_2(q_2-q_1)-\theta^\infty_3 & -\theta^\infty_3
\end{pmatrix}, \\
\hat{B}_1&=
\begin{pmatrix}
0 & \frac{(\hat{A}_0^{(0)})_{12}}{t_1} & 0 \\
\frac{(\hat{A}_0^{(0)})_{21}}{t_1} & 0 & \frac{(\hat{A}_0^{(0)})_{23}}{t_1-t_2}\\
0 & \frac{(\hat{A}_0^{(0)})_{32}}{t_1-t_2} & 0
\end{pmatrix},\quad
\hat{B}_2=
\begin{pmatrix}
0 & 0 & \frac{(\hat{A}_0^{(0)})_{13}}{t_2}\\
0 & 0 & \frac{(\hat{A}_0^{(0)})_{23}}{t_2-t_1}\\
\frac{(\hat{A}_0^{(0)})_{31}}{t_2} & \frac{(\hat{A}_0^{(0)})_{32}}{t_2-t_1} & 0
\end{pmatrix}.
\end{align*}

The Hamiltonians are given by
\begin{align}
t_1H_{\mathrm{Gar},t_1}^{\frac{3}{2}+1+1+1}
\left({\theta^\infty_2,\theta^\infty_3 \atop \theta^\infty_1};{t_1 \atop t_2};
{q_1,p_1 \atop q_2,p_2}\right)
&=
t_1H_{\rm III}(D_6)\left({-\theta^\infty_2,\theta^\infty_1-\theta^\infty_2};t_1;q_1,p_1\right)\\
&+q_1(q_1p_1p_2-\theta^\infty_2 p_2)
+\frac{t_1}{t_1-t_2}(p_1(q_1-q_2)-\theta^\infty_2)(p_2(q_2-q_1)-\theta^\infty_3),\nonumber\\
t_2H_{\mathrm{Gar},t_2}^{\frac{3}{2}+1+1+1}
\left({\theta^\infty_2,\theta^\infty_3 \atop \theta^\infty_1};{t_1 \atop t_2};
{q_1,p_1 \atop q_2,p_2}\right)
&=
t_2H_{\rm III}(D_6)\left({-\theta^\infty_3,\theta^\infty_1-\theta^\infty_3};t_2;q_2,p_2\right)\\
&+q_2(q_2p_1p_2-\theta^\infty_3 p_1)
+\frac{t_2}{t_2-t_1}(p_1(q_1-q_2)-\theta^\infty_2)(p_2(q_2-q_1)-\theta^\infty_3).\nonumber
\end{align}

The gauge parameters $u,v$ satisfy
\begin{align}
t_1(t_1-t_2)\frac{1}{u}\frac{\partial u}{\partial t_1}&=
(t_1-t_2)(1-2p_1)q_1+p_2(t_2q_1-t_1q_2)+\theta^\infty_3t_1,\\
t_1(t_1-t_2)\frac{1}{v}\frac{\partial v}{\partial t_1}&=
t_1(\theta^\infty_2-2p_1q_1)+p_1(t_2q_1+t_1q_2),\\
t_2(t_2-t_1)\frac{1}{u}\frac{\partial u}{\partial t_2}&=
t_2(\theta^\infty_3-2p_2q_2)+p_2(t_2q_1+t_1q_2),\\
t_2(t_2-t_1)\frac{1}{v}\frac{\partial v}{\partial t_2}&=
(t_2-t_1)(1-2p_2)q_2+p_1(t_1q_2-t_2q_1)+\theta^\infty_2t_2.
\end{align}

\subsubsection*{Singularity pattern: 4}
\noindent
\underline{Spectral type: (((1)(1)))(((1)))}

The Riemann scheme is given by
\[
\left(
\begin{array}{c}
 x=\infty \\
\overbrace{\begin{array}{cccc}
   0 & 0 & 0   & \theta^\infty_1\\
   1 & 0 & t_1 & \theta^\infty_2\\
   1 & 0 & t_2 & \theta^\infty_3
        	\end{array}}
\end{array}
\right) ,
\]
and the Fuchs-Hukuhara relation is written as
$\theta_1^\infty +\theta_2^\infty +\theta_3^\infty =0$.

The Lax pair is expressed as
\begin{equation}
\left\{
\begin{aligned}
\frac{\partial Y}{\partial x}&=
\left(
A_\infty^{(-3)}x^2+A_\infty^{(-2)}x+A_\infty^{(-1)}
\right)Y ,\\
\frac{\partial Y}{\partial t_1}&=(-E_2x+B_1)Y ,\qquad
\frac{\partial Y}{\partial t_2}=(-E_3x+B_2)Y ,
\end{aligned}
\right.
\end{equation}
where
\begin{align*}
A_\infty^{(k)}&=
U^{-1}\hat{A}_\infty^{(k)}U,\quad
B_i=
U^{-1}\hat{B}_iU,\quad
U={\rm diag}(1,u,v),\quad (i=1,2),\\
\hat{A}_\infty^{(-3)}&=
\begin{pmatrix}
0 &   &  \\
 & -1 &  \\
 &   & -1
\end{pmatrix},\quad
\hat{A}_\infty^{(-2)}=
\begin{pmatrix}
0    & -1 & -1 \\
-p_1 & 0  & 0 \\
-p_2 & 0  & 0
\end{pmatrix},\quad
\hat{A}_\infty^{(-1)}=
\begin{pmatrix}
-p_1-p_2 & -q_1 & -q_2 \\
p_1q_1-\theta^\infty_2 & p_1-t_1 & p_1 \\
p_2q_2-\theta^\infty_3 & p_2 & p_2-t_2
\end{pmatrix}. \\
\hat{B}_1&=
\begin{pmatrix}
0 & -1 & 0 \\
-p_1 & \frac{1}{t_1-t_2}\{p_2(q_1-q_2)+\theta^\infty_3\}+q_1 & \frac{1}{t_1-t_2}\{p_1(q_1-q_2)-\theta^\infty_2\} \\
0 & \frac{1}{t_1-t_2}\{p_2(q_2-q_1)-\theta^\infty_3\} & \frac{1}{t_1-t_2}\{p_1(q_2-q_1)+\theta^\infty_2\}
\end{pmatrix},\\
\hat{B}_2&=
\begin{pmatrix}
0 & 0 & -1 \\
0 & \frac{1}{t_2-t_1}\{p_2(q_1-q_2)+\theta^\infty_3\} & \frac{1}{t_2-t_1}\{p_1(q_1-q_2)-\theta^\infty_2\} \\
-p_2 & \frac{1}{t_2-t_1}\{p_2(q_2-q_1)-\theta^\infty_3\} & \frac{1}{t_2-t_1}\{p_1(q_2-q_1)+\theta^\infty_2\}+q_2
\end{pmatrix},\\
E_2&=
{\rm diag}(0,1,0),\qquad
E_3=
{\rm diag}(0,0,1).
\end{align*}

The Hamiltonians are given by

\begin{align}
&H_{\mathrm{Gar},t_1}^{\frac{5}{2}+1+1}
\left({\theta^\infty_2,\theta^\infty_3};{t_1 \atop t_2};
{q_1,p_1 \atop q_2,p_2}\right)\\ &\quad=
H_{\rm II}\left(-\theta^\infty_2;t_1;q_1,p_1\right)
+p_1p_2+\frac{1}{t_1-t_2}(p_1(q_1-q_2)-\theta^\infty_2)(p_2(q_2-q_1)-\theta^\infty_3),\nonumber\\
&H_{\mathrm{Gar},t_1}^{\frac{5}{2}+1+1}
\left({\theta^\infty_2,\theta^\infty_3};{t_1 \atop t_2};
{q_1,p_1 \atop q_2,p_2}\right)\\ &\quad=
H_{\rm II}\left(-\theta^\infty_3;t_2;q_2,p_2\right)
+p_1p_2+\frac{1}{t_2-t_1}(p_1(q_1-q_2)-\theta^\infty_2)(p_2(q_2-q_1)-\theta^\infty_3).\nonumber
\end{align}

The gauge parameters $u,v$ satisfy

\begin{equation}
\frac{\partial u}{\partial t_1}=0,\quad
\frac{\partial v}{\partial t_1}=0,\quad
\frac{\partial u}{\partial t_2}=0,\quad
\frac{\partial v}{\partial t_2}=0.
\end{equation}

\subsection{Sasano system}
The next family of Painlev\'e-type equations are the Sasano system and its degenerated systems.
The Sasano system is derived from the Fuchsian equation of type $31,22,22,1111$.
They are given by the Lax pairs of $4\times 4$ matrices.

\subsubsection*{Singularity pattern: 1+1+1+1}
\noindent
\underline{Spectral type: 31,22,22,1111}

The Riemann scheme is given by
\[
 \left(\begin{array}{cccc}
  x=0 & x=1 & x=t & x=\infty \\
  0 & 0 & 0 & \theta^{\infty}_1 \\
  0 & 0 & 0 & \theta^{\infty}_2 \\
  0 & \theta^1 & \theta^t & \theta^{\infty}_3 \\
  \theta^0 & \theta^1 & \theta^t & \theta^{\infty}_4
       \end{array}\right) ,
\]
and the Fuchs relation is written as
$\theta^0 +2\theta^1+ 2\theta^t +\theta_1^\infty +\theta_2^\infty
+\theta_3^\infty +\theta_4^\infty =0$.

The Lax pair is expressed as
\begin{equation}
\left\{
\begin{aligned}
\frac{\partial Y}{\partial x}&=
\left(
\frac{A_0}{x}+\frac{A_1}{x-1}+\frac{A_t}{x-t}
\right)Y, \\
\frac{\partial Y}{\partial t}&=-\frac{A_t}{x-t}Y.
\end{aligned}
\right.
\end{equation}
Here $A_0$, $A_1$, and $A_t$ is given as follows:
\begin{align*}
A_{\xi}&=
U^{-1}
P^{-1}\hat{A}_{\xi}P
U,\quad
(\xi=0,1,t),\quad U={\rm diag}(1,u,v,w)\\
\hat{A}_0&=
\begin{pmatrix}
1 \\
0 \\
0 \\
0
\end{pmatrix}
\left(\theta^0,\ (\hat{A}_0)_{12},\ -f_1+\frac{q_1f_1+q_2-q_1}{t},\ -1+\frac{q_1}{t}\right),\\
\hat{A}_{\xi}&=
\begin{pmatrix}
I_2 \\
\hat{B}_{\xi}
\end{pmatrix}
(\theta^{\xi}I_2-\hat{C}_{\xi}\hat{B}_{\xi},\ \hat{C_{\xi}}), \quad(\xi=1,t),\quad
P=
\begin{pmatrix}
1 & 0 & 0 & 0\\
\frac{\hat{a}^{\infty}_2}{\theta^{\infty}_1-\theta^{\infty}_2} & 1 & 0 & 0 \\
\frac{\hat{a}^{\infty}_3}{\theta^{\infty}_1-\theta^{\infty}_3} & 0 & 1 & 0 \\
\frac{\hat{a}^{\infty}_4}{\theta^{\infty}_1-\theta^{\infty}_4} & 0 & 0 & 1
\end{pmatrix},
\end{align*}
\begin{align*}
\hat{B}_1&=
{\small
\begin{pmatrix}
p_2q_2-\theta^{\infty}_3 & -p_2(q_2-q_1)+f_2+\theta^{\infty}_3\\
(1-f_1)(p_2q_2-\theta^{\infty}_4)+p_1q_2 & -(q_2-q_1)(p_1+p_2(1-f_1))-\theta^{\infty}_4f_1+f_3
\end{pmatrix},
}\\
\hat{C}_1&=
\begin{pmatrix}
f_1 & 1 \\
f_1-1 & 1
\end{pmatrix},\quad
\hat{B}_t=
\begin{pmatrix}
tp_2 & f_2 \\
t(p_1+p_2(1-f_1)) & f_3
\end{pmatrix},\quad
\hat{C}_t=
\begin{pmatrix}
\frac{q_1(1-f_1)-q_2}{t} & -\frac{q_1}{t}\\
1-f_1 & -1
\end{pmatrix},
\end{align*}
\begin{align*}
(\theta^{\infty}_3-\theta^{\infty}_4)f_1&=p_1(q_2-q_1)+\theta^1+\theta^t+\theta^{\infty}_2+\theta^{\infty}_3,\\
(\theta^{\infty}_3-\theta^{\infty}_2)f_2&=(p_2(q_2-q_1)-\theta^{\infty}_3)(p_2(q_1(1-f_1)-q_2)
+\theta^1+\theta^{\infty}_3)\\
&\hspace{15mm}-p_2q_1((p_1+p_2(1-f_1))(q_2-q_1)+\theta^{\infty}_4f_1),\\
(\theta^{\infty}_4-\theta^{\infty}_2)f_3&=(p_1+p_2(1-f_1))
\{(\theta^{\infty}_3-\theta^{\infty}_4)q_1f_1\\
&\hspace{10mm}+(q_1-q_2)(q_1p_1+q_2p_2-\theta^1-\theta^{\infty}_3-\theta^{\infty}_4)\}
+\theta^{\infty}_4(\theta^1+\theta^{\infty}_4)f_1,\\
(\hat{A_0})_{12}&=-(p_1+p_2)(q_2-q_1)+(\theta^{\infty}_3-\theta^{\infty}_4)f_1
+f_2\left(f_1-\frac{q_1f_1+q_2-q_1}{t}\right)+f_3\left(1-\frac{q_1}{t}\right),
\end{align*}
\begin{align*}
\hat{a}^{\infty}_2&=(\hat{C_1}\hat{B_1}+\hat{C_t}\hat{B_t}-(\theta^1+\theta^t)I_2)_{21},\quad
\hat{a}^{\infty}_3=(\hat{B_1}(\hat{C_1}\hat{B_1}-\theta^1I_2)+\hat{B_t}(\hat{C_t}\hat{B_t}-\theta^tI_2))_{11},\\
\hat{a}^{\infty}_4&=(\hat{B_1}(\hat{C_1}\hat{B_1}-\theta^1I_2)+\hat{B_t}(\hat{C_t}\hat{B_t}-\theta^tI_2))_{21}.
\end{align*}

The Hamiltonian is given by
\begin{align}
&t(t-1)H_{\mathrm{Ss}}^{D_6}
 \left({\theta^{\infty}_{1}-\theta^{\infty}_{2},\theta^{\infty}_{2}-\theta^{\infty}_{3},\theta^{\infty}_{3}-\theta^{\infty}_{4} \atop\theta^{\infty}_{4},
 \theta^{1},\theta^{t}};t;
{q_1,p_1\atop q_2,p_2}\right)\\
&=t(t-1)H_{\rm VI} \left({-\theta^{0}-\theta^{1}-\theta^{t}-\theta^{\infty}_{1}-\theta^{\infty}_{3} ,-\theta^{t}-\theta^{\infty}_{2}+\theta^{\infty}_{3}\atop -\theta^{1}-\theta^{\infty}_{2}-\theta^{\infty}_{3},\theta^{0}+\theta^{1}+\theta^{t}+\theta^{\infty}_{2}+\theta^{\infty}_{3}+1
 };t;
q_1,p_1\right)\nonumber\\
&
\quad+t(t-1)H_{\rm VI} \left({ \theta^{\infty}_{3},\theta^{1}\atop \theta^{t},-\theta^{\infty}_{1}+\theta^{\infty}_{2}-\theta^{\infty}_{3}+\theta^{\infty}_{4}+1
 };t;
q_2,p_2\right)\nonumber\\
&\quad+2(q_{1}-1)p_{2}q_{2}(p_{1}(q_{1}-t)+\theta^{0}+\theta^{1}+\theta^{t}+\theta^{\infty}_{1}+\theta^{\infty}_{3}).\nonumber
\end{align}

The gauge parameters $u, v, w$ satisfy
\begin{align}
\frac{1}{u}\frac{du}{dt}&=\frac{1}{\hat{a}_{2}^{\infty}}\left(\frac{d\hat{a}_{2}^{\infty}}{dt}+(\theta^{\infty}_{1}-\theta^{\infty}_{2})p_{1}\right),\\
\frac{1}{v}\frac{dv}{dt}&=\frac{1}{\hat{a}_{3}^{\infty}}\left(\frac{d\hat{a}_{3}^{\infty}}{dt}+(\theta^{\infty}_{1}-\theta^{\infty}_{3})(f_{2}p_{1}+p_{2}(p_{1}q_{1}+p_{2}q_{2}+\theta^{t})\right),\\
\frac{1}{w}\frac{dw}{dt}&=\frac{1}{\hat{a}_{4}^{\infty}}\left(\frac{d\hat{a}_{4}^{\infty}}{dt}+(\theta^{\infty}_{1}-\theta^{\infty}_{4})(f_{3}p_{1}+(p_{1}+p_{2}(1-f_{1}))(p_{1}q_{1}+p_{2}q_{2}+\theta^{t})\right).
\end{align}

\subsubsection*{Singularity pattern: 2+1+1}
\noindent
\underline{Spectral type: $(2)(2),31,1111$}

The Riemann scheme is given by
\[
\left(
\begin{array}{ccc}
  x=0 & x=1 & x=\infty \\
\begin{array}{c} 0 \\ 0 \\ 0 \\ \theta^0 \end{array} &
\overbrace{\begin{array}{cc}
0 & 0\\
0 & 0\\
t & \theta^1\\
t & \theta^1
      \end{array}}&
\begin{array}{c} \theta^{\infty}_1 \\ \theta^{\infty}_2 \\ \theta^{\infty}_3\\ \theta^{\infty}_4 \end{array} 
\end{array}
\right) ,
\]
and the Fuchs-Hukuhara relation is written as
$\theta^0 +2\theta^1 +\theta_1^\infty +\theta_2^\infty +\theta_3^\infty
+\theta_4^\infty =0$.

The Lax pair is expressed as
\begin{equation}
\left\{
\begin{aligned}
\frac{\partial Y}{\partial x}&=
\left(
\frac{A_1^{(-1)}}{(x-1)^2}+\frac{A_1^{(0)}}{x-1}+\frac{A_0^{(0)}}{x}
\right)Y, \\
\frac{\partial Y}{\partial
 t}&=-\frac{1}{x-1}\left(\frac{A_1^{(-1)}}{t}\right)Y,
\end{aligned}
\right.
\end{equation}
where
\begin{align*}
A_{\xi}^{(k)}&=
U^{-1}
P^{-1}\hat{A_{\xi}}^{(k)}PU, \quad
P=\left(
\begin{array}{llll}
  \ 1 & 0 & 0 & 0 \\
 -\frac{a_{21}}{\theta^{\infty}_{1}-\theta^{\infty}_{2}} & 1 & 0 & 0 \\
 -\frac{a_{31}}{\theta^{\infty}_{1}-\theta^{\infty}_{3}} & 0 & 1 & 0 \\
 -\frac{a_{41}}{\theta^{\infty}_{1}-\theta^{\infty}_{4}} & 0 & 0 & 1 
\end{array}
\right),\quad
U={\rm diag}(1,u,v,w) ,
\end{align*}
\begin{align*}
\hat{A_1}^{(-1)}&=
\begin{pmatrix}
I_2 \\
\hat{B}_{1}
\end{pmatrix}
\begin{pmatrix}
t I_2+\hat{C}_{1}\hat{B}_{1} &-\hat{C}_{1}
\end{pmatrix},\quad
\hat{C}_{1}=
\begin{pmatrix}
f_{1} & 1\\
f_{1}-1 & 1
\end{pmatrix},\
\hat{B}_{1}=
\begin{pmatrix}
p_{2} & f_{2}\\
p_{1}+(1-f_{1})p_{2} & f_{3}
\end{pmatrix},
\\
\hat{A_1}^{(0)}&=
\begin{pmatrix}
-\theta^{0}-\theta^{\infty}_{1} & a_{12} & f_{1}(1-q_{1})+q_{1}-q_{2} & 1-q_{1}\\
\theta^{1}+\theta^{\infty}_{2}+\theta^{\infty}_{4}-p_{1}(q_{1}-1) & -\theta^{\infty}_{2} & 0 & 0\\
a_{31} & 0 & -\theta^{\infty}_{3}& 0 \\
a_{41} & 0 & 0 & -\theta^{\infty}_{4} 
\end{pmatrix},\\
\hat{A_0}^{(0)}&=
\begin{pmatrix}
1\\
0\\
0\\
0
\end{pmatrix}
\begin{pmatrix}
\theta^{0} & f_{4} & f_{1}(q_{1}-1)+q_{2}-q_{1} & q_{1}-1
\end{pmatrix},
\end{align*}
\begin{align*}
a_{12}&= (q_{1}-1)(f_{1}f_{2}+f_{3})-(q_{1}-q_{2})(f_{2}+p_{2})-(\theta^{1}+\theta^{\infty}_{2}+\theta^{\infty}_{3}),\\
a_{31}&=(p_{2}+f_{2})p_{1}(1-q_{1})+p_{2}((t+p_{2})(1-q_{2})-\theta^{0}-\theta^{\infty}_{1}+\theta^{\infty}_{3})+t\theta^{\infty}_{3}\\
&\quad +f_{2}(\theta^{1}+\theta^{\infty}_{2}+\theta^{\infty}_{4}),\\
a_{41}&=(p_{1}-(f_{1}-1)p_{2})((p_{2}+t)(1-q_{2})+p_{1}(1-q_{1})-\theta^{0}-\theta^{\infty}_{1}+\theta^{\infty}_{4})\\
&\quad -f_{3}(p_{1}(q_{2}-1)-(f_{1}-1)(\theta^{\infty}_{3}-\theta^{\infty}_{4}))-(f_{1}-1)t\theta^{\infty}_{4},
\end{align*}
\begin{align*}
f_{1}&=\frac{p_{1}(q_{2}-q_{1})+\theta^{1}+\theta^{\infty}_{2}+\theta^{\infty}_{3}}{\theta_{3}^{\infty}-\theta_{4}^{\infty}},\quad
f_{2}=\frac{p_{2}((p_{2}+t)(q_{1}-q_{2})+\theta^{1}+\theta^{\infty}_{2}+\theta^{\infty}_{3})+t\theta^{\infty}_{3}}{\theta^{\infty}_{3}-\theta^{\infty}_{2}},\\
f_{3}&=\frac{(p_{1}+p_{2}(1-f_{1}))((p_{2}+t)(q_{1}-q_{2})+\theta^{1}+\theta^{\infty}_{2}+\theta^{\infty}_{3})-t \theta^{\infty}_{4}f_{1}}{\theta^{\infty}_{4}-\theta^{\infty}_{2}},\\
f_{4}&=(p_{2}+f_{2})(q_{1}-q_{2})+(f_{1}f_{2}+f_{3})(1-q_{1})+(\theta^{1}+\theta^{\infty}_{2}+\theta^{\infty}_{3}).
\end{align*}

The Hamiltonian is given by
\begin{align}
&tH_{\mathrm{Ss}}^{D_5}\left({\theta^{0}+\theta^{1}+\theta^{\infty}_{2}+\theta^{\infty}_{3},-\theta^{\infty}_{1}
+\theta^{\infty}_{4}+1,\theta^{\infty}_{1}-\theta^{\infty}_{2}-\theta^{\infty}_{3}\atop
 -\theta^{0}-\theta^{1}-\theta^{\infty}_{2}-\theta^{\infty}_{4},-\theta^{0}-\theta^{1}-\theta^{\infty}_{3}-\theta^{\infty}_{4}-1};t;
{q_1,p_1\atop q_2,p_2}\right)\\
&
\quad=tH_{\rm V}\left({-\theta^{0}-\theta^{1}-\theta^{\infty}_{3}-\theta^{\infty}_{4}-1,\theta^{0}+\theta^{1}+\theta^{\infty}_{1}+\theta^{\infty}_{2}+\theta^{\infty}_{3}-\theta^{\infty}_{4}-1\atop
-\theta^{\infty}_{1}+\theta^{\infty}_{4}+1};t;q_1,p_1\right)\nonumber\\
&
\quad+tH_{\rm V}\left({\theta^{0}+\theta^{1}+2\theta^{\infty}_{1}+\theta^{\infty}_{3}-1,\theta^{0}+\theta^{1}+\theta^{\infty}_{1}+\theta^{\infty}_{2}-\theta^{\infty}_{3}+\theta^{\infty}_{4}-1\atop
-\theta^{0}-\theta^{1}-2\theta^{\infty}_{1}+1};t;q_2,p_2\right)\nonumber\\
&\quad+2p_{2}q_{1}(p_{1}(q_{1}-1)+\theta^{0}+\theta^{1}+\theta^{\infty}_{1}+\theta^{\infty}_{3}).\nonumber
\end{align}

The gauge parameters $u, v, w$ satisfy
\begin{align}
-\frac{t}{u}\frac{du}{dt}&=(t+2 p_{2})(1-q_{1})+p_{1},\quad
-\frac{t}{v}\frac{dv}{dt}=(t+2p_{2})(1-q_{2})+p_{1}+\theta^{1}+2\theta^{\infty}_{3},\\
-\frac{t}{w}\frac{dw}{dt}&=(t+2p_{1}+2p_{2})(1-q_{1})-\theta^{1}-2\theta^{\infty}_{4}.\nonumber
\end{align}

\medskip
\noindent
\underline{Spectral type: $(11)(11),31,22$}

The Riemann scheme is given by
\[
\left(
\begin{array}{ccc}
  x=0 & x=1 & x=\infty \\
\begin{array}{c} 0 \\ 0 \\ 0 \\ \theta^{0} \end{array} &
\begin{array}{c} 0 \\ 0 \\ \theta^{1} \\ \theta^{1} \end{array} &
\overbrace{\begin{array}{cc}
0 & \theta^\infty_1\\
0 & \theta^\infty_2\\
t & \theta^\infty_3\\
t & \theta^\infty_4
      \end{array}}
\end{array}
\right) ,
\]
and the Fuchs-Hukuhara relation is written as
$\theta^0 +2\theta^1 +\theta_1^\infty +\theta_2^\infty +\theta_3^\infty
+\theta_4^\infty =0$.

The Lax pair is expressed as
\begin{equation}
\left\{
\begin{aligned}
\frac{\partial Y}{\partial x}&=
\left(
\frac{A_0}{x}+\frac{A_1}{x-1}+A_{\infty}
\right)Y, \\
\frac{\partial Y}{\partial t}&=(-E_2\otimes I_2x+B_1)Y,
\end{aligned}
\right.
\end{equation}
where
\[
A_\xi=U^{-1}\hat{A}_\xi U,\ B_1=U^{-1}\hat{B}_1 U,\ U=\mathrm{diag}(1,u,v,w),
\]
\begin{align*}
\hat{A}_0&=
\begin{pmatrix}
-f_1-p_1q_1 \\
-f_2-p_2q_2+f_5(f_1+p_1q_1) \\
1 \\
\frac{f_2}{\theta^\infty_4-\theta^\infty_3}
\end{pmatrix}\times\\
&\quad\begin{pmatrix}
\frac{-f_2+p_2(q_2-1)q_1-\theta^1-\theta^\infty_1-\theta^\infty_4}{(\theta^\infty_1-\theta^\infty_2)q_1}+1
& 1-\frac{1}{q_2} & f_1-\theta^1-\theta^\infty_3+\frac{f_2(f_3+\theta^1+\theta^\infty_3)}{\theta^\infty_3-\theta^\infty_4} & f_3+\theta^1+\theta^\infty_4
\end{pmatrix},\\
\hat{A}_1&=
\begin{pmatrix}
-f_1 & -f_3 \\
\frac{1}{\theta^\infty_1-\theta^\infty_2}\left\{(1-\frac{1}{q_1})f_4+\theta^\infty_2 f_2\right\}&
\frac{1}{\theta^\infty_1-\theta^\infty_2}\left\{(\frac{1}{q_1}-\frac{1}{q_2})f_4
+(\theta^1+\theta^\infty_1+\theta^\infty_4)\theta^\infty_2\right\}\\
1 & 0 \\
\frac{f_2}{\theta^\infty_4-\theta^\infty_3} & 1
\end{pmatrix}\times\\
&\quad\begin{pmatrix}
-f_5-1 & -1 & \frac{f_2(f_3+\theta^1+\theta^\infty_3)}
{\theta^\infty_4-\theta^\infty_3}-f_1+\theta^1& -f_3-\theta^1-\theta^\infty_4\\
-f_5 & -1 & \frac{\theta^\infty_3 f_2}{\theta^\infty_4-\theta^\infty_3}&-\theta^\infty_4 
\end{pmatrix},\\
A_{\infty}&=
\begin{pmatrix}
O_2 & O_2 \\
O_2 & -t I_2
\end{pmatrix}, \quad
\hat{B}_1=
\begin{pmatrix}
O & \frac{1}{t}[\hat{A}_0+\hat{A}_1]_{1,2}\\
\frac{1}{t}[\hat{A}_0+\hat{A}_1]_{2,1} & O
\end{pmatrix},
\end{align*}
\begin{align*}
f_1&=p_1q_1(q_1-1)-\theta^\infty_1q_1,\quad
f_2=p_2q_2(q_2-1)-(\theta^1+\theta^\infty_2+\theta^\infty_4)q_2,\quad
f_3=\frac{q_1}{q_2}(p_1(q_2-q_1)+\theta^\infty_1),\\
f_4&=p_2q_2f_1-p_1q_1f_2,\quad
f_5=\frac{1}{(\theta^\infty_1-\theta^\infty_2)q_1}
(p_2q_2(q_1-q_2)+(\theta^1+\theta^\infty_2+\theta^\infty_4)q_2).
\end{align*}

The Hamiltonian is given by
\begin{align}
&tH_{\mathrm{NY}}^{A_5}
 \left({\theta^1,\theta^\infty_1,-\theta^1-\theta^\infty_1-\theta^\infty_4
\atop \theta^1+\theta^\infty_2+\theta^\infty_4, -\theta^1-\theta^\infty_2-\theta^\infty_3};t;
{q_1,p_1\atop q_2,p_2}\right)\\
&=tH_{\rm V}
 \left({\theta^1+\theta^\infty_1, \theta^0+\theta^1 \atop
-\theta^1};t;q_1,p_1\right)
+tH_{\rm V}\left({\theta^1-\theta^\infty_1+\theta^\infty_2,\theta^0+\theta^1 \atop
\theta^\infty_1+\theta^\infty_4};t;q_2,p_2\right)
+2p_1p_2q_1(q_2-1).\nonumber
\end{align}

The gauge parameters $u, v, w$ satisfy
\begin{align}
-tq_1q_2\frac{1}{u}\frac{du}{dt}&=
q_1q_2(p_2-2p_1+\theta^\infty_2-\theta^\infty_1)+q_1(2p_1q_1-\theta^{\infty}_1-\theta^{\infty}_4)-\theta^1q_2,\\
tq_1\frac{1}{v}\frac{dv}{dt}&=
q_1p_2(1-2q_2)-2p_1q_1(q_1-1)+(t+\theta^1+2\theta^\infty_1+\theta^\infty_2+\theta^\infty_4)q_1
+\theta^1,\\
tq_1q_2\frac{1}{w}\frac{dw}{dt}&=
q_1q_2(2p_1+t+\theta^\infty_1-\theta^\infty_4)-q_1(2p_1q_1-\theta^\infty_1-\theta^\infty_4)
+\theta^1q_2.
\end{align}

\medskip
\noindent
\underline{Spectral type: $(111)(1),22,22$}

The Riemann scheme is given by
\[
\left(
\begin{array}{ccc}
  x=0 & x=1 & x=\infty \\
\begin{array}{c} 0 \\ 0 \\ \theta^{0} \\ \theta^{0} \end{array} &
\begin{array}{c} 0 \\ 0 \\ \theta^{1} \\ \theta^{1} \end{array} &
\overbrace{\begin{array}{cc}
t & \theta^{\infty}_{1}\\
0 &  \theta^{\infty}_{2}\\
0 & \theta^{\infty}_{3}\\
0 & \theta^{\infty}_{4}
      \end{array}}
\end{array}
\right) ,
\]
and the Fuchs-Hukuhara relation is written as
$2\theta^0 +2\theta^1+\theta_1^\infty +\theta_2^\infty +\theta_3^\infty
+\theta_4^\infty =0$.

The Lax pair is expressed as
\begin{equation}
\left\{
\begin{aligned}
\frac{\partial Y}{\partial x}&=
\left(
\frac{A_0}{x}+\frac{A_1}{x-1}+A_{\infty}
\right)Y, \\
\frac{\partial Y}{\partial t}&=(-E_{1}x+B)Y.
\end{aligned}
\right.
\end{equation}
Here
\begin{align*}
A_{\xi}&=
U^{-1}
\hat{A_{\xi}}
U, \quad
B=U^{-1}
\hat{B}
U,\quad
U={\rm diag}(1,u,v,w),\\
\hat{A_{\xi}}&=
\begin{pmatrix}
I_2 \\
\hat{B_{\xi}}
\end{pmatrix}
(\theta^{\xi}I_2-\hat{C_{\xi}}\hat{B_{\xi}},\ \hat{C_{\xi}}) \quad(\xi=0,1),\quad
A_{\infty}=-tE_{1}, \quad 
E_{1}={\rm diag}(1,0,0,0),
\end{align*}
\begin{align*}
\hat{B_0}&=
\begin{pmatrix}
\mu_{2} & f_{2}\\
\mu_{1}-(f_{1}-1)\mu_{2} & f_{3}
\end{pmatrix}, \quad
\hat{C_0}=
\begin{pmatrix}
(1-f_{1})\lambda_{1}-\lambda_{2} & -\lambda_{1} \\
1-f_1 & -1
\end{pmatrix},\\
\hat{B_1}&=
\begin{pmatrix}
\mu_{2}\lambda_{2}-\theta^{\infty}_{3} & f_{2}+\mu_{2}(\lambda_{1}-\lambda_{2})+\theta^{\infty}_{3} \\
\mu_{1}\lambda_{2}+(1-f_{1})(\mu_{2}\lambda_{2}-\theta^{\infty}_{4}) & f_{3}+(\lambda_{1}-\lambda_{2})(\mu_{2}(1-f_{1})+\mu_{1})-f_{1}\theta^{\infty}_{4}
\end{pmatrix},
\hat{C_1}=
\begin{pmatrix}
f_{1} & 1\\
f_{1}-1 & 1
\end{pmatrix},\\
\hat{B}&=\frac{1}{t} 
\begin{pmatrix}
0 & (\hat{A}_{0}+\hat{A}_{1})_{12} & (\hat{A}_{0}+\hat{A}_{1})_{13} & (\hat{A}_{0}+\hat{A}_{1})_{14} \\
(\hat{A}_{0}+\hat{A}_{1})_{21} & 0 & 0 & 0\\
(\hat{A}_{0}+\hat{A}_{1})_{31} & 0 & 0 & 0\\
(\hat{A}_{0}+\hat{A}_{1})_{41} & 0 & 0 & 0
\end{pmatrix},
\end{align*}
\begin{align*}
f_{1}&=\frac{\mu_{1}(\lambda_{2}-\lambda_{1})+\theta^{0}+\theta^{1}+\theta^{\infty}_{2}+\theta^{\infty}_{3}}{\theta^{\infty}_{3}-\theta^{\infty}_{4}},\\
f_{2}&=\frac{\mu_{2}\lambda_{2}(\mu_{2}(\lambda_{2}-\lambda_{1})-\theta^{1}-2\theta^{\infty}_{3})-\mu_{2}\lambda_{1}(\theta^{0}+\theta^{\infty}_{2}-\theta^{\infty}_{3})+\theta^{\infty}_{3}(\theta^{1}+\theta^{\infty}_{3})}{\theta^{\infty}_{2}-\theta^{\infty}_{3}},\\
f_{3}&=\frac{(\lambda_{2}(\mu_{2}(\lambda_{1}-\lambda_{2})+\theta^{1}+\theta^{\infty}_{3}+\theta^{\infty}_{4})+(\theta^{0}+\theta^{\infty}_{2}-\theta^{\infty}_{4})\lambda_{1})((f_{1}-1)\mu_{2}-\mu_{1})-f_{1}\theta^{\infty}_{4}(\theta^{1}+\theta^{\infty}_{4})}{\theta^{\infty}_{2}-\theta^{\infty}_{4}}.
\end{align*}

The Hamiltonian is given by
\begin{align}
&t\tilde{H}_{\mathrm{Ss}}^{D_5}
 \left({\theta^{0},\theta^{1},\theta^{\infty}_{1},\atop\theta^{\infty}_{2},
 \theta^{\infty}_{3}, };t;
{\lambda_1,\mu_1\atop \lambda_1,\mu_1}\right)\\
&=t\tilde{H}_{\rm V}
 \left({-\theta^{1}-\theta^{\infty}_{2}-\theta^{\infty}_{3},\theta^{0}+\theta^{\infty}_{2}-\theta^{\infty}_{3}\atop -2\theta^{0}-\theta^{1}-\theta^{\infty}_{1}-\theta^{\infty}_{2}};t;\lambda_1,\mu_1\right)
+t\tilde{H}_{\rm V}\left({\theta^{0},\theta^{1}\atop
\theta^{\infty}_{3}};t;\lambda_2,\mu_2\right)\nonumber\\
&\quad+2\mu_{2}\lambda_{2}(\lambda_{1}-1)(\mu_{1}(\lambda_{1}-1)+\theta^{0}+\theta^{1}+\theta^{\infty}_{1}+\theta^{\infty}_{3}).\nonumber
\end{align}

The gauge parameters $u, v, w$ satisfy
\begin{align}
-\frac{t}{u}\frac{du}{dt}&=(\lambda_{1}-1)(2\mu_{2}\lambda_{2}+\mu_{1}(\lambda_{1}-1)-\theta^{1}-\theta^{\infty}_{3}-\theta^{\infty}_{4})+t-\theta^{\infty}_{1}+\theta^{\infty}_{2},\\
-\frac{t}{v}\frac{dv}{dt}&=(\lambda_{1}-1)(\mu_{1}(\lambda_{1}-1)+\theta^{0}+\theta^{1}+\theta^{\infty}_{1}+\theta^{\infty}_{3})+\lambda_{2}(2\mu_{2}(\lambda_{2}-1)-\theta^{1}-2\theta^{\infty}_{3})+t-\theta^{0}-\theta^{\infty}_{1}+\theta^{\infty}_{3},\\
-\frac{t}{w}\frac{dw}{dt}&=2(\lambda_{1}-1)(\mu_{1}\lambda_{1}+\mu_{2}\lambda_{2})+\lambda_{1}(\theta^{0}+\theta^{\infty}_{1}-\theta^{\infty}_{4})+t-2\theta^{0}-\theta^{1}-2\theta^{\infty}_{1}.
\end{align}

When we change the canonical variables as
\begin{align*}
&\lambda_1 = 1-\frac{1}{q_{2}},\ \mu_1 = q_{2}(q_{2}p_{2}+\theta^{0}+\theta^{1}+\theta^{\infty}_{1}+\theta^{\infty}_{3}),\
\lambda_2 = 1-\frac{1}{q_{1}},\ \mu_{2} =
 q_{1}(q_{1}p_{1}-\theta^{1}-\theta^\infty_{3}) ,
\end{align*}
then we obtain
\begin{align}
tH_{\mathrm{Ss}}^{D_5}
 \left({\theta^{0},\theta^{1},\theta^{\infty}_{1},\atop\theta^{\infty}_{2},
 \theta^{\infty}_{3}, };t;
{q_1,p_1\atop q_2,p_2}\right)
&=tH_{\rm V}\left({\theta^{\infty}_{3},\theta^{0}-\theta^{1}\atop
\theta^{1}};t;q_1,p_1\right)\\
&\quad +tH_{\rm V}
 \left({-2\theta^{0}-3\theta^{1}-\theta^{\infty}_{1}-\theta^{\infty}_{2}-2\theta^{\infty}_{3},-\theta^{0}-\theta^{1}-2\theta^{\infty}_{2}\atop \theta^{0}+2\theta^{1}+\theta^\infty_{2}+\theta^{\infty}_{3}};t;q_2,p_2\right)\nonumber\\
&
\quad
+2p_{2}q_{1}(p_{1}(q_{1}-1)-\theta^{1}-\theta^{\infty}_{3}).\nonumber
\end{align}

\subsubsection*{Singularity pattern: 3+1}
\noindent
\underline{Spectral type: $((11))((11)),31$}

The Riemann scheme is given by
\[
\left(
\begin{array}{cc}
  x=0  & x=\infty \\
\begin{array}{c} 0 \\ 0 \\ 0 \\ \theta^{0} \end{array} &
\overbrace{\begin{array}{ccc}
0 & 0 & \theta^{\infty}_{1}\\
0 & 0 &\theta^{\infty}_{2}\\
1 & -t & \theta^{\infty}_{3}\\
1 & -t & \theta^{\infty}_{4}
      \end{array}}
\end{array}
\right) ,
\]
and the Fuchs-Hukuhara relation is written as
$\theta^0 +\theta_1^\infty +\theta_2^\infty +\theta_3^\infty
+\theta_4^\infty =0$.

The Lax pair is expressed as
\begin{equation}
\left\{
\begin{aligned}
\frac{\partial Y}{\partial x}&=
\left(
A_\infty^{(-2)}x+A_\infty^{(-1)}+\frac{A_0^{(0)}}{x}
\right)Y, \\
\frac{\partial Y}{\partial t}&=(E_2\otimes I_2x+B_1)Y.
\end{aligned}
\right.
\end{equation}
Here
\[
A_\xi^{(-k)}=U^{-1}\hat{A}_\xi^{(-k)} U, \quad B_1=U^{-1}\hat{B}_1 U, \quad U=\mathrm{diag}(1,u,v,w),
\]
\begin{align*}
\hat{A}_\infty^{(-2)}&=
-E_{2}\otimes I_{2}, \quad E_{2}={\rm diag}(0,1),\\
\hat{A}_\infty^{(-1)}&=
\begin{pmatrix}
0 & 0 & f_4
\left(f_3+\theta^\infty_3-\frac{\theta^\infty_3}{\theta^\infty_3-\theta^\infty_4}f_1\right)&
-f_4(f_2-p_2q_1+\theta^\infty_4)\\
0 & 0 & -f_4(f_3+\theta^\infty_3)+\frac{\theta^\infty_3}{\theta^\infty_3-\theta^\infty_4}
f_1(f_4+1) & f_4(f_2-p_2q_1+\theta^\infty_4)+\theta^\infty_4\\
1+\frac{1}{f_4} & 1 & t & 0 \\
\frac{f_1(1+f_4^{-1})}{\theta^\infty_4-\theta^\infty_3}+1 & 
\frac{f_1}{\theta^\infty_4-\theta^\infty_3}+1 & 0 & t
\end{pmatrix},\\
\hat{A}_0^{(0)}&=
\begin{pmatrix}
-p_2f_4 \\
-p_1+p_2f_4 \\
1 \\
\frac{f_1}{\theta^\infty_4-\theta^\infty_3}
\end{pmatrix}
\begin{pmatrix}
q_1+\frac{q_2}{f_4} & q_1 & f_3 & p_2q_1-f_2
\end{pmatrix},\quad
\hat{B}_1=
\begin{pmatrix}
O & -[A_\infty^{(-1)}]_{1,2}\\
-[A_\infty^{(-1)}]_{2,1} & O
\end{pmatrix},
\end{align*}
\begin{align*}
f_1&=p_1q_1-\theta^\infty_2-\theta^\infty_4,\quad
f_2=p_2q_2-\theta^\infty_1-\theta^\infty_4,\\
f_3&=\frac{(p_1q_1-\theta^\infty_2-\theta^\infty_3)(p_2q_2-\theta^\infty_1
-\theta^\infty_3)-p_2q_1(p_1q_1-\theta^\infty_2-\theta^\infty_4)}{\theta^\infty_4-\theta^\infty_3},\quad
f_4=\frac{f_1-p_1q_2}{\theta^\infty_2-\theta^\infty_1}.
\end{align*}

The Hamiltonian is given by
\begin{align}
&H_{\mathrm{NY}}^{A_4}
 \left({-\theta^\infty_2-\theta^\infty_3,\theta^\infty_2+\theta^\infty_4
\atop -\theta^\infty_1-\theta^\infty_4,\theta^\infty_1};t;
{q_1,p_1\atop q_2,p_2}\right)\\
&\quad=H_{\rm IV}
 \left(\theta^\infty_2+\theta^\infty_4,-\theta^\infty_2-\theta^\infty_3;t;q_1,p_1\right)
+H_{\rm IV}\left(\theta^\infty_1,\theta^0;t;q_2,p_2\right)+2q_1p_1p_2.\nonumber
\end{align}

The gauge parameters $u, v, w$ satisfy
\begin{align}
\frac{1}{u}\frac{du}{dt}=\frac{(\theta^\infty_1-\theta^\infty_2)p_1}{p_1q_2-f_1},\quad
\frac{1}{v}\frac{dv}{dt}=-q_1-t+\frac{(\theta^\infty_1-\theta^\infty_2)p_1}{p_1q_2-f_1},\quad
\frac{1}{w}\frac{dw}{dt}=p_1-t+\frac{(\theta^\infty_1-\theta^\infty_2)p_1}{p_1q_2-f_1}.
\end{align}

\subsubsection*{Singularity pattern: 2+2}
\noindent
\underline{Spectral type: $(2)(2),(111)(1)$}

The Riemann scheme is given by
\[
\left(
\begin{array}{cc}
  x=0 & x=\infty \\
\overbrace{
\begin{array}{cc} 0 & 0 \\ 0 & 0 \\ 1 & \theta^0\\ 1 & \theta^0\end{array}} &
\overbrace{\begin{array}{cc}
t & \theta^{\infty}_1\\
0 & \theta^{\infty}_2\\
0 & \theta^{\infty}_3\\
0 & \theta^{\infty}_4
      \end{array}}
\end{array}
\right) ,
\]
and the Fuchs-Hukuhara relation is written as
$2\theta^0 +\theta_1^\infty +\theta_2^\infty +\theta_3^\infty
+\theta_4^\infty =0$.

The Lax pair is expressed as
\begin{equation}
\left\{
\begin{aligned}
\frac{\partial Y}{\partial x}&=
\left(
\frac{A_0^{(-1)}}{x^2}+\frac{A_{0}^{(0)}}{x}+A_{\infty}
\right)Y, \\
\frac{\partial Y}{\partial t}&=(-E_{1}x+B)Y.
\end{aligned}
\right.
\end{equation}
Here
\begin{align*}
A_{\xi}^{(k)}&=
U^{-1}
\hat{A_{\xi}}^{(k)}
U, \quad
A_{\infty}=-tE_{1}, \quad
B=U^{-1}\hat{B}U, \quad
 U={\rm diag}(1,u,v,w), \quad E_{1}={\rm diag}(1,0,0,0),
\end{align*}
\begin{align*}
\hat{A_{0}}^{(0)}&=
\begin{pmatrix}
-\theta^{\infty}_{1} & a_{12} & (1-f_{1})q_{1}-q_{2} & -q_{1} \\
-p_{1}q_{1}+\theta^{0}+\theta^{\infty}_{2}+\theta^{\infty}_{4} & -\theta^{\infty}_{2} & 0 & 0 \\
a_{31} & 0 & -\theta^{\infty}_{3} & 0 \\
a_{41} & 0 & 0 & -\theta^{\infty}_{4}
\end{pmatrix},\\
\hat{A_{0}}^{(-1)}&=
\begin{pmatrix}
I_{2} \\ B_{0}^{(-1)}
\end{pmatrix}
\begin{pmatrix}
I_{2}-C_{0}^{(-1)}B_{0}^{(-1)} & C_{0}^{(-1)}
\end{pmatrix}, \quad 
B_{0}^{(-1)}=\begin{pmatrix}
p_{2} & f_{2} \\
p_{1}+(1-f_{1})p_{2} & f_{3}
\end{pmatrix}, \\
C_{0}^{(-1)}&=\begin{pmatrix}
f_{1} & 1 \\
f_{1}-1 & 1
\end{pmatrix}, \quad
\hat{B}=\frac{1}{t}
\begin{pmatrix}
0 & (A_{0}^{(-1)})_{12} & (A_{0}^{(-1)})_{13} & (A_{0}^{(-1)})_{14}\\
(A_{0}^{(-1)})_{21} & 0 & 0 & 0\\
(A_{0}^{(-1)})_{31} & 0 & 0 & 0 \\
(A_{0}^{(-1)})_{41} & 0 & 0 & 0
\end{pmatrix},
\end{align*}

\begin{align*}
&(\theta^{\infty}_{3}-\theta^{\infty}_{4})f_{1}=-p_{1}(q_{1}-q_{2})+\theta^{0}+\theta^{\infty}_{2}+\theta^{\infty}_{3},\\
&(\theta^{\infty}_{2}-\theta^{\infty}_{3})f_{2}=(q_{1}-q_{2})p_{2}(1-p_{2})-p_{2}(\theta^{0}+\theta^{\infty}_{2}+\theta^{\infty}_{3})+\theta^{\infty}_{3},\\
&(\theta^{\infty}_{4}-\theta^{\infty}_{2})f_{3}=((1-f_{1})p_{2}+p_{1})((q_{1}-q_{2})(p_{2}-1)+\theta^{0}+\theta^{\infty}_{2}+\theta^{\infty}_{3})+f_{1}\theta^{\infty}_{4},\\
&a_{12}=f_{3}q_{1}+p_{2}(q_{2}-q_{1})+f_{2}((f_{1}-1)q_{1}+q_{2})-\theta^{0}-\theta^{\infty}_{2}-\theta^{\infty}_{3},\\
&a_{31}=p_{2}(\theta^{\infty}_{3}-\theta^{\infty}_{1})+p_{2}q_{2}(1-p_{2})-p_{1}q_{1}(f_{2}+p_{2})-\theta^{\infty}_{3}+f_{2}(\theta^{0}+\theta^{\infty}_{2}+\theta^{\infty}_{4}),\\
&a_{41}=(p_{1}q_{1}+(p_{2}-1)q_{2}+\theta^{\infty}_{1}-\theta^{\infty}_{4})((f_{1}-1)p_{2}-p_{1})+ f_{3}(-p_{1}q_{1}+\theta^{0}+\theta^{\infty}_{2}+\theta^{\infty}_{4})+(f_{1}-1)\theta^{\infty}_{4} .
\end{align*}

The Hamiltonian is given by
\begin{align}
&tH_{\mathrm{Ss}}^{D_4}
 \left({\theta^0,\theta^\infty_1,\atop \theta^\infty_2, 
 \theta^\infty_3, };t;
{q_1,p_1\atop q_2,p_2}\right)\\
&=tH_{\rm III}(D_{6})
 \left(\theta^{0}+\theta^{\infty}_{1}+\theta^{\infty}_{3}, -\theta^{0}-2\theta^{\infty}_{4} 
;t;q_1,p_1\right)
+tH_{\rm III}(D_{6})\left(-\theta^\infty_3, 
-\theta^{0}-2\theta^{\infty}_{3};t;q_2,p_2\right)\nonumber\\
&\quad+2p_{2}q_{1}(p_{1}q_{1}+\theta^{0}+\theta^{\infty}_{1}+\theta^{\infty}_{3}).\nonumber
\end{align}

The gauge parameters $u, v, w$ satisfy
\begin{align}
\frac{t}{u}\frac{du}{dt}&=(1-2 p_{2})q_{1}+\theta^{\infty}_{1}-\theta^{\infty}_{2},\quad
\frac{t}{v}\frac{dv}{dt}=(1-2 p_{2})q_{2}+\theta^{0}+\theta^{\infty}_{1}+\theta^{\infty}_{3},\\
\frac{t}{w}\frac{dw}{dt}&=(1-2 p_{1}-2 p_{2})q_{1}-\theta^{0}-\theta^{\infty}_{2}-\theta^{\infty}_{3}.
\end{align}

\subsection{Matrix Painlev\'e systems}
The next family of Painlev\'e-type equations are the sixth matrix Painlev\'e system and its degenerated systems.
The  sixth matrix Painlev\'e system is derived from the Fuchsian equation of type $22,22,22,211$.
They are given by the Lax pairs of $4\times 4$ matrices.

\subsubsection*{Singularity pattern: 1+1+1+1}
\noindent
\underline{Spectral type: 22,22,22,211}

The Riemann scheme is given by
\[
 \left(\begin{array}{cccc}
  x=0 & x=1 & x=t & x=\infty \\
  0 & 0 & 0 & \theta^{\infty}_1 \\
  0 & 0 & 0 & \theta^{\infty}_1 \\
  \theta^0 & \theta^1 & \theta^t & \theta^{\infty}_2 \\
  \theta^0 & \theta^1 & \theta^t & \theta^{\infty}_3
       \end{array}\right) ,
\]
and the Fuchs relation is written as
$2\theta^0 +2\theta^1 +2\theta^t +2\theta_1^\infty +\theta_2^\infty
+\theta_3^\infty =0$.

The Lax pair is expressed as
\begin{equation}
\left\{
\begin{aligned}
\frac{\partial Y}{\partial x}&=
\left(
\frac{A_0}{x}+\frac{A_1}{x-1}+\frac{A_t}{x-t}
\right)Y ,\\
\frac{\partial Y}{\partial t}&=-\frac{A_t}{x-t}Y .
\end{aligned}
\right.
\end{equation}
Here $A_0$, $A_1$, and $A_t$ are given as
\begin{align*}
A_{\xi}&=
(U\oplus {\rm diag}(v,1))^{-1}X^{-1}\hat{A}_{\xi}X(U\oplus {\rm diag}(v,1)),
\quad(\xi=0,1,t),\\
\hat{A}_0&=
\begin{pmatrix}
I_2 \\
O
\end{pmatrix}
\left(\theta^0I_2,\frac{1}{t}Q-I_2\right),\quad
\hat{A}_1=
\begin{pmatrix}
I_2 \\
PQ-\Theta
\end{pmatrix}
(\theta^1I_2-PQ+\Theta,I_2),\\
\hat{A}_t&=
\begin{pmatrix}
I_2 \\
tP
\end{pmatrix}
\left(\theta^tI_2+QP,-\frac{1}{t}Q\right), \quad U \in \mathrm{GL}(2),\quad v \in \mathbb{C}^{\times}.
\end{align*}
Here, $Q$, $P$, and $\Theta$ are
\begin{align*}
Q=
\begin{pmatrix}
q_1 & 1\\
-q_2 & q_1
\end{pmatrix},
\quad
P=
\begin{pmatrix}
p_1/2 & -p_2\\
p_2 q_2-\theta-\theta^\infty_1-\theta^\infty_2 & p_1/2
\end{pmatrix},
\quad
\Theta=
\begin{pmatrix}
\theta^{\infty}_2 & \\
 & \theta^{\infty}_3
\end{pmatrix},
\end{align*}
where $\theta=\theta^0+\theta^1+\theta^t$.

The matrix $X$ is given by
$X=
\begin{pmatrix}
I_2 & O \\
Z & I_2
\end{pmatrix}$,
where
\begin{align*}
Z&=(\theta^{\infty}_1-\Theta)^{-1}
[-\theta^1(QP+\theta+\theta^{\infty}_1)
+(QP+\theta+\theta^{\infty}_1)^2
-t(PQ+\theta^t)P].
\end{align*}

The Hamiltonian is given by
\begin{multline}
t(t-1)H^{\mathrm{Mat}}_{\rm VI}\left({-\theta^0-\theta^t-\theta^\infty_1,-\theta^1,\theta^t \atop
\theta^0+1,\theta^1+\theta^\infty_3};t;
{q_1, p_1\atop q_2, p_2}\right)\\=
\mathrm{tr}\left[Q(Q-1)(Q-t)P^2+
\{(\theta^0+1 -(\theta+\theta^\infty_1+\theta^\infty_2)K)Q(Q-1)\right.\\
\left. +\theta^t(Q-1)(Q-t)+(\theta+2\theta^\infty_1-1)Q(Q-t)\}P
+(\theta+\theta^\infty_1)(\theta^0+\theta^t+\theta^\infty_1)Q\right].
\end{multline}

The gauge parameters $U$ and $v$ satisfy
\begin{align}
t(t-1)\frac{dU}{dt}&=MU,\\
M_{11}&=p_1(2q_1-t)(1-q_1)-(\theta^0+\theta^t+\theta^\infty_1-\theta^\infty_2)q_1
+2p_2q_2\nonumber\\
&\quad+2q_1p_2(q_1(q_1-t-1)+t-q_2)+(\theta^\infty_1-\theta^\infty_3-1)t
+\theta^0+\theta^t-\theta^\infty_2+\theta^\infty_3+1,\nonumber\\
M_{12}&=p_1(2q_1-t)+2p_2q_2+2q_1p_2(t-q_1)+\theta^0+\theta^t+\theta^\infty_1-\theta^\infty_2,\nonumber\\
M_{21}&=2(\theta+\theta^\infty_1+\theta^\infty_2)q_1(t-q_1)
-(2p_2q_2+\theta^0+\theta^t+\theta^\infty_1-\theta^\infty_2)q_2\nonumber\\
&\quad+p_1q_2(t-2q_1)+2q_1p_2q_2(q_1-t),\nonumber\\
M_{22}&=((q_1-t)p_1+\theta^0+\theta^t+\theta^\infty_1-\theta^\infty_2)q_1-2tp_2q_2
+q_2(4p_2q_1-p_1)+(\theta^0+\theta^1+\theta^\infty_1+\theta^\infty_2)t,\nonumber\\
t(t-1)\frac{1}{v}\frac{dv}{dt}&=
2q_1((t+1)p_1+\theta^1+2\theta^\infty_2)-p_1(3{q_1}^2+t)+2(t+1)p_2q_2
+2q_1p_2(q_1(q_1-t-1)+t-3q_2)\\
&\quad +p_1q_2+(\theta^0+\theta^1+2\theta^t+2\theta^\infty_1-1)t
+\theta^0+\theta^t-\theta^\infty_2+\theta^\infty_3+1.\nonumber
\end{align}

\subsubsection*{Singularity pattern: 2+1+1}
\noindent
\underline{Spectral type: (2)(2),22,211}

The Riemann scheme is given by
\[
\left(
\begin{array}{ccc}
  x=0 & x=1 & x=\infty \\
\begin{array}{c} 0 \\ 0 \\ \theta^0 \\ \theta^0 \end{array} &
\overbrace{\begin{array}{cc}
0 & 0\\
0 & 0\\
-t & \theta^1\\
-t & \theta^1
      \end{array}}&
\begin{array}{c} \theta^{\infty}_1 \\ \theta^{\infty}_1 \\ \theta^{\infty}_2 \\ \theta^{\infty}_3 \end{array} 
\end{array}
\right) ,
\]
and the Fuchs-Hukuhara relation is written as
$2\theta^0+2\theta^1 +2\theta_1^\infty +\theta_2^\infty +\theta_3^\infty =0$.

The Lax pair is expressed as
\begin{equation}
\left\{
\begin{aligned}
\frac{\partial Y}{\partial x}&=
\left(
\frac{A_0^{(0)}}{x}+\frac{A_1^{(-1)}}{(x-1)^2}+\frac{A_1^{(0)}}{x-1}
\right)Y ,\\
\frac{\partial Y}{\partial
 t}&=-\frac{1}{x-1}\left(\frac{A_1^{(-1)}}{t}\right)Y .
\end{aligned}
\right.
\end{equation}
Here
$A_0^{(0)}$, $A_1^{(-1)}$, and $A_1^{(0)}$ are given as follows:
\begin{align*}
A_{\xi}&=
(U\oplus {\rm diag}(v,1))^{-1}\hat{A}_{\xi}(U\oplus {\rm diag}(v,1)),\\
\hat{A}_1^{(-1)}&=
\begin{pmatrix}
I_2 \\
-Z+Q
\end{pmatrix}
(-t(I_2-Q)-tZ,-tI_2)
=G_1
\begin{pmatrix}
O_2 & O_2\\
O_2 & -tI_2
\end{pmatrix}
G_1^{-1},\\
\hat{A}_1^{(0)}&=
\begin{pmatrix}
-(\theta^0+\theta^{\infty}_1)I_2+(P+t)Z & P+t \\
\theta^0Z-Z(P+t)Z & -Z(P+t)-\Theta
\end{pmatrix}\\
&=G_1
\begin{pmatrix}
O_2 & (P+t)Q-\theta^0-\theta^1-\theta^\infty_1\\
(P+t)Q-P-\theta^0-\theta^\infty_1-t & \theta^1 I_2
\end{pmatrix}
G_1^{-1},\\
\hat{A}_0^{(0)}&=
\begin{pmatrix}
I_2 \\
-Z
\end{pmatrix}
(\theta^0I_2-(P+t)Z,-tI_2-P),\\
Z&=(\theta^{\infty}_1-\Theta)^{-1}
[-(P+t)Q(Q-1)+(2\theta^0+\theta^1+2\theta^{\infty}_1)Q
-\theta^0-\theta^1-\theta^{\infty}_1],
\end{align*}
\[
G_1=
\begin{pmatrix}
I_2 & I_2\\
Q-Z-I_2 & Q-Z
\end{pmatrix} ,\quad
\left(G_1^{-1}=
\begin{pmatrix}
Q-Z & -I_2\\
I_2-Q+Z & I_2
\end{pmatrix}\right).
\]
Here, $Q$, $P$, and $\Theta$ are
\begin{align*}
Q=
\begin{pmatrix}
q_1 & 1\\
-q_2 & q_1
\end{pmatrix},\quad
P=
\begin{pmatrix}
p_1/2 & -p_2\\
p_2 q_2-\theta^0-\theta^1-\theta^\infty_1-\theta^\infty_2 & p_1/2
\end{pmatrix},\quad
\Theta=
\begin{pmatrix}
\theta^{\infty}_2 & \\
 & \theta^{\infty}_3
\end{pmatrix}.
\end{align*}

The Hamiltonian is given by
\begin{align}
&tH^{\mathrm{Mat}}_{\rm V}\left({\theta^\infty_1-1,-2\theta^0-\theta^1-2\theta^\infty_1 \atop
\theta^0+\theta^1+\theta^\infty_1,-\theta^0-\theta^1-\theta^\infty_2-1}
;t;{q_1, p_1 \atop q_2, p_2}\right)\\
&=\mathrm{tr}[P(P+t)Q(Q-1)
+(\theta^0+\theta^1+\theta^{\infty}_1)P
-(\theta^0+\theta^1+2\theta^{\infty}_1-1)tQ
-(2\theta^0+\theta^1+2\theta^{\infty}_1)PQ].\nonumber
\end{align}

The gauge parameters $U$ and $v$ satisfy
\begin{align}
t\frac{dU}{dt}&=MU,\\
M_{11}&=(1-2q_1)\left(\frac{3}{2}p_1+2t\right)+2(q_1(q_1-1)p_2-2p_2q_2+\theta^0+\theta^\infty_1-\theta^\infty_3), \nonumber\\
M_{12}&=(2q_1-1)p_2-p_1-2t,\nonumber\\
M_{21}&=(p_2q_2-\theta^0-\theta^1-\theta^\infty_1-\theta^\infty_2)(1-2q_1)+(p_1+2t)q_2,\nonumber\\
M_{22}&=(1-2q_1)\left(\frac{1}{2}p_1+t\right)-2p_2q_2+4\theta^0+3\theta^1+4\theta^\infty_1+2\theta^\infty_2,\nonumber\\
t\frac{1}{v}\frac{dv}{dt}&=
p_1(1-2q_1)-2tq_1-2p_2q_2+2p_2q_1(q_1-1)+t+4\theta^0+3\theta^1+4\theta^\infty_1+2\theta^\infty_2.
\end{align}

\medskip
\noindent
\underline{Spectral type: (2)(11),22,22}

The Riemann scheme is given by
\[
\left(
\begin{array}{ccc}
  x=0 & x=1 & x=\infty \\
\begin{array}{c} 0 \\ 0 \\ \theta^0 \\ \theta^0 \end{array} &
\begin{array}{c} 0 \\ 0 \\ \theta^1 \\ \theta^1 \end{array} &
\overbrace{\begin{array}{cc}
     0   &   \theta^\infty_2 \\
     0   &   \theta^\infty_3 \\
     t & \theta^\infty_1 \\
     t & \theta^\infty_1
           \end{array}}
\end{array}
\right) ,
\]
and the Fuchs-Hukuhara relation is written as
$2\theta^0 +2\theta^1 +2\theta_1^\infty +\theta_2^\infty
+\theta_3^\infty =0$.

The Lax pair is expressed as
\begin{equation}
\left\{
\begin{aligned}
\frac{\partial Y}{\partial x}&=
\left(A_\infty+\frac{A_0^{(0)}}{x}+\frac{A_1^{(0)}}{x-1}
\right)Y, \\
\frac{\partial Y}{\partial t}&=(-E_2\otimes I_2x+B_1)Y ,
\end{aligned}
\right.
\end{equation}
where
\begin{align*}
A_{\xi}&=
({\rm diag}(1,v)\oplus U)^{-1}\hat{A}_{\xi}({\rm diag}(1,v)\oplus
 U),\\
A_\infty&=
\begin{pmatrix}
O_2 & O_2  \\
O_2 & -t I_2
\end{pmatrix},\quad
\hat{A}_0^{(0)}=
\begin{pmatrix}
QP+\theta^0+\theta^\infty_1 \\
tI_2
\end{pmatrix}
\left(I_2-Q,\ \frac{1}{t}\{(Q-I_2)QP+(\theta^0+\theta^\infty_1)Q-\theta^\infty_1\}\right),\\
\hat{A}_1^{(0)}&=
\begin{pmatrix}
(QP+\theta^0+\theta^\infty_1)(Q-I_2)-\Theta \\
tQ
\end{pmatrix}
\left(I_2,\ \frac{1}{t}
\{(QP+\theta^0+\theta^1+\theta^\infty_1+\Theta)Q^{-1}-QP-\theta^0-\theta^\infty_1\}\right).
\end{align*}
Furthermore
\begin{align*}
B_1&=
({\rm diag}(1,v)\oplus U)^{-1}
\begin{pmatrix}
O_2 & \frac{[\hat{A}_0^{(0)}+\hat{A}_1^{(0)}]_{12}}{t}\\
\frac{[\hat{A}_0^{(0)}+\hat{A}_1^{(0)}]_{21}}{t} & O_2
\end{pmatrix}
({\rm diag}(1,v)\oplus U) ,
\end{align*}
where $[\hat{A}_0^{(0)}+\hat{A}_1^{(0)}]_{ij}$
is the $(i,j)$-block of the matrix $\hat{A}_0^{(0)}+\hat{A}_1^{(0)}$.
Here, $Q$, $P$, and $\Theta$ are
\begin{align*}
Q=
\begin{pmatrix}
q_1 & 1\\
-q_2 & q_1
\end{pmatrix},\quad
P=
\begin{pmatrix}
p_1/2 & -p_2\\
p_2 q_2-\theta^0-\theta^1-\theta^\infty_1-\theta^\infty_2 & p_1/2
\end{pmatrix},\quad
\Theta=
\begin{pmatrix}
\theta^{\infty}_2 & \\
 & \theta^{\infty}_3
\end{pmatrix}.
\end{align*}

The Hamiltonian is given by
\begin{equation}
tH^{\mathrm{Mat}}_{\rm V}\left({-\theta^0-\theta^1-\theta^\infty_1,\theta^0-\theta^1 \atop \theta^1, \theta^\infty_3};t;
{q_1, p_1 \atop q_2, p_2}\right)=
\mathrm{tr}[Q(Q-1)P(P+t)+(\theta^0-\theta^1)QP+\theta^1 P+(\theta^0+\theta^\infty_1)tQ].
\end{equation}

The gauge parameters $U$ and $v$ satisfy
\begin{align}
t\frac{dU}{dt}&=
\begin{pmatrix}
tq_1-\theta^\infty_1+\theta^\infty_2+1 & t\\
-tq_2 & 2p_1q_1+3tq_1-p_1+2p_2q_2-2q_1p_2(q_1-1)-\eta-t
\end{pmatrix}U, \\
t\frac{1}{v}\frac{dv}{dt}&=(2q_1-1)p_1-2q_1p_2(q_1-1)+2tq_1+2p_2q_2-t+\theta^0-\theta^1, 
\quad \eta=\theta^0+3\theta^1+3\theta^\infty_1+\theta^\infty_2-1.
\end{align}

\subsubsection*{Singularity pattern: 3+1}
\noindent
\underline{Spectral type: ((2))((2)),211}

The Riemann scheme is given by
\[
\left(
\begin{array}{cc}
  x=0 & x=\infty \\
\overbrace{\begin{array}{ccc}
   0 & 0  & 0 \\
   0 & 0  & 0 \\
  -1 & t  & \theta^0 \\
  -1 & t & \theta^0
        \end{array}}
& \begin{array}{c} \theta^{\infty}_1 \\ \theta^{\infty}_1 \\ \theta^{\infty}_2 \\ \theta^{\infty}_3 \end{array}
\end{array}
\right) ,
\]
and the Fuchs-Hukuhara relation is written as
$2\theta^0 +2\theta_1^\infty +\theta_2^\infty +\theta_3^\infty =0$.

The Lax pair is expressed as
\begin{equation}
\left\{
\begin{aligned}
\frac{\partial Y}{\partial x}&=
\left(
\frac{A_0^{(-2)}}{x^3}+\frac{A_0^{(-1)}}{x^2}+\frac{A_0^{(0)}}{x}
\right)Y ,\\
\frac{\partial Y}{\partial t}&=\frac{A_0^{(-2)}}{x}Y .
\end{aligned}
\right.
\end{equation}
Here
$A_0^{(-2)}$, $A_0^{(-1)}$, and $A_0^{(0)}$ are given as follows:
\begin{align*}
A_{\xi}^{(k)}&=
(U\oplus {\rm diag}(v,1))^{-1}\hat{A}_{\xi}^{(k)}(U\oplus {\rm diag}(v,1)),\\
\hat{A}_0^{(-2)}&=
\begin{pmatrix}
I_2 \\
-Z
\end{pmatrix}
(-I_2-Z,-I_2),\quad
\hat{A}_0^{(-1)}=
\begin{pmatrix}
PZ+Q+t & P \\
-ZPZ-QZ-ZQ-tZ-Q & -ZP-Q
\end{pmatrix},\\
\hat{A}_0^{(0)}&=-
\begin{pmatrix}
\theta^{\infty}_1I_2 & O \\
O & \Theta
\end{pmatrix},\quad
Z=(\theta^{\infty}_1-\Theta)^{-1}
[(P-Q-t)Q-\theta^0-\theta^{\infty}_1].
\end{align*}
Here, $Q$, $P$, and $\Theta$ are
\begin{align*}
Q=
\begin{pmatrix}
q_1 & 1\\
-q_2 & q_1
\end{pmatrix},
\quad
P=
\begin{pmatrix}
p_1/2 & -p_2\\
p_2q_2-\theta^0-\theta^\infty_1-\theta^\infty_2 & p_1/2
\end{pmatrix},\quad
\Theta=
\begin{pmatrix}
\theta^{\infty}_2 & \\
 & \theta^{\infty}_3
\end{pmatrix}.
\end{align*}

The Hamiltonian is given by
\begin{align}
&H^{\mathrm{Mat}}_{\rm IV}
\left({\theta^0+2\theta^\infty_1-1,-\theta^0-\theta^\infty_1,\theta^\infty_1-\theta^\infty_2-1};t;
  {q_1,p_1\atop q_2,p_2}\right)\\
&\quad=\mathrm{tr}[PQ(P-Q-t)-(\theta^0+\theta^{\infty}_1)P+(\theta^0+2\theta^{\infty}_1-1)Q].\nonumber
\end{align}

The gauge parameters $U$ and $v$ satisfy
\begin{align}
 \frac{dU}{dt}=
\begin{pmatrix}
-\frac{3}{2}p_1+2(p_2+2)q_1+2t & p_2+2 \\
-(p_2+2)q_2+\theta^0+\theta^\infty_1+\theta^\infty_2 & 2q_1+t-\frac{p_1}{2}
\end{pmatrix}U,\qquad
\frac{1}{v}\frac{dv}{dt}=2(p_2+1)q_1-p_1+t.
\end{align}

\medskip
\noindent
\underline{Spectral type: ((2))((11)),22}

The Riemann scheme is given by
\[
\left(
\begin{array}{cc}
  x=0 & x=\infty \\
\begin{array}{c} 0 \\ 0 \\ \theta^0 \\ \theta^0\end{array}
&\overbrace{\begin{array}{ccc}
   0 & 0  & \theta^\infty_2 \\
   0 & 0  & \theta^\infty_3 \\
   1 & -t & \theta^\infty_1 \\
   1 & -t & \theta^\infty_1
        	\end{array}}
\end{array}
\right) ,
\]
and the Fuchs-Hukuhara relation is written as
$2\theta^0 +2\theta_1^\infty +\theta_2^\infty +\theta_3^\infty =0$.

The Lax pair is expressed as
\begin{equation}
\left\{
\begin{aligned}
\frac{\partial Y}{\partial x}&=
\left(
A_\infty^{(-2)}x+A_\infty^{(-1)}+\frac{A_0^{(0)}}{x}
\right)Y ,\\
\frac{\partial Y}{\partial t}&=(E_2\otimes I_2x+B_1)Y ,
\end{aligned}
\right.
\end{equation}
where
\begin{align*}
A_{\xi}^{(k)}&=
({\rm diag}(1,v)\oplus U)^{-1}\hat{A}_{\xi}^{(k)}({\rm diag}(1,v)\oplus U),\\
\hat{A}_\infty^{(-2)}&=
\begin{pmatrix}
O_2 & O_2\\
O_2 & -I_2
\end{pmatrix},\quad
\hat{A}_\infty^{(-1)}=
\begin{pmatrix}
O_2 & PQ-\Theta\\
I_2 & tI_2
\end{pmatrix},\quad
\hat{A}_0^{(0)}=
\begin{pmatrix}
-P \\
I_2
\end{pmatrix}
(Q,\ QP+\theta^0I_2),\\
B_1&=
\begin{pmatrix}
O_2 & -[A_\infty^{(-1)}]_{12}\\
-[A_\infty^{(-1)}]_{21} & O_2
\end{pmatrix},\quad
E_2=
\begin{pmatrix}
0 &  \\
 & 1
\end{pmatrix}.
\end{align*}
Here, $Q$, $P$, and $\Theta$ are
\begin{align*}
Q=
\begin{pmatrix}
q_1 & 1 \\
-q_2 & q_1
\end{pmatrix},\quad
P=
\begin{pmatrix}
p_1/2 & -p_2 \\
p_2q_2-\theta^0-\theta^\infty_1-\theta^\infty_2 & p_1/2
\end{pmatrix},\quad
\Theta=
\begin{pmatrix}
\theta^\infty_2 & \\
 & \theta^\infty_3
\end{pmatrix}.
\end{align*}

The Hamiltonian is given by
\begin{equation}
H^{\mathrm{Mat}}_{\rm IV}\left(-\theta^0-\theta^\infty_1, \theta^0, \theta^\infty_3 ;t;{q_1, p_1\atop q_2, p_2}\right)
 =\mathrm{tr}[QP(P-Q-t)+\theta^0 P-(\theta^0+\theta^\infty_1)Q].
\end{equation}

The gauge parameters $U$ and $v$ satisfy
\begin{align}
 \frac{dU}{dt}=
\begin{pmatrix}
-q_1-t & -1\\
q_2 & p_1-(2p_2+3)q_1-2t
\end{pmatrix}U,\quad
\frac{1}{v}\frac{dv}{dt}=p_1-2(p_2+1)q_1-t.
\end{align}

\subsubsection*{Singularity pattern: 2+2}
\noindent
\underline{Spectral type: (2)(2),(2)(11)}

The Riemann scheme is given by
\[
\left(
\begin{array}{cc}
  x=0 & x=\infty \\
\overbrace{\begin{array}{cc}
     0 & 0 \\
     0 & 0 \\
     t & \theta^0 \\
     t & \theta^0
           \end{array}}
& 
\overbrace{\begin{array}{cc}
     1 & \theta^\infty_1 \\
     1 & \theta^\infty_1 \\
     0 & \theta^\infty_2 \\
     0 & \theta^\infty_3
           \end{array}} 
\end{array}
\right) ,
\]
and the Fuchs-Hukuhara relation is written as
$2\theta^0 +2\theta_1^\infty +\theta_2^\infty +\theta_3^\infty =0$.

The Lax pair is expressed as
\begin{equation}
\left\{
\begin{aligned}
\frac{\partial Y}{\partial x}&=
\left(
\frac{A_0^{(-1)}}{x^2}+\frac{A_0^{(0)}}{x}+A_{\infty}
\right)Y ,\\
\frac{\partial Y}{\partial
 t}&=-\frac{1}{x}\left(\frac{A_0^{(-1)}}{t}\right)Y .
\end{aligned}
\right.
\end{equation}
Here
$A_0^{(-1)}$, $A_0^{(0)}$, and $A_{\infty}$ are given as follows:
\begin{align*}
A_{\xi}^{(k)}&=
(U\oplus {\rm diag}(v,1))^{-1}\hat{A}_{\xi}^{(k)}(U\oplus {\rm diag}(v,1)),\\
\hat{A}_0^{(-1)}&=
\begin{pmatrix}
I_2 \\
P
\end{pmatrix}
(t(1-P),tI_2),\quad
\hat{A}_0^{(0)}=
\begin{pmatrix}
-\theta^\infty_1I_2 & -Q \\
-Z & -\Theta
\end{pmatrix},\quad
\hat{A}_{\infty}=
\begin{pmatrix}
-I_2 & O \\
O & O
\end{pmatrix},\\
Z&=(QP+\theta^0+2\theta^{\infty}_1)P-(QP+\theta^0+\theta^{\infty}_1).
\end{align*}
Here, $Q$, $P$, and $\Theta$ are
\begin{align*}
Q=
\begin{pmatrix}
q_1 & 1\\
-q_2 & q_1
\end{pmatrix},
\quad
P=
\begin{pmatrix}
p_1/2 & -p_2\\
p_2q_2-\theta^0-\theta^\infty_1-\theta^\infty_2 & p_1/2
\end{pmatrix},\quad
\Theta=
\begin{pmatrix}
\theta^{\infty}_2 & \\
 & \theta^{\infty}_3
\end{pmatrix}.
\end{align*}

The Hamiltonian is given by
\begin{align}
&tH^{\mathrm{Mat}}_{\mathrm{III}(D_6)}
\left({\theta^0+\theta^\infty_1,\theta^0+2\theta^\infty_1,-\theta^\infty_2};t;
  {q_1,p_1\atop q_2,p_2}\right)
=\mathrm{tr}[Q^2P^2-(Q^2-(\theta^0+2\theta^\infty_1)Q-t)P-(\theta^0+\theta^\infty_1)Q].
\end{align}

The gauge parameters $U$ and $v$ satisfy
\begin{align}
t\frac{dU}{dt}&=
\begin{pmatrix}
(2p_2q_1-p_1+1)q_1 & p_1-2p_2q_1-1 \\
2q_1(p_2q_2-\theta^0-\theta^\infty_1-\theta^\infty_2)-(p_1-1)q_2 & 2p_2q_2+(p_1-1)q_1+\theta^0+2\theta^\infty_1
\end{pmatrix}U,\\
t\frac{1}{v}\frac{dv}{dt}&=2p_2({q_1}^2-q_2)-2(p_1-1)q_1+\theta^0+2\theta^\infty_2.
\end{align}

\subsubsection*{Singularity pattern: 4}
\noindent
\underline{Spectral type: (((2)))(((11)))}

The Riemann scheme is given by
\[
\left(
\begin{array}{c}
 x=\infty \\
\overbrace{\begin{array}{cccc}
   0 &  0 & 0  & \theta^\infty_1\\
   0 &  0 & 0  & \theta^\infty_1\\
   -1 & 0 & -t & \theta^\infty_2\\
   -1 & 0 & -t & \theta^\infty_3
        	\end{array}}
\end{array}
\right) ,
\]
and the Fuchs-Hukuhara relation is written as
$2\theta_1^\infty +\theta_2^\infty +\theta_3^\infty =0$.

The Lax pair is expressed as
\begin{equation}
\left\{
\begin{aligned}
\frac{\partial Y}{\partial x}&=
\left(
A_\infty^{(-3)}x^2+A_\infty^{(-2)}x+A_\infty^{(-1)}
\right)Y ,\\
\frac{\partial Y}{\partial t}&=(A_\infty^{(-3)}x+B_1)Y ,
\end{aligned}
\right.
\end{equation}
where
\begin{align*}
A_\infty^{(k)}&=
(U\oplus {\rm diag}(v,1))^{-1}\hat{A}_\infty^{(k)}(U\oplus {\rm diag}(v,1)),\\
\hat{A}_\infty^{(-3)}&=
\begin{pmatrix}
O & O \\
O & I_2
\end{pmatrix},\quad
\hat{A}_\infty^{(-2)}=
\begin{pmatrix}
O    & I_2 \\
-P+Q^2+t & O  
\end{pmatrix},\quad
\hat{A}_\infty^{(-1)}=
\begin{pmatrix}
-P+Q^2+t & Q \\
(P-Q^2-t)Q-\Theta & P-Q^2
\end{pmatrix},\\
B_1&=
(U\oplus {\rm diag}(v,1))^{-1}
\begin{pmatrix}
Q & I_2 \\
-P+Q^2+t & O
\end{pmatrix}
(U\oplus {\rm diag}(v,1)).
\end{align*}
Here $Q$, $P$, and $\Theta$ are
\begin{align*}
Q=
\begin{pmatrix}
q_1 & 1 \\
-q_2 & q_1
\end{pmatrix},\quad
P=
\begin{pmatrix}
p_1/2 & -p_2 \\
p_2q_2-\theta^\infty_1-\theta^\infty_2 & p_1/2
\end{pmatrix},\quad
\Theta=
\begin{pmatrix}
\theta^\infty_2 &  \\
 & \theta^\infty_3
\end{pmatrix}.
\end{align*}

The Hamiltonian is given by
\begin{equation}
H^{\mathrm{Mat}}_{\rm II}\left({-\theta^\infty_1+1,\theta^\infty_3+1};t;
  {q_1,p_1\atop q_2,p_2}\right)=\mathrm{tr}[P^2-(Q^2+t)P+(\theta^\infty_1-1)Q].
\end{equation}

The gauge parameters $U$ and $v$ satisfy
\begin{align}
\frac{dU}{dt}=
\begin{pmatrix}
2(q_1+p_2) & 0\\
0 & 0
\end{pmatrix}U,\quad
\frac{1}{v}\frac{dv}{dt}=2(q_1+p_2).
\end{align}

\section{Derivation and Calculations}\label{sec:comment}
Thus far, we have listed results without showing the derivation.
Here, we explain with examples how we have derived them. 
In this section, we treat four topics.
In subsection 4.1, we see derivation of Hamiltonians in Fuchsian case.
In subsection 4.2, we see procedure of degeneration with an example.
In subsection 4.3, we review Laplace transform.
In subsection 4.4, we treat degeneration of Fuchsian equations with three singular points.

\subsection{How to derive Hamiltonians for the Painlev\'e-type equations}
Our way is to compute isomonodromic deformation of Fuchsian equations and derive all the other equations by the degeneration  process.
The isomonodromic deformation of Fuchsian equations is treated minutely in Part 1.  

However, we can also directly consider the isomonodromic deformation of non-Fuchsian equations.
Let us take a simple example to see the computation.

%


Consider a system of linear differential equations which has
regular singularities on $\mathbb{C}$ and Poincar\'e rank (at most) 1 singularity
at infinity:
\begin{equation}
\frac{\mathrm{d}Y}{\mathrm{d}x}=
\left( \sum_{i=1}^n\frac{A_i}{x-u_i}+A_\infty \right)Y.
\end{equation}
Here $A_i$'s are $m \times m$ matrices and $A_\infty=\mathrm{diag}(a_1,\ldots,a_m)$.

The ``monodromy data'' of the above system stay constant if and only if $Y$ satisfies
\begin{equation}\label{eq:deformation}
\left\{
\begin{aligned}
\frac{\partial Y}{\partial u_i}&=-\frac{A_i}{x-u_i}Y,\\
\frac{\partial Y}{\partial a_j}&=(E_jx+B_j)Y.
\end{aligned}
\right.
\end{equation}
Here $E_j=\mathrm{diag}(0,\ldots,1,\ldots,0)$. $B_j$ are given by
\begin{equation}
(B_j)_{kl}=
\left\{
\begin{aligned}
&\frac{-(A_\infty^{(0)})_{jl}}{a_j-a_l} \quad (k=j,\ l\neq j),\\
&\frac{(A_\infty^{(0)})_{kj}}{a_k-a_j} \quad (k\neq j,\ l=j),\\
&0 \quad (\text{the otherwise}),
\end{aligned}
\right.
\end{equation}
where $A_\infty^{(0)}=-\sum_{i=1}^n A_i$.
The compatibility condition for (\ref{eq:deformation})
is the following equations
\begin{equation}\label{eq:gsch}
\left\{
\begin{aligned}
\frac{\partial A_i}{\partial u_k}&=\frac{[A_k,A_i]}{u_k-u_i} \quad (k \neq i),\\
\frac{\partial A_i}{\partial u_i}&=-\sum_{k\neq i}\frac{[A_k,A_i]}{u_k-u_i}-[A_i,A_\infty],\\
\frac{\partial A_i}{\partial a_j}&=-[A_i,u_iE_j+B_j] \quad (j=1,\ldots,m).
\end{aligned}
\right.
\end{equation}

If $A_\infty=0$, by ignoring the terms involving $a_j$'s, we have
\begin{equation}
\left\{
\begin{aligned}
\frac{\partial A_i}{\partial u_k}&=\frac{[A_k,A_i]}{u_k-u_i} \quad (k \neq i),\\
\frac{\partial A_i}{\partial u_i}&=-\sum_{k\neq i}\frac{[A_k,A_i]}{u_k-u_i}.
\end{aligned}
\right.
\end{equation}
This system is called the Schlesinger system.

Equation (\ref{eq:gsch}) can be written in Hamiltonian form.
First define the Poisson structure over $M_m(\mathbb{C}^n)$ by
\begin{equation}
\{ (A_p)_{i,j},(A_q)_{k,l} \}=\delta_{p,q}(\delta_{i,l}(A_p)_{k,j}-\delta_{k,j}(A_p)_{i,l}).
\end{equation}
Next define Hamiltonians by
\begin{equation}\label{eq:hamiltonian_derivation_nonfuchsian}
\left\{
\begin{aligned}
H_{u_i}&=\sum_{k\neq i}\frac{\mathrm{tr}(A_iA_k)}{u_i-u_k}+\mathrm{tr}(A_iA_\infty),\\
H_{a_j}&=\sum_{i=1}^n\mathrm{tr}(u_iA_iE_j)-\frac{1}{2}\mathrm{tr}\left(B_jA_\infty^{(0)}\right).
\end{aligned}
\right.
\end{equation}
Then (\ref{eq:gsch}) is rewritten into
\begin{equation}
\left\{
\begin{aligned}
\frac{\partial A_i}{\partial u_k}=\{ A_i,H_{u_k} \},\\
\frac{\partial A_i}{\partial a_j}=\{ A_i,H_{a_j} \}.
\end{aligned}
\right.
\end{equation}
In general, the entries of the matrices $A_i$'s are not canonical variables.
To take the canonical variables, find matrices $B$ and $C$ such that
\[
\sum_{i=1}^n \frac{A_i}{x-u_i}=B(x I-T)^{-1}C
\]
for
\[
T=
\begin{pmatrix}
u_1 I_{\mathrm{rank}A_1} &        & \\
& \ddots & \\
&        & u_n I_{\mathrm{rank}A_n}
\end{pmatrix}.
\]
In terms of entries of $B$ and $C$, (\ref{eq:gsch}) is written in Hamiltonian form
\begin{align}
\left\{
\begin{aligned}
\frac{\partial b_{ij}}{\partial u_k}&=\frac{\partial H_{u_k}}{\partial c_{ji}},\\
\frac{\partial c_{ij}}{\partial u_k}&=-\frac{\partial H_{u_k}}{\partial b_{ji}},
\end{aligned}\right.\qquad
\left\{
\begin{aligned}
\frac{\partial b_{ij}}{\partial a_k}&=\frac{\partial H_{a_k}}{\partial c_{ji}},\\
\frac{\partial c_{ij}}{\partial a_k}&=-\frac{\partial H_{a_k}}{\partial b_{ji}}.
\end{aligned}\right.
\end{align}

As an example, let us see the linear equation of type $(1)(1)(1),21,21$,
and how to derive Hamiltonians and canonical variables.
In this case, by a suitable gauge transformation, we can normalize the Riemann scheme
as follows:
\[
\left(
\begin{array}{ccc}
x=0 & x=1 & x=\infty \\
\begin{array}{c} 0 \\ 0 \\ \theta^0 \end{array}&
\begin{array}{c} 0 \\ 0 \\ \theta^1 \end{array}
&
\overbrace{\begin{array}{cc}
	0 & \theta^{\infty}_1 \\
	-t_1 & \theta^{\infty}_2 \\
	-t_2 & \theta^{\infty}_3
	\end{array}}
\end{array}
\right),
\]
Here the characteristic exponents satisfy the Fuchs-Hukuhara relation
\[
\theta^0 +\theta^1 +\theta_1^\infty +\theta_2^\infty
+\theta_3^\infty =0
\]
We consider a parameterization of the system
\[
\frac{d}{dx}Y=A(x)Y,\quad
A(x)=\frac{A_0}{x}+\frac{A_1}{x-1}+A_\infty,\quad
A_\infty ={\rm diag}(0,t_1,t_2)
\]
with the above Riemann scheme.

Here we factorize $A_\xi$，$\xi =0,1$:
\[
A_0=\left(\begin{array}{c}
b_1^0 \\
b_2^0 \\
b_3^0
\end{array}\right) (c_1^0, c_2^0, c_3^0),\qquad
A_1=\left(\begin{array}{c}
b_1^1 \\
b_2^1\\
b_3^1
\end{array}\right) (c_1^1, c_2^1, c_3^1).
\]
When $b_1^0b_1^1$ is not equal to zero, $b_1^0$ and $b_1^1$ can be normalized to 1.
Then the symplectic form becomes
$$\omega =dc_2^0\wedge db_2^0+dc_3^0\wedge db_3^0+dc_2^1\wedge db_2^1+dc_3^1\wedge db_3^1.$$

Now we have parameterized the linear system
by $b_2^0, b_2^1, b_3^0, b_3^1, c_1^0, c_1^1, c_2^0, c_2^1, c_3^0$, and
$c_3^1$.
However they have some relations.
First, from the trace of $A_0$ and $A_1$, we have
\[
\theta^0=c_1^0+b_2^0c_2^0+b_3^0c_3^0,\quad
\theta^1=c_1^1+b_2^1c_2^1+b_3^1c_3^1.
\]
Moreover, from the relation $A_0+A_1+A^{(0)}_\infty =0$,
and the diagonal part of $A_\infty^{(0)}$
is $\text{diag}(\theta_1^\infty,\theta_2^\infty,\theta_3^\infty)$, we have
\[
c_1^0+c_1^1+\theta_1^\infty  =0,\quad
b_2^0c_2^0+b_2^1c_2^1+\theta_2^\infty =0,
\quad b_3^0c_3^0+b_3^1c_3^1+\theta^{\infty}_3 =0.
\]
By the last relation, the symplectic form reduces to
\begin{eqnarray*}
	\omega&=&\frac{dc_2^0}{c_2^0}\wedge d(b_2^0c_2^0)
	+\frac{dc_2^1}{c_2^1}\wedge d(b_2^1c_2^1)
	+\frac{dc_3^0}{c_3^0}\wedge d(b_3^0c_3^0)
	+\frac{dc_3^1}{c_3^1}\wedge d(b_3^1c_3^1)\\
	&=&\left(\frac{dc_2^1}{c_2^1}-\frac{dc_2^0}{c_2^0}\right)
	\wedge d(b_2^1c_2^1)
	+\left(\frac{dc_3^1}{c_3^1}-\frac{dc_3^0}{c_3^0}\right)
	\wedge d(b_3^1c_3^1)\\
	&=&d\left(\frac{c_2^1}{c_2^0}\right)\wedge d(b_2^1c_2^0)
	+d\left(\frac{c_3^1}{c_3^0}\right)\wedge d(b_3^1c_3^0).
\end{eqnarray*}
Here we define $\lambda_j$ and $\mu_j$, $j=1,2$ as
\begin{equation*}
\mu_1 =b_2^0c_2^1,\quad \lambda_1 =-\frac{c_2^0}{c_2^1},\quad
\mu_2 =b_3^0c_3^1,\quad \lambda_2 =-\frac{c_3^0}{c_3^1}.
\end{equation*}

The Hamiltonian $H_{t_1}$ is calculated from the
relation (\ref{eq:hamiltonian_derivation_nonfuchsian})
as
\begin{align*}
H_{t_1}=&\text{tr}(0\cdot A_0E_2+1\cdot A_1E_2)\\
& +\frac{1}{t_1}(A_0+A_1)_{12}(A_0+A_1)_{21}
+\frac{1}{t_1-t_2}(A_0+A_1)_{23}(A_0+A_1)_{32}.
\end{align*}
$H_{t_2}$ is also obtained in the same manner.
Finally the Hamiltonians are given as
\begin{align}
&t_1\tilde {H}^{2+1+1+1}_{\mathrm{Gar},t_{1}}
\left({\theta^0, \theta^1 \atop \theta^\infty_2,\theta^\infty_3};{t_1 \atop t_2};
{\lambda_1,\mu_1 \atop \lambda_2,\mu_2}\right)\\
&=
t_1\tilde{H}_{\rm V}\left({\theta^0, \theta^1+\theta^\infty_3 \atop \theta^\infty_2};t_1;\lambda_1,\mu_1\right)
+(1-\lambda_1)\mu_2\lambda_2(\mu_1-\mu_1\lambda_1+\theta^\infty_2)\nonumber\\
&\quad+\frac{t_1}{t_1-t_2}(\mu_1(\lambda_1-\lambda_2)-\theta^\infty_2)(\mu_2(\lambda_2-\lambda_1)-\theta^\infty_3),\nonumber\\
&t_2\tilde{H}^{2+1+1+1}_{\mathrm{Gar},t_{2}}
\left({\theta^0, \theta^1 \atop \theta^\infty_1,\theta^\infty_2,\theta^\infty_3};{t_1 \atop t_2};
{\lambda_1,\mu_1 \atop \lambda_2,\mu_2}\right)\\
&=
t_2\tilde{H}_{\rm V}\left({\theta^0, \theta^1+\theta^\infty_2 \atop \theta^\infty_3};t_2;\lambda_2,\mu_2\right)
+(1-\lambda_2)\mu_1\lambda_1(\mu_2-\mu_2\lambda_2+\theta^\infty_3)\nonumber\\
&\quad+\frac{t_2}{t_2-t_1}(\mu_1(\lambda_1-\lambda_2)-\theta^\infty_2)(\mu_2(\lambda_2-\lambda_1)-\theta^\infty_3).\nonumber
\end{align}
These are Hamiltonians (\ref{eq:Gar2111_3}) and (\ref{eq:Gar2111_4}).

The general theory for the construction of symplectic structure get more complicated for the general cases.
Still, it is not difficult to compute individual cases.
Especially, when the Lax pair is given, it is easy to derive the deformation equations themselves, since we merely have to compute the compatibility conditions.

Most systems can be derived from degeneration so that we do not need the direct method explained in this section.
However in some cases, the number of deformation parameters increases in deformation.
In such cases, merely degenerating is insufficient.
We explain such cases in the end of Section~\ref{subsec:deg}.

\subsection{Procedure of degeneration}\label{subsec:deg}
In this subsection, we explain confluence of singularities of linear differential equations
and the way how they induce the degeneration of Painlev\'e-type equations.

As an example, we treat the degeneration of $A_5^{(1)}$-type Fuji-Suzuki system to $A_5^{(1)}$-type
Noumi-Yamada system.
It corresponds to the confluence of singularities of $21,21,111,111$-type Fuchsian system
that leads to $(2)(1),111,111$-type equation.

Consider the following system of linear differential equations
\begin{align}
\frac{dY}{dx}&=A(x)Y,\\
A(x)&=\frac{A_0}{x}+\frac{A_1}{x-1}+\frac{A_t}{x-t},
\end{align}
whose Riemann scheme is given by
\[
 \left(\begin{array}{cccc}
  x=0 & x=1 & x=t & x=\infty \\
  0 & 0 & 0 & \theta^{\infty}_1 \\
  \theta_1^0 & 0 & 0 & \theta^{\infty}_2 \\
  \theta_2^0 & \theta^1 & \theta^t & \theta^{\infty}_3
       \end{array}\right).
\]
Thus its spectral type is $21,21,111,111$.

Now we put $t=1+\varepsilon \tilde{t}$ and take the limit $\varepsilon \to 0$
so that the singular point $x=t$ merges to $x=1$.
Suppose the coefficient matrix
\begin{equation}\label{eq:before}
 A(x)=\frac{A_0}{x}+\frac{-\varepsilon \tilde{t} A_1}{(x-1)(x-1-\varepsilon \tilde{t})}
+\frac{A_{1}+A_t}{x-1-\varepsilon \tilde{t}} 
\end{equation}
tends to
\begin{equation}\label{eq:after}
 \tilde{A}(x)=\frac{\tilde{A}_0^{(0)}}{x}+\frac{\tilde{A}_{1}^{(-1)}}{(x-1)^2}+\frac{\tilde{A}_{1}^{(0)}}{x-1}
\end{equation}
as $\varepsilon \to 0$.
Here the coefficient in (\ref{eq:after}) corresponds to a system of linear differential equations
with spectral type $(2)(1),111,111$, whose Riemann scheme is
\[
\left(\begin{array}{ccc}
  x=0 & x=1 & x=\infty \\
\begin{array}{c}0 \\ \tilde{\theta}^0_1 \\ \tilde{\theta}^0_2  \end{array}
& \overbrace{\begin{array}{cc}
     0 & 0 \\
     0 & 0 \\
     \tilde{t} & \tilde{\theta}^1
        \end{array}}\ 
& \begin{array}{c} \tilde{\theta}^{\infty}_1 \\ \tilde{\theta}^{\infty}_2  \\ \tilde{\theta}^{\infty}_3 \end{array} 
       \end{array}\right).
\]
By comparing (\ref{eq:before}) and (\ref{eq:after}), we obtain
\begin{align*}
&\lim_{\varepsilon \to 0}(-\varepsilon \tilde{t}\theta^1)=
\lim_{\varepsilon \to 0}\mathrm{tr}(-\varepsilon\tilde{t}A_1)=\mathrm{tr}(\tilde{A}_{1}^{(-1)})=\tilde{t},\\
&\lim_{\varepsilon \to 0}(\theta^1+\theta^{t})=\lim_{\varepsilon \to 0}\mathrm{tr}(A_{1}+A_{t})=
\mathrm{tr}(\tilde{A}_{1}^{(0)})=\tilde{\theta}^{1},\\
&\lim_{\varepsilon \to 0}(\theta^0_1+\theta^0_2)=
\lim_{\varepsilon \to 0}\mathrm{tr}(A_0)=\mathrm{tr}(\tilde{A}_{0}^{(0)})=\tilde{\theta}^0_1+\tilde{\theta}^0_2.
\end{align*}
Taking the above into account, we require the following relations:
\[\theta^{0}_{i}=\tilde{\theta}^{0}_{i} \quad(i=1,2), \quad \theta^{1}=-\varepsilon^{-1},
\quad \theta^{t}=\tilde{\theta}^{1}+\varepsilon^{-1},
\quad \theta^{\infty}_{j}=\tilde{\theta}^{\infty}_{j} \quad(j=1,2,3).
\]
Furthermore, we define a canonical transformation in the following manner:
\begin{align*}
&q_{1}=1+\varepsilon \tilde{t} \tilde{q}_{1},\quad p_{1}=\varepsilon^{-1}\tilde{t}^{-1}\tilde{p}_{1},
\quad q_{2}=1+\varepsilon \tilde{t} \tilde{q}_{2},\quad p_{2}=\varepsilon^{-1}\tilde{t}^{-1}\tilde{p}_{2},\\
&H_{\mathrm{FS}}^{A_{5}}=\varepsilon^{-1}(\tilde{H}+t^{-1}(\tilde{p}_{1}\tilde{q}_{1}+\tilde{p}_{2}\tilde{q}_{2})).
\end{align*}
It is easy to check that $\lim_{\varepsilon \to 0}\tilde{H}=H_{\mathrm{NY}}^{A_{5}}$.
In this way, we obtain $A_5^{(1)}$-type Noumi-Yamada system.
\begin{rem}\label{rem:transf_abbreviation}
Instead of using tilde, we often express such a transformation as
\begin{align*}
&\theta^1 \to -\varepsilon^{-1},\ \theta^t \to \theta^1+\varepsilon^{-1},\ t \to 1+\varepsilon t,\ 
H\to \varepsilon^{-1}(H+t^{-1}(p_1 q_1+p_2 q_2)),\\
&q_1 \to 1+\varepsilon tq_1,\ p_1 \to \varepsilon^{-1}t^{-1}p_1,\ 
q_2 \to 1+\varepsilon tq_2,\ p_2 \to \varepsilon^{-1}t^{-1}p_2.
\end{align*}
\qed
\end{rem}

Here we would like to make some important remarks. 
It is impossible to merge two singular points without changing
the number of accessory parameters unless one of the spectral types of singular points is
a refinement of the other.
Therefore, for example,
the linear equation with only one irregular singularity (i.e. corresponds to 4) does not appear
in the degeneration scheme of the $31,22,22,1111$ system since singularity of
type $22$ and $31$ cannot be merged.

If the singularity of type $22$ and $31$ should merge, then we have a linear equation with a singularity of Poincar\'e rank one and two regular singularity of type $22$ and $1111$.
The number of accessory parameters of this linear equation is not four according to the following formula.  

Recall the formula of rigidity index for non-Fuchsian case \cite{Tk}.
We write the spectral type of non-Fuchsian equation
\[
[p_{0}^{(r_{0})}, p_{0}^{(r_{0}-1)},\dots,p_{0}^{(0)}],
[p_{1}^{(r_{1})},\dots,p_{1}^{(0)}],
[p_{n}^{(r_{n})},\dots,p_{n}^{(0)}]
\]
where $p_{i}^{(s)}$ is a partition
\[
p_{i}^{(s)}=m_{i_{1}}^{(s)}\dots \, m_{i,l_{is}}^{(s)}.
\]
Then the rigidity index is given by
\[
{\rm idx}=-\left(\sum_{i=0}^{n}r_{i}+n-1\right)m^{2}
+\sum_{i=0}^{n}\sum_{s=0}^{r_{i}}\sum_{j=1}^{l_{is}}\left(m_{ij}^{(s)}\right)^{2}.
\]


We should also mention that the degenerate linear equations may admit more freedoms of deformation than the source equation.
In our degeneration scheme, such phenomena can be only seen in the degeneration of Fuji-Suzuki system.

For non-Fuchsian case, the deformation parameters come from both configuration of singularities and 
canonical forms at irregular singularities except the residue parts.
The residue parts of the canonical forms are  constants and they are not deformation parameters. 
If we exclude parameters that can be normalized to $0$ or $1$ or other constants by automorphisms of $\mathbb{P}^{1}$, that is to say, linear fractional transformations and gauge transformations, then the rest is deformation parameters.

Within the degenerations of Fuji-Suzuki system, four Painlev\'e-type equations associated with linear equations of spectral types $(1)(1)(1),21,21$ and $((1)(1))((1)),21$ and $(2)(1),(1)(1)(1)$ and $(((1)(1)))(((1)))$ are expressed as partial differential equations with two independent variables.
For example, when the spectral type is $(1)(1)(1),21,21$, we can assume that the singular points are at $0,1,\infty$, and that the Riemann scheme can be reduced by gauge transformation of scalar matrices as
\[
\left(
\begin{array}{ccc}
  x=0 & x=1 & x=\infty \\
\begin{array}{c} 0 \\ 0 \\ \theta^0  \end{array} &
\begin{array}{c} 0 \\ 0 \\ \theta^1  \end{array} &
\overbrace{\begin{array}{cc}
0 & \theta^{\infty}_{1}\\
-t_{1} & \theta^{\infty}_{2}\\
-t_{2} & \theta^{\infty}_{3}
      \end{array}}
\end{array}
\right).
\]
In this case, $t_{1}$ and $t_{2}$ become deformation parameters.
Its deformation equation is expressed by $H_{\mathrm{Gar}}^{2+1+1+1}$.
This is the reason why degenerate Garnier systems appear in the degeneration scheme of Fuji-Suzuki system.
On the other hand, we can also consider usual degeneration from $21,21,111,111$:
\begin{align*}
&x\to (1-\varepsilon tx)/(1-\varepsilon t),\ t \to 1/(1-\varepsilon t),\\
&\theta^0_i \to \eta_i \varepsilon^{-1},\
\theta^t \to \theta^0,\ \theta^1 \to \theta^1,\\
&\theta^\infty_1\to -\eta_1 \varepsilon^{-1}+\theta^\infty_1,\ 
\theta^\infty_2\to -\eta_2 \varepsilon^{-1}+\theta^\infty_2,\ 
\theta^\infty_3\to \theta^\infty_3, \\
&q_1 \to \frac{(1-\varepsilon t)q_1}{(1-\varepsilon t)q_1-1},\ 
p_1 \to -(1-\varepsilon t)q_1(p_1q_1-\theta_2^\infty)
\left( 1-\frac{1}{(1-\varepsilon t)q_1}\right) ,\\
&q_2 \to \frac{(1-\varepsilon t)q_2}{(1-\varepsilon t)q_2-1},\ 
p_2 \to 
-(1-\varepsilon t)q_2(p_2q_2-\theta_3^\infty)
\left( 1-\frac{1}{(1-\varepsilon t)q_2}\right) .
\end{align*}
This yields Hamiltonian
\begin{align*}
t&H
\left({\theta^{\infty}_1,\theta^{\infty}_2,\theta^{\infty}_3
 \atop\theta^0, \theta^1};\eta_1,\eta_2;t;
  {q_1,p_1\atop q_2,p_2}\right)\\
=&tH_{\rm V}\left({\theta^0, \theta^\infty_2,
\theta^1}
;(\eta_2-\eta_1)t;q_1,p_1\right) 
+tH_{\rm V}\left({\theta^0, \theta^\infty_3,
\theta^1}
;-\eta_1t;q_2,p_2\right) \\
&+p_1p_2(2q_1q_2-q_1-q_2).
\end{align*}
This Hamiltonian is a restriction of $H_{\mathrm{Gar}}^{2+1+1+1}$.
Similar stories are true for other three systems, and they correspond to $H_{\mathrm{Gar}}^{3+1+1}, H_{\mathrm{Gar}}^{\frac{3}{2}+1+1+1}$, and $H_{\mathrm{Gar}}^{\frac{5}{2}+1+1}$, respectively.

In order to obtain the Hamiltonians for $(1)(1)(1),21,21$ and $((1)(1))((1)),21$, and $(2)(1),(1)(1)(1)$, and $(((1)(1)))(((1)))$, we calculate the Hamiltonian of $(1)(1)(1),21,21$ directly, and obtain the others by degeneration. 
Here the deformation equation of $(1)(1)(1),21,21$ in
Hamiltonian form
is calculated in Section~\ref{subsec:deg} by the use of the technique described in an appendix of \cite{JMMS}.

\subsection{Laplace transform}

As we have already mentioned, different linear equations sometimes induce the same Painlev\'e-type equations.
In this sense, there is no one-to-one correspondence of linear equations in the scheme to 4-dimensional Painlev\'e-type equations.
However, as far as we know, cases when the same nonlinear equations appear are only those cases when corresponding linear equations transforms one another by the Laplace transforms.
It is an interesting question whether or not linear equations with the same Painlev\'e-type equation always transform one another by the Laplace transform or some other transformations.
Concerning the correspondences of linear equations by the Laplace transform, see \cite{Hrd,B2}.

Within the degeneration scheme, there are seven deformation equations that have two different linear equations.

Let us see the Laplace transform in a case when there are one irregular singularity of Poincar\'e rank 1 and some regular singularities.
In this case, the linear equation can be expressed as
\begin{equation}
 \frac{d}{dx}Y=\left[ Q\left(xI_l-T\right)^{-1}P+S\right] Y,
\end{equation}
where $l=\sum_{i=1}^n{\rm rank}A_i$, $Q$ is $m\times l$ matrix, and $P$ is $l\times m$ matrix.
Matrices $T$ and $S$ are diagonal.

This equation can be rewritten as
\begin{equation}
 \left(\begin{array}{cc}
  \frac{d}{dx}-S & -Q \\
  -P & xI_l-T
       \end{array}\right)
 \left(\begin{array}{c}
  Y \\
  Z
       \end{array}\right) =0.
\end{equation}
Let us apply the Laplace transform $(x, d/dx)\mapsto (-d/d\xi , \xi)$.
We express the transform of a dependent variable using $\hat{~}$.
If we regard this equation as equation of $\hat{Z}$, the equation reads
\begin{equation}
  \frac{d}{d\xi}\hat{Z}=-\left[ P\left(\xi I_m-S\right)^{-1}Q+T\right]
  \hat{Z},
\end{equation}
which is similar to the original one.
Such calculation tells us four correspondences:
\begin{align*}
 &(1)(1),11,11,11\leftrightarrow (1)(1)(1),21,21\qquad
 (2)(1),111,111\leftrightarrow (11)(11),31,21\\
 &(2)(2),31,1111\leftrightarrow (111)(1),22,22\qquad
 (2)(2),22,211\leftrightarrow (2)(11),22,22.
\end{align*}

\begin{eg}
	Here we see the first one of the above correspondences $(1)(1),11,11,11\leftrightarrow (1)(1)(1),21,21$.
	The isomonodromic deformation equations of these equations are both $H_{\rm Gar}^{2+1+1+1}$.
	
	Define the matrices $S$, $T$, $P$, and $Q$ by
	\begin{align*}
	&S=
	\begin{pmatrix}
	0 & & \\
	& t_1 & \\
	& & t_2
	\end{pmatrix}, \quad
	T=
	\begin{pmatrix}
	0 & \\
	& 1
	\end{pmatrix}, \\
	&Q=
	\begin{pmatrix}
	1 & 1 \\
	\mu_1 & \mu_1\lambda_1-\theta^\infty_2 \\
	\mu_2 & \mu_2\lambda_2-\theta^\infty_3
	\end{pmatrix}, \quad
	P=
	\begin{pmatrix}
	\mu_1\lambda_1+\mu_2\lambda_2+\theta^0 & -\lambda_1 & -\lambda_2 \\
	-\mu_1\lambda_1-\mu_2\lambda_2+\theta^1+\theta^\infty_2+\theta^\infty_3 & 1 & 1
	\end{pmatrix}.
	\end{align*}
	Then the $(1)(1)(1), 21, 21$ system can be written as follows:
	\begin{align*}
	\frac{dY}{dx}&=\left(
	\frac{\hat{A}_0^{(0)}}{x}+\frac{\hat{A}_1^{(0)}}{x-1}+A_\infty \right)Y \\
	&=\left[ Q(x I_2-T)^{-1}P+S \right]Y.
	\end{align*}
	Applying the Laplace transformation, we obtain
	\begin{align*}
	\frac{d\hat{Z}}{d\xi}&=[-P(\xi I_3-S)^{-1}Q-T]\hat{Z}\\
	&=\left[
		\begin{pmatrix}
		-\mu_1\lambda_1-\mu_2\lambda_2-\theta^0 \\
		\mu_1\lambda_1+\mu_2\lambda_2-\theta^1-\theta^\infty_2-\theta^\infty_3
		\end{pmatrix}
		\begin{pmatrix}
		1 & 1
		\end{pmatrix}\frac{1}{\xi}\right.\\
	&\qquad\left.
	+\begin{pmatrix}
	  \lambda_1 \\
	  -1
	  \end{pmatrix}
	  \begin{pmatrix}
	  \mu_1 & \mu_1\lambda_1-\theta^\infty_2
	  \end{pmatrix}\frac{1}{\xi-t_1}
	+\begin{pmatrix}
	  \lambda_2 \\
	  -1
	  \end{pmatrix}
	  \begin{pmatrix}
	  \mu_2 & \mu_2\lambda_2-\theta^\infty_3
	  \end{pmatrix}\frac{1}{\xi-t_2}
	-T
	\right]\hat{Z}.
	\end{align*}
	The spectral type of the system is $(1)(1), 11, 11, 11$.
	
	Moreover, by a change of variables
	\[
	\xi \to t_1\xi, \quad
	\hat{Z} \to 
	\begin{pmatrix} -\mu_1\lambda_1-\mu_2\lambda_2-\theta^0 & 0 \\
	0 & \mu_1\lambda_1+\mu_2\lambda_2-\theta^1-\theta^\infty_2-\theta^\infty_3
	\end{pmatrix}^{-1}\hat{Z},
	\]
	we have
	\begin{align}\label{eq:Lap_(1)(1)(1),21,21}
	\frac{d\hat{Z}}{d\xi}=
	\left(
	\frac{\tilde{A}_0}{\xi}+\frac{\tilde{A}_1}{\xi-1}+\frac{\tilde{A}_t}{\xi-t_2/t_1}-t_1T
	\right)\hat{Z},
	\end{align}
	where
	\begin{align*}
	\tilde{A}_0&=
	\begin{pmatrix}
	1 \\
	1
	\end{pmatrix}
	\begin{pmatrix}
	-\mu_1\lambda_1-\mu_2\lambda_2-\theta^0 &
	\mu_1\lambda_1+\mu_2\lambda_2-\theta^1-\theta^\infty_2-\theta^\infty_3
	\end{pmatrix},\\
	\tilde{A}_1&=
	\begin{pmatrix}
	\frac{\lambda_1}{-\mu_1\lambda_1-\mu_2\lambda_2-\theta^0} \\
	-\frac{1}{\mu_1\lambda_1+\mu_2\lambda_2-\theta^1-\theta^\infty_2-\theta^\infty_3}
	\end{pmatrix}
	\begin{pmatrix}
	\mu_1(-\mu_1\lambda_1-\mu_2\lambda_2-\theta^0) &
	(\mu_1\lambda_1-\theta^\infty_2)(\mu_1\lambda_1+\mu_2\lambda_2-\theta^1-\theta^\infty_2-\theta^\infty_3)
	\end{pmatrix}\\
	&=
	\begin{pmatrix}
	1 \\
	\frac{\mu_1\lambda_1+\mu_2\lambda_2+\theta^0}{\lambda_1(\mu_1\lambda_1+\mu_2\lambda_2-\theta^1-\theta^\infty_2-\theta^\infty_3)}
	\end{pmatrix}
	\begin{pmatrix}
	\mu_1\lambda_1 &
	\lambda_1(\mu_1\lambda_1-\theta^\infty_2)
	\frac{\mu_1\lambda_1+\mu_2\lambda_2-\theta^1-\theta^\infty_2-\theta^\infty_3}{-\mu_1\lambda_1-\mu_2\lambda_2-\theta^0}
	\end{pmatrix},\\
	\tilde{A}_t&=
	\begin{pmatrix}
	\frac{\lambda_2}{-\mu_1\lambda_1-\mu_2\lambda_2-\theta^0} \\
	-\frac{1}{\mu_1\lambda_1+\mu_2\lambda_2-\theta^1-\theta^\infty_2-\theta^\infty_3}
	\end{pmatrix}
	\begin{pmatrix}
	\mu_2(-\mu_1\lambda_1-\mu_2\lambda_2-\theta^0) &
	(\mu_2\lambda_2-\theta^\infty_3)(\mu_1\lambda_1+\mu_2\lambda_2-\theta^1-\theta^\infty_2-\theta^\infty_3)
	\end{pmatrix}\\
	&=
	\begin{pmatrix}
	1 \\
	\frac{\mu_1\lambda_1+\mu_2\lambda_2+\theta^0}{\lambda_2(\mu_1\lambda_1+\mu_2\lambda_2-\theta^1-\theta^\infty_2-\theta^\infty_3)}
	\end{pmatrix}
	\begin{pmatrix}
	\mu_2\lambda_2 &
	\lambda_2(\mu_2\lambda_2-\theta^\infty_3)
	\frac{\mu_1\lambda_1+\mu_2\lambda_2-\theta^1-\theta^\infty_2-\theta^\infty_3}{-\mu_1\lambda_1-\mu_2\lambda_2-\theta^0}
	\end{pmatrix}.
	\end{align*}
	Comparing them with (3.5), we obtain the following canonical transformation:
	\begin{align*}
	&\lambda_1 \to 
	\frac{\mu_1\lambda_1+\mu_2\lambda_2+\theta^0}{\lambda_1(\mu_1\lambda_1+\mu_2\lambda_2-\theta^1-\theta^\infty_2-\theta^\infty_3)}, \quad
	\mu_1 \to \lambda_1(\mu_1\lambda_1-\theta^\infty_2)
	\frac{\mu_1\lambda_1+\mu_2\lambda_2-\theta^1-\theta^\infty_2-\theta^\infty_3}{-\mu_1\lambda_1-\mu_2\lambda_2-\theta^0},\\
	&\lambda_2 \to 
	\frac{\mu_1\lambda_1+\mu_2\lambda_2+\theta^0}{\lambda_2(\mu_1\lambda_1+\mu_2\lambda_2-\theta^1-\theta^\infty_2-\theta^\infty_3)}, \quad
	\mu_2 \to \lambda_2(\mu_2\lambda_2-\theta^\infty_3)
	\frac{\mu_1\lambda_1+\mu_2\lambda_2-\theta^1-\theta^\infty_2-\theta^\infty_3}{-\mu_1\lambda_1-\mu_2\lambda_2-\theta^0},\\
	&t_1 \to -t_1, \quad t_2 \to -t_2, \quad H_{t_1} \to -H_{t_1}, \quad H_{t_2} \to -H_{t_2},
	\end{align*}
	which maps the Hamiltonians (\ref{eq:Gar2111_1}) and (\ref{eq:Gar2111_2}) to (\ref{eq:Gar2111_3}) and (\ref{eq:Gar2111_4}).
	\qed
\end{eg}

\begin{rem}
If we put $S={\rm diag}(0,1,t)$, $T=\left(\begin{array}{cc}
					0 & 1 \\
					0 & 0
					     \end{array}\right)$
in the above correspondence, we obtain a transformation between the two by two system with the singularity pattern $\frac{3}{2}+1+1+1$ and the three by three system with the spectral type $(2)(1),(1)(1)(1)$.
The deformation equations of these linear equations are both the Garnier system of type $\frac{3}{2}+1+1+1$.
\qed
\end{rem}

When the Poincar\'e rank is 2, the calculation becomes more complicated.          
For an equation
\begin{equation}
 \frac{d}{dx}Y=\left[ Q\left(xI_l-T\right)^{-1}P+S_0+S_1x\right] Y,
\end{equation}
let us put $S_1={\rm diag}(a_1,\ldots ,a_k, 0,\ldots , 0)$, and $\tilde{S_1}={\rm diag}(a_1,\ldots , a_k)$.
We also assume that
$S_0=\left(\begin{array}{cc}
	     S_0^{11} & S_0^{12} \\
	     S_0^{21} & S_0^{22}
		  \end{array}\right)$,
$Q=\left(\begin{array}{c}
    Q_1 \\
    Q_2
	 \end{array}\right)$,
$P=(P_1, P_2)$,
$Y=\left(\begin{array}{c}
    Y_1 \\
    Y_2
	 \end{array}\right).$
Here, $S_0^{11}$ is $k\times k$ matrix, $Q_1$ is $k\times l $matrix, $Q_2$ is
$(m-k)\times l$ matrix, $P_1$ is $l\times k$ matrix, and $P_2$ is $l\times (m-k)$ matrix.

We can rewrite the equation as
\begin{equation}
 \left(\begin{array}{ccc}
  \frac{d}{dx}-S_0^{11}-\tilde{S_1}x & -S_0^{12} & -Q_1 \\
  -S_0^{21} & \frac{d}{dx}-S_0^{22} & -Q_2 \\
  -P_1 & -P_2 & xI_l-T
       \end{array}\right)
 \left(\begin{array}{c}
  Y_1 \\
  Y_2 \\
  Z
       \end{array}\right) =0.
\end{equation}
If we perform the Laplace transform $(x, d/dx)\mapsto (-d/d\xi , \xi)$,  we can eliminate
\begin{equation}
 \hat{Y_2}=(\xi I_{m-k}-S_0^{22})^{-1}(S_0^{21}\hat{Y_1}+Q_2\hat{Z}).
\end{equation}
Thus, the equation becomes
\begin{align}
 &\frac{d}{d\xi}\left(\begin{array}{c}
		 \hat{Y_1}
		 \label{eq:Lap-pr2}
		 \\
		 \hat{Z}
		      \end{array}\right) \\
 &=\left[\left(\begin{array}{c}
	  \tilde{S_1}^{-1}S_0^{12}\\
	  -P_2
	       \end{array}\right)
 (\xi I_{m-k}-S_0^{22})^{-1}
 (S_0^{21}, Q_2)
 +\left(\begin{array}{cc}
   \tilde{S_1}^{-1}S_0^{11} & \tilde{S_1}^{-1}Q_1 \\
   -P_1 &-T
	\end{array}\right)
 -\xi 
 \begin{pmatrix}
 \tilde{S}_1^{-1} & \\
 & O 
 \end{pmatrix}
 \right]
 \left(\begin{array}{c}
		 \hat{Y_1}\\
		 \hat{Z}
		      \end{array}\right). \nonumber
\end{align}
Such calculation tells us three correspondences:
\begin{align*}
 &((1))((1)),11,11\leftrightarrow ((1)(1))((1)),21\qquad
 ((11))((1)),111\leftrightarrow ((11))((11)),31\\
 &((2))((2)),211\leftrightarrow ((11))((2)),22.
\end{align*}

\begin{eg}
	Here we illustrate the correspondence between $((1)(1))((1)), 21$ and $((1))((1)), 11, 11$.
	
	First, applying a gauge transformation $Y \to \exp(\frac12 x^2-t_2x)Y$ to
	the $((1)(1))((1)), 21$-system (\ref{eq:LaxGar3+1+1}), we assume the Riemann scheme of the $((1)(1))((1)), 21$-system
	have the following form:
	\[
	\left(
	\begin{array}{cc}
	x=0 & x=\infty \\
	\begin{array}{c} 0 \\ 0 \\ \theta^0 \end{array}
	&\overbrace{\begin{array}{ccc}
		1  & -t_2 & \theta^\infty_1 \\
		0  & t_1-t_2 & \theta^\infty_2 \\
		0  & 0 & \theta^\infty_3
		\end{array}}
	\end{array}
	\right).
	\]
	Then the $((1)(1))((1)), 21$-system is written as follows:
	\begin{align*}
	\frac{dY}{dx}&=
	\left( \hat{A}_\infty^{(-2)}x+\hat{A}_\infty^{(-1)}+\frac{\hat{A}_0^{(0)}}{x} \right)Y\\
	&=\left[ Q(x-T)^{-1}P+S_0+S_1x \right]Y
	\end{align*}
	where
	\begin{align*}
	&\hat{A}_\infty^{(-2)}=
	\begin{pmatrix}
	-1 &  &  \\
	& 0 &  \\
	&  & 0
	\end{pmatrix}=S_1,\quad
	\hat{A}_\infty^{(-1)}=
	\begin{pmatrix}
	t_2                       & -1   & -1   \\
	-p_1q_1+\theta^\infty_2 & t_2-t_1 &  0   \\
	-p_2q_2+\theta^\infty_3 &  0   & 0
	\end{pmatrix}=S_0,\quad
	\hat{A}_0^{(0)}=QP,\\
	&Q=\begin{pmatrix}
	1 \\
	p_1 \\
	p_2
	\end{pmatrix},\quad
	P=
	\begin{pmatrix}
	p_1q_1+p_2q_2+\theta^0 & -q_1 & -q_2
	\end{pmatrix}, \quad T=0.
	\end{align*}
	
	In this case, (\ref{eq:Lap-pr2}) reads
		\begin{align*}
		\frac{d}{d\xi}\left(\begin{array}{c}
		\hat{Y_1}\\
		\hat{Z}
		\end{array}\right)
		&=\left[
		\begin{pmatrix}
		1 & 1 \\
		q_1 & q_2
		\end{pmatrix}
		\left( \xi I_2-\begin{pmatrix} 
		t_2-t_1 & 0 \\
		0 & 0
		\end{pmatrix} \right)^{-1}
		\begin{pmatrix}
		-p_1q_1+\theta^\infty_2 & p_1 \\
		-p_2q_2+\theta^\infty_3 & p_2
		\end{pmatrix}\right.\\
		&\qquad\left.+\begin{pmatrix}
		-t_2 & -1 \\
		-p_1q_1-p_2q_2-\theta^0 & 0
		\end{pmatrix}
		+
		\begin{pmatrix}
		1 & 0 \\
		0 & 0 
		\end{pmatrix}\xi 
		\right]
		\left(\begin{array}{c}
		\hat{Y_1}\\
		\hat{Z}
		\end{array}\right) \nonumber \\
		&=\left[
		\begin{pmatrix}
		1 \\
		q_1
		\end{pmatrix}
		\begin{pmatrix}
		-p_1q_1+\theta^\infty_2 & p_1
		\end{pmatrix}\frac{1}{\xi-(t_2-t_1)}
		+
		\begin{pmatrix}
		1 \\
		q_2
		\end{pmatrix}
		\begin{pmatrix}
		-p_2q_2+\theta^\infty_3 & p_2
		\end{pmatrix}\frac{1}{\xi}\right.\\
		&\qquad\left.+\begin{pmatrix}
		-t_2 & -1 \\
		-p_1q_1-p_2q_2-\theta^0 & 0
		\end{pmatrix}
		+
		\begin{pmatrix}
		1 & 0 \\
		0 & 0 
		\end{pmatrix}\xi 
		\right]
		\left(\begin{array}{c}
		\hat{Y_1}\\
		\hat{Z}
		\end{array}\right).
		\nonumber
		\end{align*}
	By changing the independent variable
	\[
	\begin{pmatrix}
	\hat{Y_1}\\
	\hat{Z}
	\end{pmatrix}
	\to
	\begin{pmatrix}
	0 & 1 \\
	1 & 0
	\end{pmatrix}
	\begin{pmatrix}
	\hat{Y_1}\\
	\hat{Z}
	\end{pmatrix}
	=
	\begin{pmatrix}
	\hat{Z}\\
	\hat{Y_1}
	\end{pmatrix},
	\]
	we have
	\begin{align*}
	\frac{d}{d\xi}\left(\begin{array}{c}
	\hat{Z}\\
	\hat{Y_1}
	\end{array}\right)
	&=\left[
	\begin{pmatrix}
	q_2 \\
	1
	\end{pmatrix}
	\begin{pmatrix}
	p_2 & -p_2q_2+\theta^\infty_3
	\end{pmatrix}\frac{1}{\xi}
	+
	\begin{pmatrix}
	q_1 \\
	1
	\end{pmatrix}
	\begin{pmatrix}
	p_1 & -p_1q_1+\theta^\infty_2
	\end{pmatrix}\frac{1}{\xi-(t_2-t_1)}\right.\\
	&\qquad\left.+\begin{pmatrix}
	0 & -p_1q_1-p_2q_2-\theta^0 \\
	-1 & -t_2
	\end{pmatrix}
	+\begin{pmatrix}
	0 & 0 \\
	0 & 1 
	\end{pmatrix}\xi 
	\right]
	\left(\begin{array}{c}
	\hat{Z}\\
	\hat{Y_1}
	\end{array}\right).
	\nonumber
	\end{align*}
	This coincides with the $((1))((1)), 11, 11$-system.
	\qed
\end{eg}

\begin{rem}
If we put
$T={\rm diag}(0,1)$, $S_1=\left(\begin{array}{cc}
					0 & 1 \\
					0 & 0
					     \end{array}\right),$
we can obtain a transformation 
between the two by two system with the singularity pattern $\frac{5}{2}+1+1$ and the three by three system with the spectral type $(((1)(1)))(((1)))$.
The deformation equations of these linear equations are both the Garnier system of type $\frac{5}{2}+1+1$.
\qed
\end{rem}

\subsection{Fuchsian equations with three singular points}
So far, we have omitted Fuchsian equations with only three singular points for they do not admit deformation.
There are 9 equations whose number of singularities is three.
Do not they admit deformation when they are degenerated?
The answer is yes; some of Fuchsian equations admit deformation when they are degenerated.
However, we can see that all the equations derived from the 9 equations by confluences of singular points can be transformed to one of the equations in the scheme by the Laplace transform.
Thus, if we put equations of ramified-type aside, then all the equations are included in the degeneration scheme.

According to Oshima's classification of Fuchsian equations with four accessory parameters, those with only three singular points have the following spectral types:  

211,1111,1111\qquad 221,221,11111\qquad 32,11111,11111\qquad
222,222,2211\qquad 33,2211,111111

44,2222,22211\qquad 44,332,11111111\qquad 55,3331,22222\qquad
66,444,2222211.

Let us see the degeneration schemes arising from the Fuchsian equations with three singularities.

\medskip
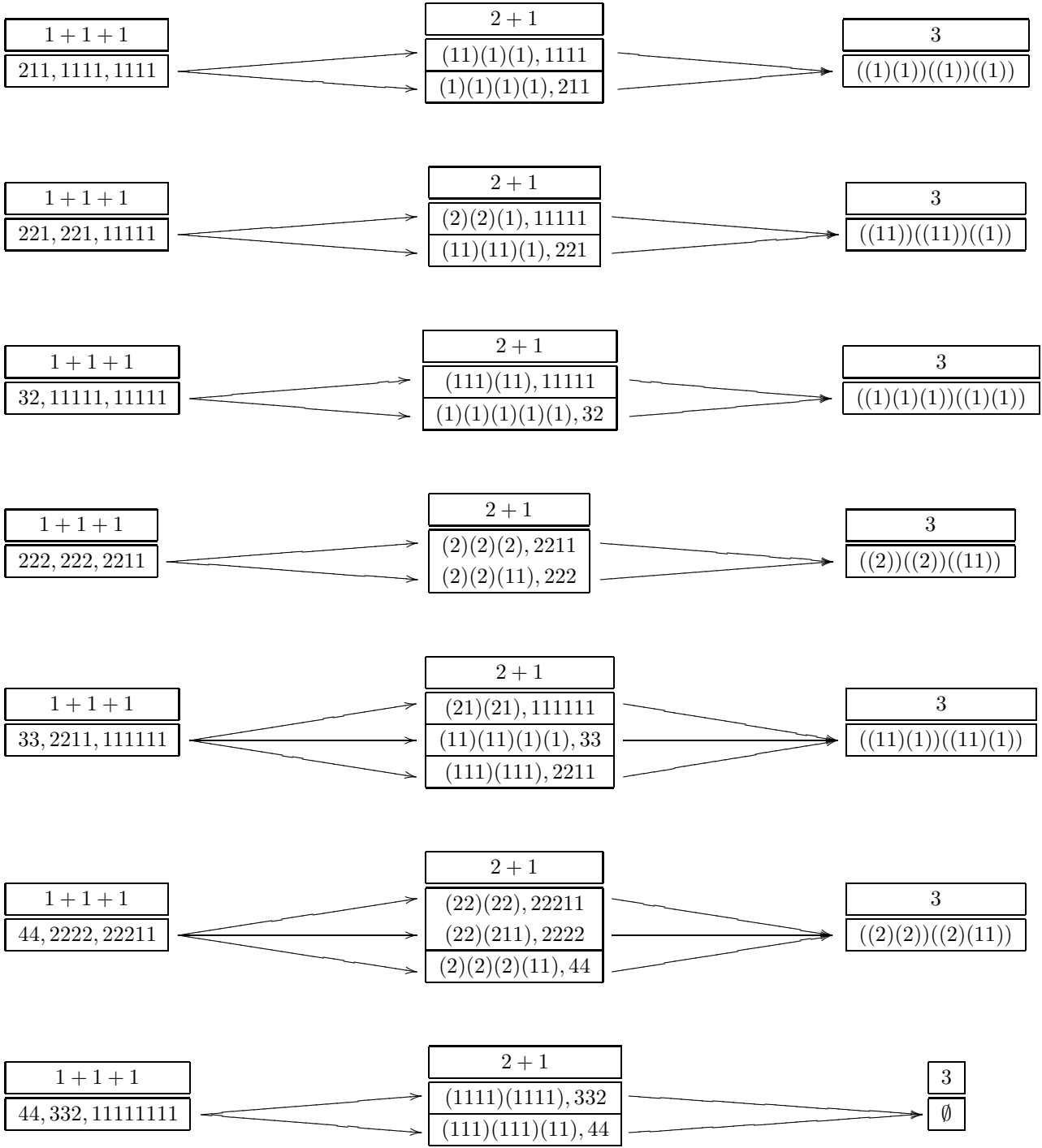
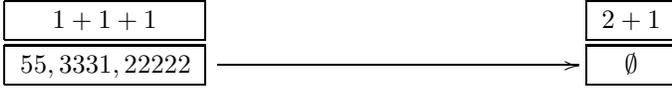
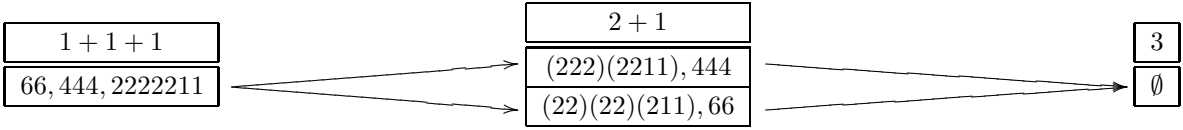
\begin{xy}
{(0,0) *{\begin{tabular}{|c|}
\hline
$1+1+1$\\
\hline\hline
$211, 1111, 1111$\\
\hline
\end{tabular}
}},
{\ar (15,-3);(54,0)},
{\ar (15,-3);(54,-6)},
{(70,0) *{\begin{tabular}{|c|}
\hline
$2+1$\\
\hline\hline
$(11)(1)(1), 1111$\\
\hline
$(1)(1)(1)(1), 211$\\
\hline
\end{tabular}}},
{\ar (87,0);(122,-3)},
{\ar (87,-6);(122,-3)},
{(139,0) *{\begin{tabular}{|c|}
\hline
3 \\
\hline\hline
$((1)(1))((1))((1))$\\
\hline
\end{tabular}}},
\end{xy}

\vspace{10mm}

\begin{xy}
{(0,0) *{\begin{tabular}{|c|}
\hline
$1+1+1$\\
\hline\hline
$221, 221, 11111$\\
\hline
\end{tabular}
}},
{\ar (15,-3);(54,0)},
{\ar (15,-3);(54,-6)},
{(70,0) *{\begin{tabular}{|c|}
\hline
$2+1$\\
\hline\hline
$(2)(2)(1), 11111$\\
\hline
$(11)(11)(1), 221$\\
\hline
\end{tabular}}},
{\ar (86,0);(123,-3)},
{\ar (86,-6);(123,-3)},
{(139,0) *{\begin{tabular}{|c|}
\hline
3 \\
\hline\hline
$((11))((11))((1))$\\
\hline
\end{tabular}}},
\end{xy}

\vspace{10mm}

\begin{xy}
{(0,0) *{\begin{tabular}{|c|}
\hline
$1+1+1$\\
\hline\hline
$32, 11111, 11111$\\
\hline
\end{tabular}
}},
{\ar (16,-3);(52,0)},
{\ar (16,-3);(52,-6)},
{(70,0) *{\begin{tabular}{|c|}
\hline
$2+1$\\
\hline\hline
$(111)(11), 11111$\\
\hline
$(1)(1)(1)(1)(1), 32$\\
\hline
\end{tabular}}},
{\ar (88,0);(121,-3)},
{\ar (88,-6);(121,-3)},
{(139,0) *{\begin{tabular}{|c|}
\hline
3 \\
\hline\hline
$((1)(1)(1))((1)(1))$\\
\hline
\end{tabular}}},
\end{xy}

\vspace{10mm}

\begin{xy}
{(0,0) *{\begin{tabular}{|c|}
\hline
$1+1+1$\\
\hline\hline
$222, 222, 2211$\\
\hline
\end{tabular}
}},
{\ar (14,-3);(55,0)},
{\ar (14,-3);(55,-6)},
{(70,0) *{\begin{tabular}{|c|}
\hline
$2+1$\\
\hline\hline
$(2)(2)(2), 2211$\\
$(2)(2)(11), 222$\\
\hline
\end{tabular}}},
{\ar (85,0);(123,-3)},
{\ar (85,-6);(123,-3)},
{(139,0) *{\begin{tabular}{|c|}
\hline
3 \\
\hline\hline
$((2))((2))((11))$\\
\hline
\end{tabular}}},
\end{xy}

\vspace{10mm}

\begin{xy}
{(0,0) *{\begin{tabular}{|c|}
\hline
$1+1+1$\\
\hline\hline
$33, 2211, 111111$\\
\hline
\end{tabular}
}},
{\ar (16,-3);(53,3)},
{\ar (16,-3);(53,-3)},
{\ar (16,-3);(53,-9)},
{(70,0) *{\begin{tabular}{|c|}
\hline
$2+1$\\
\hline\hline
$(21)(21), 111111$\\
\hline
$(11)(11)(1)(1), 33$\\
\hline
$(111)(111), 2211$\\
\hline
\end{tabular}}},
{\ar (87,3);(122,-3)},
{\ar (87,-3);(122,-3)},
{\ar (87,-9);(122,-3)},
{(139,0) *{\begin{tabular}{|c|}
\hline
3 \\
\hline\hline
$((11)(1))((11)(1))$\\
\hline
\end{tabular}}},
\end{xy}

\vspace{10mm}

\begin{xy}
{(0,0) *{\begin{tabular}{|c|}
\hline
$1+1+1$\\
\hline\hline
$44, 2222, 22211$\\
\hline
\end{tabular}
}},
{\ar (15,-3);(54,3)},
{\ar (15,-3);(54,-3)},
{\ar (15,-3);(54,-9)},
{(70,0) *{\begin{tabular}{|c|}
\hline
$2+1$\\
\hline\hline
$(22)(22), 22211$\\
$(22)(211), 2222$\\
\hline
$(2)(2)(2)(11), 44$\\
\hline
\end{tabular}}},
{\ar (86,3);(122,-3)},
{\ar (86,-3);(122,-3)},
{\ar (86,-9);(122,-3)},
{(139,0) *{\begin{tabular}{|c|}
\hline
3 \\
\hline\hline
$((2)(2))((2)(11))$\\
\hline
\end{tabular}}},
\end{xy}

\vspace{10mm}

\begin{xy}
{(0,0) *{\begin{tabular}{|c|}
\hline
$1+1+1$\\
\hline\hline
$44, 332, 11111111$\\
\hline
\end{tabular}
}},
{\ar (17,-3);(53,0)},
{\ar (17,-3);(53,-6)},
{(70,0) *{\begin{tabular}{|c|}
\hline
$2+1$\\
\hline\hline
$(1111)(1111), 332$\\
\hline
$(111)(111)(11), 44$\\
\hline
\end{tabular}}},
{\ar (87,0);(135,-3)},
{\ar (87,-6);(135,-3)},
{(139,0) *{\begin{tabular}{|c|}
\hline
3 \\
\hline\hline
$\emptyset$\\
\hline
\end{tabular}}},
\end{xy}

\vspace{10mm}

\begin{xy}
{(0,0) *{\begin{tabular}{|c|}
\hline
$1+1+1$\\
\hline\hline
$55, 3331, 22222$\\
\hline
\end{tabular}
}},
{\ar (15,-3);(63,-3)},
{(70,0) *{\begin{tabular}{|c|}
\hline
$2+1$\\
\hline\hline
$\emptyset$\\
\hline
\end{tabular}}},
\end{xy}

\vspace{10mm}

\begin{xy}
{(0,0) *{\begin{tabular}{|c|}
\hline
$1+1+1$\\
\hline\hline
$66, 444, 2222211$\\
\hline
\end{tabular}
}},
{\ar (16,-3);(54,0)},
{\ar (16,-3);(54,-6)},
{(70,0) *{\begin{tabular}{|c|}
\hline
$2+1$\\
\hline\hline
$(222)(2211), 444$\\
\hline
$(22)(22)(211), 66$\\
\hline
\end{tabular}}},
{\ar (87,0);(135,-3)},
{\ar (87,-6);(135,-3)},
{(139,0) *{\begin{tabular}{|c|}
\hline
3 \\
\hline\hline
$\emptyset$\\
\hline
\end{tabular}}},
\end{xy}
\medskip					

If we consider the equations with the singularity pattern $2+1$, then we obtain 18 equations:

(11)(1)(1),1111\qquad (1)(1)(1)(1),211\qquad
(2)(2)(1),11111\qquad (11)(11)(1),221

(111)(11),11111\qquad (1)(1)(1)(1)(1),32\qquad
(2)(2)(2),2211\qquad (2)(2)(11),222

(21)(21),111111\quad (111)(111),2211\quad (11)(11)(1)(1),33\quad
(22)(22),22211\quad (22)(211),2222\quad (2)(2)(2)(11),44

(1111)(1111),332\qquad (111)(111)(11),44\qquad
(222)(2211),444\qquad (22)(22)(211),66.

\noindent
Among these 18 equations, the following 7 equations do not admit deformations:

(111)(11),11111\qquad(21)(21),111111\qquad (111)(111),2211\qquad
(22)(22),22211\qquad (22)(211),2222

(1111)(1111),332\qquad (222)(2211),444.

\noindent
If we consider the Laplace transform of 11 equations which admit deformations, then the leading terms at irregular singularities become scalar matrix, since the original ones have only one regular singularity.
Thus, we can eliminate these leading terms at irregular singularities by gauge transformation of rmamatio.
As a result, we obtain Fuchsian equations, and the Fuchsian equations are already classified.

For example, we can see a correspondence as below:
\begin{equation*}
 (11)(1)(1),1111\leftrightarrow (111),111,21,21\sim 21,21,111,111.
\end{equation*}
Similarly, we can find correspondences for the rest 10 equations:
\begin{align*}
 &(1)(1)(1)(1),211\leftrightarrow 11,11,11,11,11\qquad
 (2)(2)(1),11111\leftrightarrow 31,22,22,1111 \\
 &(11)(11)(1),221\leftrightarrow 21,21,111,111\qquad
 (1)(1)(1)(1)(1),32\leftrightarrow 11,11,11,11,11 \\
 &(2)(2)(2),2211\leftrightarrow 22,22,22,211\qquad
 (2)(2)(11),222\leftrightarrow 22,22,22,211 \\
 &(11)(11)(1)(1),33\leftrightarrow 21,21,111,111\qquad
 (2)(2)(2)(11),44\leftrightarrow 22,22,22,211 \\
 &(111)(111)(11),44\leftrightarrow 211,1111,1111\qquad
 (22)(22)(211),66\leftrightarrow 222,222,2211.
\end{align*}
For the last two equations, their original deformations are trivially solved, and corresponding Fuchsian equations do not admit deformations.
Other equations correspond to one of the four Fuchsian equations corresponding to 4-dimensional Painlev\'e-type equations.

If we consider the confluences of three regular singular points to one points for each nine Fuchsian equations with only three singularities in Oshima's list, we obtain 6 equations:

((1)(1))((1))((1))\quad ((11))((11))((1))\quad
((1)(1)(1))((1)(1))\quad ((2))((2))((11))

((11)(1))((11)(1))\quad ((2)(2))((2)(11)).

\noindent
They are same as one of the degenerate systems of the four 4-dimensional Painlev\'e-type equations.

In this paper, we depend on the idea that all Painlev\'e-type equations of unramified non-Fuchsian systems with 4-dimensional phase space are derived from the four Painlev\'e-type equations of the Fuchsian equations by degeneration processes.
If we degenerate equations after applying middle convolution, do not we come up with new equations that are not derived from the original equations?
In fact, when the dimension of phase space is equal to or greater than six, such a case happens.
However, when the dimension of phase space is two or four, it is proved that such a case does not happen \cite{HO}.

\appendix
\section{Data on degenerations}
In this appendix, we give explicitly the canonical transformations with $\varepsilon$
that link two Hamiltonians in each degeneration.
The way of degenerations of Hamiltonians are explained in subsection \ref{subsec:deg}.
For notation see Remark~\ref{rem:transf_abbreviation}.
It is not necessary to change the data that do not appear in the table below.
Note that the terms in Hamiltonian that do not contain the canonical variables $p_i$ and $q_i$ are irrelevent to Hamiltonian system, thus we add or subtract such terms as needed.

Here we omit the following five degenerations
\begin{align*}
&21,21,111,111\to(1)(1)(1),21,21,\quad (11)(1),21,111\to ((1)(1))((1)),21, \\
&(2)(1),111,111 \to (2)(1),(1)(1)(1),\quad(11)(1),(11)(1)\to (((1)(1)))(((1))),\\
& ((11))((1)),111\to (((1)(1)))(((1)))
\end{align*}
since the number of deformation parameters increases.

\subsection{Garnier system}
\noindent
{\bf 1+1+1+1+1 $\to$ 2+1+1+1}
\begin{align*}
&\theta^{t_1} \to -\varepsilon^{-1},\ \theta^{t_2} \to \theta^t,\ 
\theta^\infty_1 \to \theta^\infty_2+\varepsilon^{-1},\ 
\theta^\infty_2 \to \theta^\infty_1,\\ 
& t_1 \to (\varepsilon t_1)^{-1},\ t_2 \to t_2/t_1,\
H_{t_1} \to -\varepsilon {t_1}^2 H_{t_1}-\varepsilon t_1 t_2 H_{t_2},\ 
H_{t_2} \to t_1H_{t_2},\\
&q_1 \to \frac{p_1(1-q_1)+p_2(1-q_2)+\theta^1+\theta^t+\theta^\infty_1}
{\varepsilon t_1(q_1-1)(p_1(q_1-1)-\theta^1)},\ 
p_1 \to \varepsilon t_1(q_1-1)(p_1(q_1-1)-\theta^1),\\
&q_2 \to -\frac{t_2(q_2-1)(p_2(q_2-1)-\theta^t)}{t_1(q_1-1)(p_1(q_1-1)-\theta^1)},\ 
p_2 \to -\frac{t_1(q_1-1)(p_1(q_1-1)-\theta^1)}{t_2(q_2-1)}.
\end{align*}
\noindent
{\bf 2+1+1+1 $\to$ 3+1+1}
\begin{align*}
&\theta^0 \to \theta^\infty_2+\varepsilon^{-2},\ \theta^t \to \theta^0,\ 
\theta^\infty_1 \to \theta^\infty_2,\ 
\theta^\infty_2 \to \theta^\infty_1-\theta^\infty_2-\varepsilon^{-2},\\
&t_1 \to -\varepsilon^{-1}t_1-\varepsilon^{-2},\ 
t_2 \to -\varepsilon^{-1}t_2-\varepsilon^{-2},\ H_{t_1} \to -\varepsilon H_{t_1},\ H_{t_2} \to -\varepsilon H_{t_2},\\
&q_1 \to \frac{1}{1-\varepsilon q_1},\ 
p_1 \to (1-\varepsilon q_1)(\varepsilon^{-1}p_1(1-\varepsilon q_1)+\theta^1),\\
&q_2 \to \frac{1}{1-\varepsilon q_2},\ 
p_2 \to (1-\varepsilon q_2)(\varepsilon^{-1}p_2(1-\varepsilon q_2)+\theta^0).
\end{align*}
\noindent
{\bf 2+1+1+1 $\to$ 2+2+1}
\begin{align*}
&\theta^0 \to -\varepsilon^{-1},\ \theta^t \to \theta^0+\varepsilon^{-1},\
 t_2 \to \varepsilon t_2,\ H_{t_1} \to H_{t_1},\ H_{t_2} \to \varepsilon^{-1}H_{t_2}, \\
&q_1 \to q_1,\ p_1 \to p_1,\ 
q_2 \to -\frac{1}{\varepsilon q_2},\ p_2 \to q_2(\varepsilon p_2q_2-\varepsilon \theta^0-1).
\end{align*}
\noindent
{\bf 2+2+1 $\to$ 3+2}
\begin{align*}
&\theta^1 \to -\varepsilon^{-2},\ \theta^\infty_2 \to \theta^\infty_2+\varepsilon^{-2},\ 
t_1 \to \varepsilon^{-1}t_2-\varepsilon^{-2},\ t_2 \to \varepsilon^{-1}t_1,\ H_{t_1} \to \varepsilon H_{t_2},\ H_{t_2} \to \varepsilon H_{t_1},\\
&q_1 \to \frac{\varepsilon(p_1q_1-p_2q_2-\theta^\infty_1)}
{\varepsilon(p_1q_1-p_2q_2-\theta^\infty_1)(1+\varepsilon q_2)+q_2},\ 
p_1 \to \frac{1+\varepsilon q_2}
{\varepsilon^2 q_2}\{\varepsilon(p_1q_1-p_2q_2-\theta^\infty_1)(1+\varepsilon q_2)+q_2\},\\
&q_2 \to -t_1\frac{p_1}{q_2},\ p_2 \to \frac{q_1q_2}{t_1}+1.
\end{align*}
\noindent
{\bf 2+2+1 $\to$ 4+1}
\begin{align*}
&\theta^0 \to -2\varepsilon^{-3},\ \theta^1 \to \theta^0,\ 
\theta^\infty_2 \to \theta^\infty_2+2\varepsilon^{-3},\ 
t_1 \to -\varepsilon^{-2}t_2+\varepsilon^{-3},\ 
t_2 \to \varepsilon^{-4}t_1+\varepsilon^{-6},\\
&H_{t_1} \to -\varepsilon^2 H_{t_2},\ H_{t_2} \to \varepsilon^4 H_{t_1},\ q_1 \to \frac{q_2}{q_2+\varepsilon p_1},\ 
p_1 \to \frac{(q_2+\varepsilon p_1)(p_2(q_2+\varepsilon p_1)-\theta^0)}{\varepsilon p_1},\\
&q_2 \to \frac{-p_1(1+\varepsilon q_1)+\varepsilon(p_2q_2+\theta^\infty_1)}
{\varepsilon^3p_1},\ 
p_2 \to 1-\varepsilon^2 p_1.
\end{align*}
\noindent
{\bf 3+1+1 $\to$ 3+2}
\begin{align*}
&\theta^0 \to -\varepsilon^{-1},\ \theta^1 \to \theta^0+\varepsilon^{-1},\ 
t_1 \to t_2-\varepsilon t_1,\ H_{t_1} \to -\varepsilon^{-1}H_{t_1},\ H_{t_2} \to \varepsilon^{-1}H_{t_1}+H_{t_2},\\
&q_1 \to \varepsilon t_1 p_1+q_2,\ p_1 \to \frac{-q_1}{\varepsilon t_1},\ 
p_2 \to p_2+\frac{q_1}{\varepsilon t_1}.
\end{align*}
\noindent
{\bf 3+1+1 $\to$ 4+1}
\begin{align*}
&\theta^1 \to \varepsilon^{-6},\ \theta^\infty_2 \to \theta^\infty_2-\varepsilon^{-6},\ 
t_1 \to \varepsilon t_1-2\varepsilon^{-3},\ t_2 \to -\varepsilon^{-1}t_2-\varepsilon^{-3},\\
&H_{t_1} \to \varepsilon^{-1} H_{t_1},\ H_{t_2} \to -\varepsilon H_{t_2},
\ q_1 \to -\varepsilon p_1,\ p_1 \to \varepsilon^{-1}q_1-\varepsilon^{-3},\ 
q_2 \to \varepsilon^{-1}q_2,\ p_2 \to \varepsilon p_2.
\end{align*}
\noindent
{\bf 3+2 $\to$ 5}
\begin{align*}
&\theta^0 \to 3\varepsilon^{-4},\ \theta^\infty_2 \to \theta^\infty_2-3\varepsilon^{-4},\ 
t_1 \to \varepsilon^{-3}t_1+\varepsilon^{-4}t_2+\varepsilon^{-6},\ 
t_2 \to t_2-3\varepsilon^{-2},\\
&H_{t_1} \to \varepsilon^3 H_{t_1},\ H_{t_2} \to -\varepsilon^{-1}H_{t_1}+H_{t_2}, 
\ q_1 \to \varepsilon^{-3}q_1+\varepsilon^{-4},\ p_1 \to \varepsilon^3 p_1-\varepsilon^2 q_2,\ 
p_2 \to p_2-\varepsilon^{-1}q_1-2\varepsilon^{-2}.
\end{align*}
\noindent
{\bf 4+1 $\to$ 5}
\begin{align*}
&\theta^0 \to -\varepsilon^{-12},\ \theta^\infty_1 \to \theta^\infty_1,\
t_1 \to -\varepsilon t_1+\varepsilon^{-2}t_2+\frac34\varepsilon^{-8},\
t_2 \to -\varepsilon^2t_2+\frac32\varepsilon^{-4},\\
&H_{t_1} \to -\varepsilon^{-1} H_{t_1},\ H_{t_2} \to -\varepsilon^{-5}H_{t_1}-\varepsilon^{-2}H_{t_2}-\varepsilon^{-2}q_2,
\ q_1 \to -\varepsilon^{-1}q_1-\varepsilon^{-4}/2,\ p_1 \to
 -\varepsilon^{-2} q_2,\\
&q_2 \to
\varepsilon^2\frac{q_2+\varepsilon^3(q_1q_2-p_1)-\varepsilon^6\theta^\infty_1}{1+\varepsilon^3q_1+\varepsilon^6(p_2-t_2)},\
p_2 \to \varepsilon^{-2}(p_2-t_2)+\varepsilon^{-5}q_1+\varepsilon^{-8}.
\end{align*}

\subsection{Fuji-Suzuki system}
\noindent
{\bf 1+1+1+1$\to$ 2+1+1}

\noindent
$21,21,111,111 \to (2)(1),111,111$
\begin{align*}
&\theta^1 \to -\varepsilon^{-1},\ \theta^t \to \theta^1+\varepsilon^{-1},\ t \to 1+\varepsilon t,\ 
H\to \varepsilon^{-1}(H+t^{-1}(p_1 q_1+p_2 q_2)),\\
&q_1 \to 1+\varepsilon tq_1,\ p_1 \to \varepsilon^{-1}t^{-1}p_1,\ 
q_2 \to 1+\varepsilon tq_2,\ p_2 \to \varepsilon^{-1}t^{-1}p_2.
\end{align*}
\noindent
$21,21,111,111 \to (11)(1),21,111$
\begin{align*}
&\theta^0_1 \to \theta^0_2-\varepsilon^{-1},\ 
\theta^0_2 \to \theta^0_1,\ \theta^t \to \varepsilon^{-1},\ t \to \varepsilon t,\ H\to \varepsilon^{-1}H,\\
&q_1 \to 1/q_1,\ p_1 \to -q_1(p_1q_1-\theta^0_1-\theta^{\infty}_2),\ 
q_2 \to 1/q_2,\ p_2 \to -q_2(p_2q_2-\theta^{\infty}_3).
\end{align*}
\noindent
{\bf 2+1+1$\to$ 3+1}

\noindent
$(2)(1),111,111 \to ((11))((1)),111$
\begin{align*}
&\theta^0_1 \to \theta^0_2+\varepsilon^{-2},\ \theta^0_2 \to \theta^0_1,\ \theta^1 \to -\varepsilon^{-2},
t \to -\varepsilon^{-1}t+\varepsilon^{-2},\ 
H \to -\varepsilon H,\\
&q_1 \to \varepsilon q_1,\ p_1 \to \varepsilon^{-1}p_1,\ 
q_2 \to \varepsilon q_2,\ p_2 \to \varepsilon^{-1}p_2.
\end{align*}
\noindent
$(11)(1),21,111\to((11))((1)),111$
\begin{align*}
&\theta^1 \to \varepsilon^{-2},\ 
\theta^0_2 \to \theta^0_2-\varepsilon^{-2},\ 
t \to -\varepsilon^{-1}t-\varepsilon^{-2},\ 
H \to -\varepsilon H,\\
&q_1 \to \frac{1}{1-\varepsilon q_2},
\ p_1 \to (1-\varepsilon q_2)\left(\theta^0_1+\theta^\infty_2+\varepsilon^{-1}p_2-p_2q_2\right),\\
&q_2 \to \frac{1}{1-\varepsilon q_1},
\ p_2 \to (1-\varepsilon q_1)\left(\theta^\infty_3+\varepsilon^{-1}p_1-p_1q_1\right).
\end{align*}
\noindent
$(1)(1)(1),21,21 \to ((1)(1))((1)),21$
\begin{align*}
&\theta^1 \to \varepsilon^{-2},\ \theta^\infty_1 \to \theta^\infty_1-\varepsilon^{-2},\ 
t_1 \to -\varepsilon^{-1}t_1-\varepsilon^{-2},\ t_2 \to -\varepsilon^{-1}t_2-\varepsilon^{-2},\\ 
&H_{t_1} \to -\varepsilon H_{t_1},\ H_{t_2} \to -\varepsilon H_{t_2},\\
&q_1 \to \frac{1}{1-\varepsilon q_1},\ 
p_1 \to (1-\varepsilon q_1)\left(\frac{p_1}{\varepsilon}(1-\varepsilon q_1)+\theta^\infty_2\right),\ 
q_2 \to \frac{1}{1-\varepsilon q_2},
\ p_2 \to (1-\varepsilon q_2)\left(\frac{p_2}{\varepsilon}(1-\varepsilon q_2)+\theta^\infty_3\right).
\end{align*}
\noindent
{\bf 2+1+1$\to$ 2+2}

\noindent
$(11)(1),21,111\to(11)(1),(11)(1)$
\begin{align*}
&\theta^1 \to \varepsilon^{-1},\ \theta^\infty_1 \to \theta^\infty_3-\varepsilon^{-1},\ 
\theta^\infty_3 \to \theta^\infty_1,\ t \to \varepsilon t,\ H \to \varepsilon^{-1}\left(H-\frac{p_1q_1+p_2q_2}{t}\right),\\
&q_1 \to -\frac{q_1}{\varepsilon t},\ p_1 \to -\varepsilon t p_1,\ 
q_2 \to -\frac{q_2}{\varepsilon t},\ p_2 \to -\varepsilon t p_2.
\end{align*}
\noindent
$(1)(1)(1),21,21\to(2)(1),(1)(1)(1)$
\begin{align*}
&\theta^0 \to -\varepsilon^{-1},\ 
\theta^1 \to \theta^0+\varepsilon^{-1},\ 
t_1 \to \varepsilon t_1,\ t_2 \to \varepsilon t_2,\ 
H_{t_1} \to \varepsilon^{-1}H_{t_1},\ H_{t_2} \to \varepsilon^{-1}H_{t_2},\\
&q_1 \to -\frac{1}{\varepsilon q_1},\  
p_1 \to \varepsilon q_1(p_1q_1-\theta^\infty_2),\ 
q_2 \to -\frac{1}{\varepsilon q_2},\  
p_2 \to \varepsilon q_2(p_2q_2-\theta^\infty_3).
\end{align*}
\noindent
{\bf 2+2$\to$ 4}

\noindent
$(2)(1),(1)(1)(1) \to (((1)(1)))(((1)))$
\begin{align*}
&\theta^\infty_1 \to 2\varepsilon^{-2},\ \theta^\infty_2 \to
 -\theta^\infty_2,\ \theta^\infty_3 \to -\theta^\infty_3 \
t_1 \to -\varepsilon^{-4} t_1-\varepsilon^{-6},\
t_2 \to -\varepsilon^{-4} t_2-\varepsilon^{-6},\ 
H_{t_1} \to -\varepsilon^4 H_{t_1},\
H_{t_2} \to -\varepsilon^4 H_{t_2}\\
&q_1 \to \varepsilon^{-3}(1+\varepsilon (q_1-\theta^\infty_2/p_1)),\
p_1 \to \varepsilon^2 p_1,\
q_2 \to \varepsilon^{-3}(1+\varepsilon (q_2-\theta^\infty_3/p_2)),\
p_2 \to \varepsilon^2 p_2.
\end{align*}
\noindent
{\bf 3+1$\to$ 4}

\noindent
$((1)(1))((1)),21 \to (((1)(1)))(((1)))$
\begin{align*}
&\theta^0 \to -\varepsilon^{-6},\ \theta^\infty_1 \to \theta^\infty_1+\varepsilon^{-6},\ 
t_1 \to \varepsilon t_1-2\varepsilon^{-3},\ t_2 \to \varepsilon t_2-2\varepsilon^{-3},\
H_{t_1} \to \varepsilon^{-1} H_{t_1},\ H_{t_2} \to \varepsilon^{-1} H_{t_2},\\
&q_1 \to \varepsilon^{-1}q_1+\varepsilon^{-3},\ p_1 \to \varepsilon p_1,\ 
q_2 \to \varepsilon^{-1}q_2+\varepsilon^{-3},\ p_2 \to \varepsilon p_2.
\end{align*}

\subsection{Sasano system}
\noindent
{\bf 1+1+1+1$\to$2+1+1}

\noindent
$31,22,22,1111\to(2)(2),31,1111$
\begin{align*}
&\theta^1 \to \theta^1- \varepsilon^{-1},\ \theta^t \to \varepsilon^{-1},\ 
t \to 1+\varepsilon t,\ 
H\to \varepsilon^{-1}(H+t^{-1}(p_1 q_1+p_2 q_2)),
\\
&q_1 \to 1+\varepsilon tq_1,\ p_1 \to \varepsilon^{-1}t^{-1}p_1,\ 
q_2 \to 1+\varepsilon tq_2,\ p_2 \to \varepsilon^{-1}t^{-1}p_2.
\end{align*}
\noindent
$31,22,22,1111\to(11)(11),31,22$
\begin{align*}
&\theta^t \to \varepsilon^{-1},\ \theta^\infty_1 \to \theta^\infty_3-\varepsilon^{-1},\ 
\theta^\infty_2 \to \theta^\infty_4-\varepsilon^{-1},\ \theta^\infty_3 \to \theta^\infty_1,\ 
\theta^\infty_4 \to \theta^\infty_2,\\
&t \to 1/\varepsilon t,\
H \to -\varepsilon t^2 H-\varepsilon t(p_1q_1+p_2q_2),\\
&q_1 \to 1/\varepsilon tq_2,
\ p_1 \to -\varepsilon tq_2(p_2q_2-\theta^1-\theta^\infty_2-\theta^\infty_4),\
q_2 \to 1/\varepsilon tq_1,\ p_2 \to -\varepsilon tq_1(p_1q_1-\theta^\infty_1).
\end{align*}
\noindent
$31,22,22,1111\to(111)(1),22,22$
\begin{align*}
&
\theta^{0} \to \varepsilon^{-1},\ \theta^{t}\to \theta^{0},\
\theta^{\infty}_{1} \to \theta^{\infty}_{1}-\varepsilon^{-1},\
\ t \to \frac{1}{1-\varepsilon t},\ H \to \varepsilon^{-1} H,\\
& q_{1} \to \frac{q_{1}-1}{q_{1}(1- \varepsilon t)},\
p_{1} \to q_{1}(1-\varepsilon t)(p_{1}q_{1}+\theta^{0}+\theta^{1}+\theta^{\infty}_{1}+\theta^{\infty}_{3}),\\
& q_{2} \to \frac{q_{2}-1}{q_{2}(1- \varepsilon t)},\
p_{2} \to q_{2}(1-\varepsilon t)(p_{2}q_{2}-\theta^{1}-\theta^{\infty}_{3}).
\end{align*}
\noindent
{\bf 2+1+1$\to$ 3+1}

\noindent
$(11)(11),31,22\to((11))((11)),31$
\begin{align*}
&\theta^1 \to \varepsilon^{-2},\ \theta^\infty_3 \to \theta^\infty_3-\varepsilon^{-2},\ 
\theta^\infty_4 \to \theta^\infty_4-\varepsilon^{-2},\
t \to -\varepsilon^{-1}t-\varepsilon^{-2},\ H \to -\varepsilon H,\\
&q_1 \to \frac{1}{1-\varepsilon q_2},
\ p_1 \to (1-\varepsilon q_2)(\varepsilon^{-1}p_2(1-\varepsilon q_2)+\theta^\infty_1),\\
&q_2 \to \frac{1}{1-\varepsilon q_1},
\ p_2 \to (1-\varepsilon q_1)(\varepsilon^{-1}p_1(1-\varepsilon q_1)+\theta^\infty_2+\theta^\infty_4).
\end{align*}
\noindent
$(2)(2),31,1111\to((11))((11)),31$
\begin{align*}
&\theta^{1} \to -\varepsilon^{-2},\ \theta^{\infty}_{1} \to \theta^{\infty}_{1}+\varepsilon^{-2},\ \theta^{\infty}_{2} \to \theta^{\infty}_{2}+\varepsilon^{-2},\ t \to -\varepsilon^{-1}t-\varepsilon^{-2},\ H \to -\varepsilon H,\\
& q_{1} \to \frac{1}{1-\varepsilon q_{1}},\ p_{1}\to (1-\varepsilon q_{1})(\varepsilon^{-1}p_{1}-p_{1} q_{1}-\theta^{0}-\theta^{\infty}_{1}-\theta^{\infty}_{3}),\\
& q_{2} \to \frac{1}{1-\varepsilon q_{2}},\ p_{2}\to (1-\varepsilon q_{2})(\varepsilon^{-1}p_{2}-p_{2} q_{2}+\theta^{\infty}_{3}).
\end{align*}
\noindent
{\bf 2+1+1$\to$ 2+2}

\noindent
$(2)(2),31,1111\to(2)(2),(111)(1)$
\begin{align*}
&\theta^0 \to \varepsilon^{-1},\ \theta^{1} \to \theta^{0},\ \theta^\infty_1 \to \theta^\infty_1- \varepsilon^{-1},\
t \to -\varepsilon t,\ H \to -\varepsilon^{-1}\left(H-\frac{p_{1}q_{1}+p_{2}q_{2}}{t}\right),\\
&q_1 \to \varepsilon^{-1} t^{-1}q_1,\ p_1 \to \varepsilon t p_1,\
q_2 \to \varepsilon^{-1} t^{-1}q_2,\ p_2 \to \varepsilon t p_2.
\end{align*}
\noindent
$(111)(1),22,22\to (2)(2),(111)(1)$
\begin{align*}
&\theta^{0} \to \theta^{0}-\varepsilon^{-1},\ \theta^{1} \to \varepsilon^{-1},\
t \to \varepsilon t,\ H \to \varepsilon^{-1} H,\\
&q_{1} \to -\frac{1}{\varepsilon q_{1}},\ p_{1} \to \varepsilon q_{1}(p_{1}q_{1}+\theta^{0}+\theta^{\infty}_{1}+\theta^{\infty}_{3}),\
q_{2} \to -\frac{1}{\varepsilon q_{2}},\ p_{2} \to q_{2}(\varepsilon (p_{2}q_{2}-\theta^{\infty}_{3})-1).
\end{align*}

\subsection{matrix Painlev\'e system}
\noindent
{\bf1+1+1+1$\to$2+1+1}

\noindent
$22,22,22,211\to(2)(2),22,211$
\begin{align*}
&\theta^1 \to \varepsilon^{-1},\ \theta^t \to \theta^1- \varepsilon^{-1},\
t \to 1+\varepsilon t,\ 
H \to \varepsilon^{-1}H,\\
&Q \to 1-\varepsilon P, \ P \to \varepsilon^{-1} Q.
\end{align*}
$22,22,22,211\to(2)(11),22,22$
\begin{align*}
&\theta^t \to \varepsilon^{-1},\ 
\theta^\infty_1 \to \theta^\infty_1-\varepsilon^{-1},\
t \to (\varepsilon t)^{-1},\ 
H \to -\varepsilon t^2H-\varepsilon t\ \mathrm{tr}(PQ),\\
&Q \to (\varepsilon t)^{-1}Q^{-1}, \ P \to -\varepsilon t(QP+\theta^0+\theta^\infty_1)Q.
\end{align*}
\noindent
{\bf 2+1+1$\to$ 3+1}

\noindent
$(2)(2),22,211\to((2))((2)),211$
\begin{align*}
& \theta^0 \to \theta^0-\varepsilon^{-2},\ 
\theta^1 \to \varepsilon^{-2},\ t \to \varepsilon^{-1}(-t+\varepsilon^{-1}),\ H \to -\varepsilon H,\\
&Q \to \varepsilon Q, \ P \to \varepsilon^{-1}P.
\end{align*}
\noindent
$(2)(2),22,211\to((2))((11)),22$
\begin{align*}
& \theta^1 \to \varepsilon^{-2},\ 
\theta^\infty_2 \to \theta^\infty_2-\varepsilon^{-2},\ \theta^\infty_3 \to \theta^\infty_3-\varepsilon^{-2},\
t \to \varepsilon^{-1}(-t-\varepsilon^{-1}),\
H \to -\varepsilon(H+\mathrm{tr}P),\\
&Q \to (1+\varepsilon P)^{-1}, \ P \to \{(P+\varepsilon^{-1})(-P+Q+t)
+\theta^0+2\theta^\infty_1+\varepsilon^{-2}-1\}(\varepsilon P+1).
\end{align*}
\noindent
$(2)(11),22,22\to((2))((11)),22$
\begin{align*}
& \theta^1 \to \varepsilon^{-2},\ 
\theta^\infty_1 \to \theta^\infty_1-\varepsilon^{-2},\
t \to \varepsilon^{-1}(-t-\varepsilon^{-1}),\
H \to -\varepsilon H,\\
&Q \to (1-\varepsilon Q)^{-1}, \ 
P \to \{(Q-\varepsilon^{-1})P+\theta^0+\theta^\infty_1-\varepsilon^{-2}\}(\varepsilon Q-1).
\end{align*}
\noindent
{\bf2+1+1$\to$2+2}

\noindent
$(2)(2),22,211\to(2)(2),(2)(11)$
\begin{align*}
&\theta^0 \to \varepsilon^{-1},\ \theta^1 \to \theta^0,\ \theta^\infty_1 \to \theta^\infty_1-\varepsilon^{-1},\
t \to \varepsilon t,\  
H \to \varepsilon^{-1}H,\ Q \to P, \ P \to -Q.
\end{align*}
\noindent
$(2)(11),22,22\to(2)(2),(2)(11)$
\begin{align*}
& \theta^0 \to -\varepsilon^{-1},\ \theta^1 \to \theta^0+\varepsilon^{-1},\
t \to \varepsilon t,\
H \to \varepsilon^{-1}H,\\
&Q \to (-\varepsilon Q)^{-1}, \ 
P \to \varepsilon(QP-\varepsilon^{-1}+\theta^\infty_1)Q.
\end{align*}
\noindent
{\bf 2+2$\to$ 4}

\noindent
$(2)(2),(2)(11)\to(((2)))(((11)))$
\begin{align*}
& \theta^0 \to -2\varepsilon^{-3},\ 
\theta^\infty_1 \to \theta^\infty_1+2\varepsilon^{-3},\
t \to -\varepsilon^{-4} t-\varepsilon^{-6},\
H \to -\varepsilon^4(H+\mathrm{tr}Q),\\
&Q \to \varepsilon^{-3}(1-\varepsilon Q), \ 
P \to \varepsilon^2(-P+Q^2+t).
\end{align*}
\noindent
{\bf 3+1$\to$ 4}

\noindent
$((2))((2)),211\to(((2)))(((11)))$
\begin{align*}
&\theta^0 \to -\varepsilon^{-6},\ 
\theta^\infty_2 \to \theta^\infty_2+\varepsilon^{-6},\ \theta^\infty_3 \to \theta^\infty_3+\varepsilon^{-6},\
t \to \varepsilon t-2\varepsilon^{-3},\ 
H \to \varepsilon^{-1}H,\\
&Q \to \varepsilon^{-1}Q+\varepsilon^{-3}, \ 
P \to \varepsilon P+\varepsilon^3(\theta^\infty_1-\varepsilon^{-6})(\varepsilon^2 Q+1)^{-1}.
\end{align*}
\noindent
$((2))((11)),22\to(((2)))(((11)))$
\begin{align*}
& \theta^0 \to -\varepsilon^{-6},\ 
\theta^\infty_1 \to \theta^\infty_1+\varepsilon^{-6},\
t \to \varepsilon t-2\varepsilon^{-3},\
H \to \varepsilon^{-1}(H+\mathrm{tr}Q),\\
&Q \to -\varepsilon^{-1}Q+\varepsilon^{-3}, \ P \to -\varepsilon(P-Q^2-t).
\end{align*}

\bibliographystyle{unsrt}

\end{document}